\documentclass[]{elsarticle}
\usepackage{lineno,hyperref}
\modulolinenumbers[5]
\makeatletter
\def\ps@pprintTitle{%
 \let\@oddhead\@empty
 \let\@evenhead\@empty
 \def\@oddfoot{}%
 \let\@evenfoot\@oddfoot}
\makeatother
\usepackage[a4paper, total={6.4in, 9.4in}]{geometry}
\usepackage[a4paper]{geometry}
 
\usepackage{color}
\usepackage[x11names]{xcolor}
\usepackage[utf8]{inputenc}
\usepackage[english]{babel}
\usepackage{amsmath}
\usepackage{pgfplots}
\usepackage{tabularx,ragged2e,booktabs,caption}

\usepackage[labelfont=bf]{caption}
\usepackage[font=footnotesize]{caption}
\usepackage{amssymb,latexsym}
\usepackage{graphicx}
\graphicspath{{./images/}}
\usepackage{booktabs}
\usepackage{amsfonts}
\usepackage{tikz}
\usepackage{array}
\usepackage{mathtools}
 \usepackage{relsize}
\usepackage[T1]{fontenc}
\usepackage{float}
\usepackage{cases}
\usepackage{multirow}
\usepackage{amsthm}
\usepackage{subcaption}
\pgfplotsset{compat=newest}
\usepackage[toc,page]{appendix}
\newcolumntype{Y}{>{\centering\arraybackslash}X}

\newtheorem{remark}{Remark}[section]
\numberwithin{equation}{section}
\usepackage{tabularx,ragged2e,booktabs,caption}

\usepackage[toc]{appendix}
\DeclareMathAlphabet{\mathpzc}{OT1}{pzc}{m}{it}

\begin{document}
\title{On local Fourier analysis of multigrid methods for PDEs with jumping and random coefficients}
\address[cwi]{CWI - Centrum Wiskunde \& Informatica, Amsterdam, The Netherlands}
\address[diam]{DIAM, Delft University of Technology, Delft, The Netherlands.}
\address[iuma]{IUMA and Applied Mathematics Department, University of Zaragoza}
\author[cwi]{Prashant~Kumar}
\ead{pkumar@cwi.nl}
\author[iuma]{Carmen~Rodrigo}
\ead{carmenr@unizar.es}
\author[cwi,iuma]{Francisco~J.~Gaspar}
\ead{gaspar@cwi.nl}
\author[cwi,diam]{Cornelis~W.~Oosterlee}
\ead{c.w.oosterlee@cwi.nl}
\begin{frontmatter}
% REQUIRED
\begin{abstract}
In this paper, we propose a novel non-standard Local Fourier Analysis (LFA) variant for accurately predicting the multigrid convergence of problems with random and jumping coefficients. This LFA method is based on a specific basis of the Fourier space rather than the commonly used Fourier modes. To show the utility of this analysis, we consider, as an example, a simple cell-centered multigrid method for solving a steady-state single phase flow problem in a random porous medium. We successfully demonstrate the prediction capability of the proposed LFA using a number of challenging benchmark problems. The information provided by this analysis helps us to estimate a-priori the time needed for solving certain uncertainty quantification problems by means of a multigrid multilevel Monte Carlo method.

%Additionally, a simple cell-centered multigrid solver for lognormal diffusion problem is presented that is highly efficient and robust with respect to the heterogeneity of the random field. Finally, this multigrid solver is ultilized to perform multilevel Monte Carlo simulations with different permeability parameters.
\end{abstract}

% REQUIRED
\end{frontmatter}
% REQUIRED
%\begin{AMS}
%  65F10, 65M22, 65M55
%\end{AMS}

\section{Introduction}\label{sec:1}

\setcounter{section}{1}
A number of problems in science and engineering involve solving partial differential equations (PDEs) with random parameters or coefficients.  The solution to such PDEs is of stochastic nature, and the aim is to compute expected values and corresponding variances of a functional of the solution. For such uncertainty quantification (UQ) problems, Monte Carlo (MC) type methods are preferred due to their dimension-independent convergence. For any sampling-based approach, the availability of a highly efficient and robust (w.r.t. the random inputs) iterative solver becomes critical. In general, the sample-wise computational cost can become highly heterogeneous, depending on the random inputs. Therefore, if the performance statistics of such solvers were known a-priori, one could utilize this information to optimize and parallelize the MC simulations efficiently. 

In this paper, we present a non-standard Local Fourier Analysis (LFA) technique to predict the convergence rate of multigrid solvers for problems involving random and jumping coefficients.  Standard LFA techniques are typically based on constant coefficient discretization stencils, whereas for stochastic PDEs we encounter varying coefficients throughout the computational domain, due to the randomness. One of the main contributions of this work is to generalize the LFA towards problems with random and jumping coefficients, with the aim of predicting, a-priori, the total time needed to solve UQ problems. Some efforts have already been done in \cite{bolten} regarding the generalization of LFA for jumping coefficients. The novelty of our approach lies in the choice of basis functions. Here, we utilize a new basis from the Fourier space rather than the standard Fourier modes. We benchmark the prediction capability of the proposed LFA technique using a set of challenging jumping coefficient problems and a number of spatially correlated random fields with varying heterogeneity.

One of the results of this paper is the ability of the aforementioned LFA technique to predict accurately the multigrid convergence rates for elliptic PDEs with heterogeneous random coefficients, as encountered in the stochastic Darcy flow problem. In stochastic subsurface flow modeling, a lognormal random field often represents the permeability of an occurring heterogeneous porous medium \cite{prob3,WRCR:WRCR1821,WRCR:WRCR3746}. Uncertainty may then be assessed by means of stochastic collocation approaches \cite{babuska-collocation-2007,Xiu:2010}, where the use of multigrid has been analyzed in~\cite{SEYNAEVE2007132}. In the current paper, we are interested in solving the PDE in the context of the Multilevel Monte Carlo (MLMC) technique~\cite{MLMC2,MLMC1,ANU:9672986}. Within Monte Carlo methods, samples of the stochastic random field are generated and the corresponding solution of the PDE with ``frozen stochastic coefficients'' is computed. Plain Monte Carlo methods may require many thousands of random samples on a fine computational mesh, before the mean and variance of the numerical solution stabilize and converge. This may cost substantial CPU time. Multilevel Monte Carlo methods are based on the decomposition of the expectation and variance operators into sequences of differences of these quantities on fine and coarse problem scales. The number of iterations required on each scale is different, where fewer samples (fewer random fields) are typically required on the finer scales, and larger numbers of samples are needed on the coarser, cheaper to compute, scales.

Combining MLMC methods with multigrid seems natural, where the multigrid method is employed for the numerical solution of a PDE which is based on a sample of the stochastic quantity on a certain (fine or coarse) scale. Due to the stochasticity, however, we deal with PDE problems with jumping coefficients, where different jump patterns are encountered each time a new random field is generated. The generalization of the LFA towards these PDEs will provide us insight in the average number of multigrid iterations, the spread of the convergence factors, amongst other things. This information helps to estimate the total CPU time needed for the multiple multigrid computations in a multilevel Monte Carlo setting. We consider it useful to construct a technique to assess the quality of the choice of the multigrid components in the context of the PDEs with random problems, before the actual multigrid computation has taken place.

In this work, we will employ a basic cell-centered multigrid (CCMG) algorithm for solving elliptic PDEs with a variable coefficient field. The components of this algorithm include a simple Gauss-Seidel iteration as the smoother, a piecewise constant prolongation operator and its adjoint as the restriction and a direct discretization technique to define the discrete operators on the coarse grids. We show that for this special combination, the coarse grid discretization operators are equivalent to the ones obtained from commonly used Galerkin operators \cite{brandt2011multigrid}. We utilize this CCMG method to perform MLMC simulations with different permeability parameters. It may be surprising that such a basic algorithm converges well in the context of the generated random fields, where computation takes place for thousands of different samples. Despite we restrict ourselves to this basic CCMG to demonstrate the accuracy of the predictions of the novel LFA technique, we emphasize that this approach can be used for a wider range of problems, discretizations and multigrid methods. In fact, we think that the proposed LFA allows to deal in a easy way with challenging problems for which a standard LFA is very difficult to apply or even impracticable.

The paper is organized as follows. In Section \ref{sec:2} we introduce the context of PDEs with jumping and random coefficients, together with their discretization by a cell-centered finite volume scheme. 
A discussion on multigrid methods for this type of problems is also included, and the multigrid components that will be considered in this work are also defined.
Section \ref{sec:3} is devoted to the generalization of the LFA to deal with jumping coefficients and problems with random fields. In Section \ref{sec:4}, we present results obtained by this analysis for different benchmark problems with jumping coefficients. Section \ref{sec:5} presents the LFA results for problems with random coefficients, and in Section \ref{sec:6} multilevel Monte Carlo computations for PDEs with random coefficients are presented. Finally, in Section \ref{sec:7} conclusions are drawn.

\section{Jumping coefficients, random coefficients, multigrid}\label{sec:2}

\setcounter{section}{2}
Robust and efficient iterative solution methods are very relevant for partial differential equations with variable coefficients. For PDEs with jumping coefficients, multigrid methods have already shown to be this type of solvers. When using Monte Carlo methods, in the case of elliptic PDEs with random coefficient fields, many samples of the random field are generated and for each field the numerical solution should be computed. This can take substantial CPU time if very many samples are required. For a fixed sample of the random field, we deal with an elliptic PDE with varying coefficients, due to the randomness. Multigrid comes in naturally as a highly efficient solution method for the resulting PDEs. In this section, we explain this setting and we briefly describe an efficient multigrid method based on a cell-centered grid and a finite volume discretization.

\subsection{PDEs with jumping and with random coefficients}
We start with classical PDE problems with jumping coefficients.  In particular, we deal with the following two-dimensional diffusion equation on the square domain ${ D} = (0,\ell)^2$ 
\begin{eqnarray}\label{diffusion_problem}
- \nabla \cdot \left(k({\mathbf x}) \nabla u ({\mathbf x})\right) & = & f({\mathbf x}), \quad {\mathbf x} \in {D}, \\
u({\mathbf x}) & = & g({\mathbf x}), \quad {\mathbf x} \in \partial {D}, 
\end{eqnarray}
where $k({\mathbf x})$ is a function which may be discontinuous across internal boundaries. 

To discretize this problem, we use a {\em cell-centered finite volume method} based on the harmonic average of the diffusion coefficient $k({\mathbf x})$. We consider a uniform grid $D_h$ with the same step size $h = \ell/M, M \in \mathbb{N}$ in both directions, 
\begin{equation}
\label{Grid_fine}
D_h = \{(x_{i_1},x_{i_2}); x_{i_\alpha} = (i_{\alpha}-1/2) h, i_{\alpha} = 1,\ldots,M, \alpha=1,2\}. 
\end{equation}
This gives, for each interior cell with center $(x_{i_1},x_{i_2})$, denoted by $D_h^{i_1,i_2}$, a five-point scheme
\begin{equation}
c^h_{i_1,i_2} u_{i_1,i_2} + w^h_{i_1,i_2} u_{i_1-1,i_2} + e^h_{i_1,i_2} u_{i_1+1,i_2} + s^h_{i_1,i_2} u_{i_1,i_2-1} + n^h_{i_1,i_2} u_{i_1,i_2+1} = f^h_{i_1,i_2},
\label{five_stencil}
\end{equation}
where 
\begin{eqnarray}
w^h_{i_1,i_2} &=& -\frac{2}{h^2} \frac{k_{i_1,i_2} k_{i_1-1,i_2}}{k_{i_1,i_2} + k_{i_1-1,i_2}}, \quad
e^h_{i_1,i_2} = -\frac{2}{h^2} \frac{k_{i_1,i_2} k_{i_1+1,i_2}}{k_{i_1,i_2} + k_{i_1+1,i_2}}, 	\nonumber \\
s^h_{i_1,i_2} &=& -\frac{2}{h^2} \frac{k_{i_1,i_2} k_{i_1,i_2-1}}{k_{i_1,i_2} + k_{i_1,i_2-1}}, \quad
n^h_{i_1,i_2} = -\frac{2}{h^2} \frac{k_{i_1,i_2} k_{i_1,i_2+1}}{k_{i_1,i_2} + k_{i_1,i_2+1}}, \nonumber \\
c^h_{i_1,i_2} &=& -(w^h_{i_1,i_2} + e^h_{i_1,i_2} + n^h_{i_1,i_2}+ s^h_{i_1,i_2}) \nonumber,
\end{eqnarray}
with, for instance, $k_{i_1,i_2}$ the diffusion coefficient associated with the cell $D_h^{i_1,i_2}$. By
{\em interior cell} we mean a cell for which none of its edges lies at the boundary of the domain. This scheme is changed appropriately for the cells close to the boundary. \\

As mentioned, we also consider elliptic PDEs with random coefficient fields. The PDE of our interest describes the steady-state single-phase flow in a random porous medium. Denoting by $\omega$  an event in the probability space $(\Omega,\mathbb{F},\mathbb{P})$, with sample space $\Omega$, $\sigma$-algebra $\mathbb{F}$ and probability measure $\mathbb{P}$, the permeability in the porous medium is described by $k({\mathbf x}, \omega): \overline{{D}} \times \Omega   \rightarrow  \mathbb{R}_+$.  The PDE is then given by 
\begin{equation}
-\nabla \cdot (k({\mathbf x}, \omega)\nabla u({\mathbf x}, \omega)) = f({\mathbf x}), \quad  {\mathbf x} \in { D}, \omega \in \Omega, 
\label{1a}
\end{equation}
with $f$ as a source term. 

The engineering interest in the solution of (\ref{1a}) is typically found in expected values of linear functionals of the solution $u$, denoted by $Q := Q(u)$. 

To discretize these problems, we use the same cell-centered finite volume method based on the harmonic average of the random diffusion coefficient as previously described for problems with jumping coefficients. We make the common assumption that the permeability random field is constant over each cell of the grid. 

\subsection{Multigrid for PDEs with jumping and with random coefficients}\label{MGalgo}

\begin{figure}[htb]
\begin{center}
\includegraphics[width =0.9\textwidth]{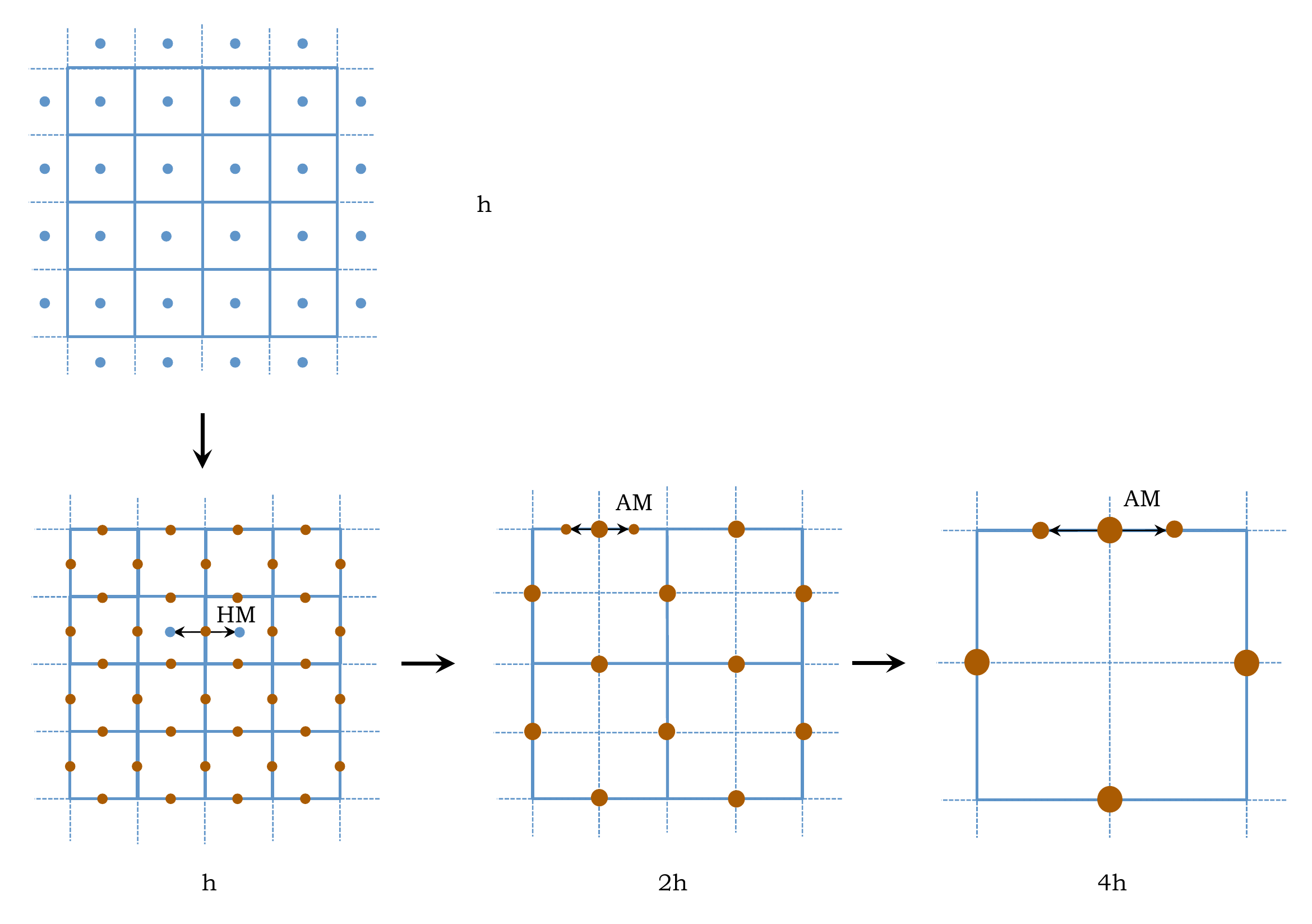} 
\end{center}
\caption{Schematic representation of permeability upscaling used in the multigrid hierarchy (h-2h-4h). (Top) Permeability values generated at cell-centres (blue dots). (Bottom left) Permeability values at face-centres (red dots) obtained from the harmonic mean (HM) of permeabilities from two adjacent cell-centres. (Bottom middle) Permeability at face-centres (bigger red dots) of 2h-grid is the arithmetic mean (AM) of permeabilities from face-centres of the h-grid. (Bottom right) Permeability at face-centres (biggest red dots) of 4h-grid is the arithmetic mean (AM) of permeabilities from face-centres of the 2h-grid.}
\label{Perm_coarse}
\end{figure}

In this work, the multigrid components for the above cell-centered discrete problems are chosen as follows. We use a lexicographic Gauss-Seidel iteration as the multigrid smoother and we consider standard coarsening obtained by doubling the mesh size in both directions. Each coarse cell is the union of four fine cells, and, since the unknowns are located in the cell-centres, this results in a non-nested hierarchy of grids. We consider a simple prolongation operator $P_{2h}^h$, that is, the piecewise constant interpolation operator. In stencil notation, it is given by
\begin{equation}
\label{cp}
P_{2h}^h = 
\left]
\begin{array}{ccc}
1 &  & 1 \\ 
& \star & \\ 
1 & & 1
\end{array}
\right[_{2h}^{h},
\end{equation}
where $\star$ denotes the position of a coarse grid unknown. The classical stencil notation shows the contribution of the coarse grid node to the neighbouring fine grid nodes. The restriction operator ${R}_h^{2h}$ is chosen as the scaled adjoint of the prolongation, given in stencil form by
\begin{equation}
\label{cr}
R_{h}^{2h} = \frac{1}{4}
\left[
\begin{array}{ccc}
1 & & 1 \\ & \star & \\ 1 & & 1
\end{array}
\right]_{h}^{2h}.
\end{equation}
 The coarse grid operators are constructed by direct discretization defining the diffusion coefficients at the edges of the coarse cells appropriately, which we will describe in more detail. We assume that the diffusion coefficient $k({\mathbf x})$ is piecewise constant on the fine grid.  The flux over an edge, dependent on the solution in the two adjacent cells, is calculated based on the harmonic average. 
The values of the diffusion coefficients at a coarse edge located between two coarse cells, however, are calculated as the arithmetic average of the corresponding fine grid coefficients, see Figure \ref{Perm_coarse} for a more detailed description. 
As it was pointed in \cite{MOLENAAR199625}, this direct discretization procedure is equivalent to the often used Galerkin approach, i.e.,  $L_{2h} = \frac{1}{2}{R}_h^{2h} L_h P_{2h}^h$, but  computationally more efficient. The factor $1/2$ in the previous expression is due to the lack of consistency of the operator ${R}_h^{2h} L_h P_{2h}^h$ with the differential operator \cite{Yav98}. In the next result, we prove that both discretizations are indeed equivalent. \\

%\begin{proposition}
\vspace{0.2cm}
PROPOSITION 1. \emph{Let $L_h$ be the fine-grid operator based on the cell-centered finite volume discretization of problem \eqref{1a} on a uniform grid of mesh size $h=\ell/M$ with $M$ even. Let $P_{2h}^h$ be the piecewise constant prolongation operator and ${R}_h^{2h}$ its adjoint. Then, the Galerkin coarse grid operator $L_{2h} = \frac{1}{2}{R}_h^{2h} L_h P_{2h}^h$ is equivalent to a direct discretization on the coarse grid based on the arithmetic average of the corresponding fine grid coefficients.}
%\end{proposition}
\begin{proof}
We prove the equivalence for a coarse grid cell $D_{2h}^{i_1,i_2}$ such that none of its edges lies on the boundary of the domain. The equivalence for coarse cells close to boundaries with Dirichlet or Neumann boundary conditions can be proven similarly. By applying the restriction operator ${R}_h^{2h}$ in \eqref{cr}, the equation associated with the cell $D_{2h}^{i_1,i_2}$ by using the Galerkin approach is given by
\begin{eqnarray}
\label{formu1}
\frac{1}{2}({R}_h^{2h} L_h P_{2h}^h u)_{i_1,i_2} & = & \frac{1}{8} ((L_h P_{2h}^h u)_{2i_1,2i_2}  + (L_h P_{2h}^h u)_{2i_1-1,2i_2} + (L_h P_{2h}^h u)_{2i_1,2i_2-1} + \nonumber \\
& & (L_h P_{2h}^h u)_{2i_1-1,2i_2-1}) .
\end{eqnarray}
Taking into account that the prolongation operator is piecewise constant, we obtain 
\begin{eqnarray}
(L_h P_{2h}^h u)_{2i_1,2i_2} &=& e^h_{2i_1,2i_2} u_{i_1+1,i_2} + w^h_{2i_1,2i_2}  u_{i_1,i_2} + n^h_{2i_1,2i_2} u_{i_1,i_2+1} +  \nonumber \\
\quad & & s^h_{2i_1,2i_2}  u_{i_1,i_2} + c^h_{2i_1,2i_2} u_{i_1,i_2}, \nonumber
\end{eqnarray}
and similar expressions for the other terms in \eqref{formu1}. By substituting these expressions in \eqref{formu1}, the following discretization for the coarse cell $D_{2h}^{i_1,i_2}$ is obtained
\begin{align}
\label{formu3}
\frac{1}{2}({R}_h^{2h} L_h P_{2h}^h u)_{i_1,i_2} = &c^{2h}_{i_1,i_2} u_{i_1,i_2} + w^{2h}_{i_1,i_2} u_{i_1-1,i_2} + e^{2h}_{i_1,i_2} u_{i_1+1,i_2} + \nonumber\\
&s^{2h}_{i_1,i_2} u_{i_1,i_2-1} + n^{2h}_{i_1,i_2} u_{i_1,i_2+1},
\end{align}
where
\begin{eqnarray}
w^{2h}_{i_1,i_2} &=& \displaystyle\frac{1}{8}\left(w^h_{2i_1-1,2i_2} + w^h_{2i_1-1,2i_2-1}\right), \quad
e^{2h}_{i_1,i_2} = \displaystyle\frac{1}{8}\left(e^h_{2i_1,2i_2} + e^h_{2i_1,2i_2-1}\right),	\nonumber \\
s^{2h}_{i_1,i_2} &=& \displaystyle\frac{1}{8}\left(s^h_{2i_1,2i_2-1} + s^h_{2i_1-1,2i_2-1}\right), \quad
n^{2h}_{i_1,i_2} = \displaystyle\frac{1}{8}\left(n^h_{2i_1-1,2i_2} + n^h_{2i_1,2i_2}\right),	\nonumber \\
c^{2h}_{i_1,i_2} &=& -(w^{2h}_{i_1,i_2} + e^{2h}_{i_1,i_2} + n^{2h}_{i_1,i_2}+ s^{2h}_{i_1,i_2}).\nonumber
\end{eqnarray}
\end{proof}
We observe that this scheme is equivalent to a direct discretization on the coarse grid where the diffusion coefficients on the edges are the arithmetic averages of the corresponding fine grid coefficients.  
\begin{remark} In order to achieve a mesh-independent multigrid convergence following the analysis from \cite{hackbusch85}, the next condition must be satisfied:
\begin{equation}\label{orders}
m_p +m_r > M_{pde},
\end{equation}
where the orders $m_p$ and $m_r$ are the highest degree plus one of the polynomials that are exactly interpolated by $P_{2h}^h$ and $R^{2h}_h$, respectively, and $M_{pde}$ is the order of the PDE to be solved. For PDE \eqref{1a}, we have $M_{pde}=2$ and for the considered operators \eqref{cp} and \eqref{cr}, we get $m_p=m_r=1$, which is does not satisfy the inequality \eqref{orders}. 
In \cite{Bramble1996}  is shown however, that this condition is not needed to prove uniform convergence. The idea is to use the power of the theoretical approach provided in \cite{Bramble1991}.
\end{remark}

\subsection{Discussion about other multigrid methods for jumping coefficients}
In the context of {\em algebraic multigrid} methods for the numerical solution of partial differential equations, basically two prevailing methods have proved their use for multiple engineering problems, i.e., algebraic multigrid and aggregation-based multigrid methods,  \cite{Brezina,stuben2001introduction,Van2001,Braess1995,Van1996}. These methods converge remarkably well, for example, for scalar PDEs with jumping coefficients. It is not always easily understood why these methods, and particularly the aggregation-based method, converge so well.

The origin of these algebraic methods may be found already in the early days of multigrid, where black-box multigrid with operator-dependent transfer operators (restriction and prolongation) and Galerkin coarse grid operators for structured vertex-centered Cartesian grids was proposed in \cite{Alcouffe:1981:MGM,dendy0,DENDY1983261}. This can be seen as a predecessor of classical AMG, where these components were essentially enhanced by a flexible coarsening strategy.

The aggregation-based multigrid methods, with their origin in the work by Mandel \cite{Van1996} (smoothed aggregation), may be related to the cell-centered multigrid methods as proposed in 
\cite{MG2,WESSELING198885}. In~\cite{WESSELING198885}, it was shown that constant, i.e., operator-{\it independent}, transfer operators, in combination with Galerkin coarse grid discretization provided highly efficient multigrid results for cell-centered discretizations of elliptic PDEs that included jumping coefficients. These cell-centered multigrid components were augmented with robust smoothing, like Incomplete Lower-Upper decomposition (ILU) relaxation. The individual contributions of the coarse grid correction and the smoothing parts were difficult to distinguish. Also, a cell-centered multigrid based on coarsening by a factor of three together with operator-dependent 
interpolations was also explored in \cite{dendy2}.

%So, in the early days of multigrid, among others, vertex-centered multigrid with operator-induced transfer operators as well as cell-centered multigrid with constant coefficient transfer operators existed, both working very well based on evidence from numerous experiments, as stand-alone iterative solvers (without any Krylov subspace acceleration) for elliptic PDEs with jumping coefficients.
%We will analyze below the cell-centered case in more detail, by means of the generalization of the LFA, as described in Section \ref{sec:2}.

\section{Local Fourier Analysis for variable coefficients}\label{sec:3}
\setcounter{section}{3}
In this section we describe LFA in a setting which allows us to estimate the multigrid convergence factors for problems with jumping coefficients and problems with random fields. A discrete linear operator with constant coefficients, which is formally defined on an infinite grid,  is usually assumed for carrying out a standard local Fourier analysis.  As we will show, this assumption can be relaxed by considering a discrete operator with constant coefficients {\em in appropriate infinite subgrids}. This allows us to generalize the analysis to problems for which the discrete operator consists of different stencils. A key point in this improved analysis is to consider a specific basis of the Fourier space, rather than the standard basis which is based on the Fourier modes. The use of this new basis will simplify  the analysis.

\begin{figure}[h]
\centering
\includegraphics[width = 0.5\textwidth]{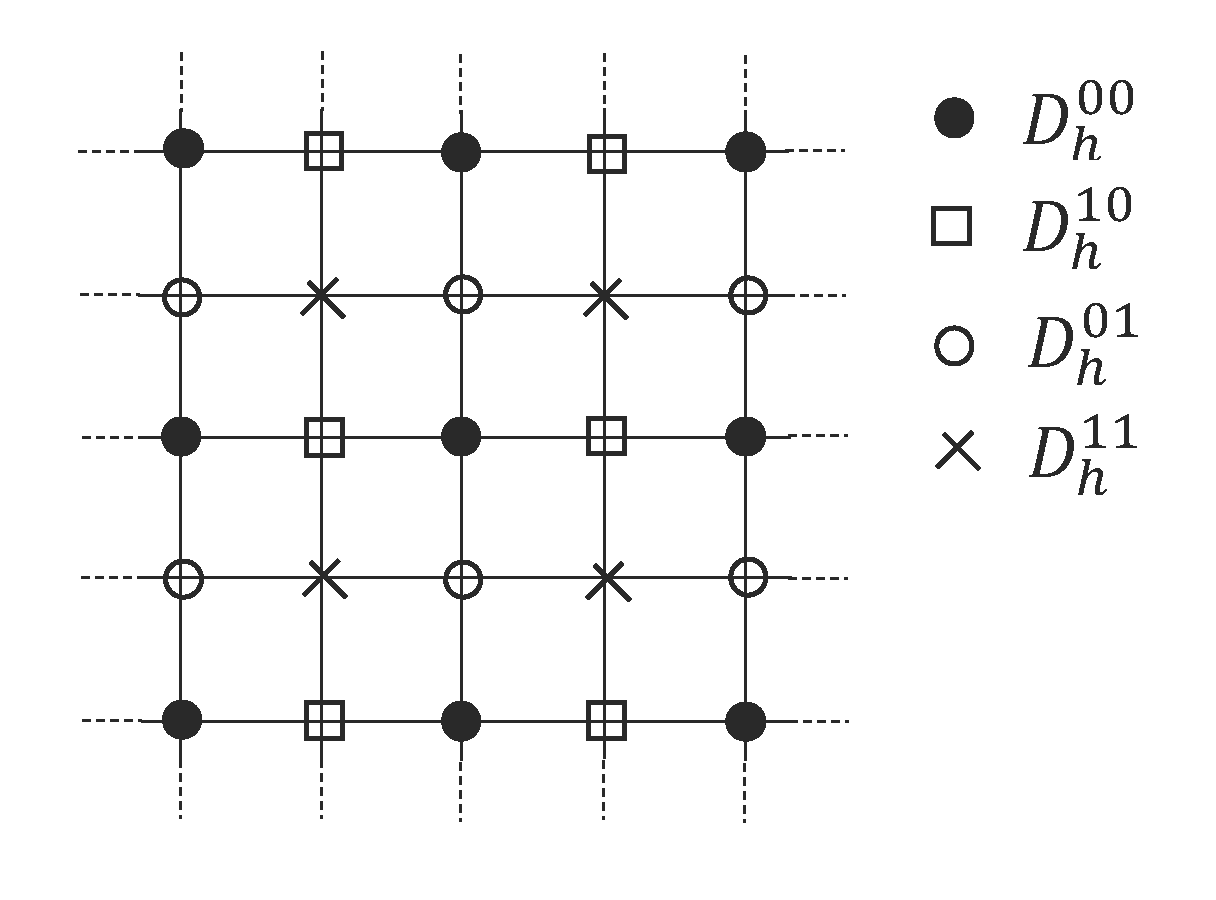}
\caption{Infinite grid $D_h$ divided into the corresponding subgrids for $n=2$.}\label{LFA_grid}
\end{figure}

We start from a regular infinite grid $D_h$ with grid size $h$ in both directions. Such an infinite grid will be {\em split into $n\times n$ subgrids} in the following way. First of all, a window comprising $n\times n$ cells of the original grid is {\it adequately} chosen, and, subsequently, we consider its periodic extension. The choice of the size of the $n\times n$ window is made such that the variability of the discrete operator in the computational grid can be appropriately represented, as will be explained by means of examples of different nature. Once $n$ is fixed, the infinite subgrids are defined as follows (see Figure \ref{LFA_grid} for an example with $n=2$),
\begin{equation}\label{subgrids}
D_h^{kl} = \left\{(k,l)h+(nk_1,nk_2)h \; | \; k_1,k_2\in {\mathbb Z}\right\},\; k,l=0,\ldots,n-1.
\end{equation}

For each low frequency, ${\boldsymbol \theta}^{00}\in{\boldsymbol \Theta}_{nh} = (-\pi/nh, \pi/nh]^2$, we introduce the grid-functions:
\begin{equation}\label{new_modes}
\psi_h^{kl}({\boldsymbol \theta^{00}},{\mathbf x}) = \varphi_h({\boldsymbol \theta}^{00},{\mathbf x})\chi_{D_h^{kl}}({\mathbf x}), \; k,l=0,\ldots,n-1,\; {\mathbf x}\in D_h,
\end{equation} 
where $\displaystyle  \varphi_h({\boldsymbol \theta}^{00},{\mathbf x}) = e^{\imath {\boldsymbol \theta}^{00}\cdot {\mathbf x}}$ is the standard Fourier mode on $D_h$
corresponding to the frequency ${\boldsymbol \theta}^{00}$. It is easy to see that the subspace generated by these $n^2$ grid-functions,
\begin{equation}\label{nndimsub}
{\mathcal F}_h^{n^2}({\boldsymbol \theta}^{00}) = span\{\psi_h^{kl}({\boldsymbol \theta}^{00},\cdot),\; k,l=0,\ldots,n-1 \}
\end{equation}
is the same as the one spanned by the $n^2$ Fourier modes $\displaystyle  \varphi_h({\boldsymbol \theta}^{00}_{kl},\cdot)$ associated with the frequencies:
\begin{equation}\label{frecuencies}
{\boldsymbol \theta}^{00}_{kl} = {\boldsymbol \theta}^{00} + (k,l)\frac{2\pi}{nh},\; k,l=0,\ldots,n-1.
\end{equation}
%\begin{example}
In the case $n=2$, the basis $\{\psi_h^{00}({\boldsymbol \theta}^{00},\cdot), \psi_h^{11}({\boldsymbol \theta}^{00},\cdot), \psi_h^{10}({\boldsymbol \theta}^{00},\cdot),\psi_h^{01}({\boldsymbol \theta}^{00},\cdot)\}$ is related to the standard basis of Fourier modes $\{\varphi_h({\boldsymbol \theta}^{00},\cdot), \varphi_h({\boldsymbol \theta}^{11},\cdot), \varphi_h({\boldsymbol \theta}^{10},\cdot),\varphi_h({\boldsymbol \theta}^{01},\cdot)\}$ in the following way:
\begin{equation}\label{rel_psi_phi}
\left(\begin{array}{c} 
\psi_h^{00}({\boldsymbol \theta}^{00},\cdot)\\ 
\psi_h^{11}({\boldsymbol \theta}^{00},\cdot)\\ 
\psi_h^{10}({\boldsymbol \theta}^{00},\cdot)\\ 
\psi_h^{01}({\boldsymbol \theta}^{00},\cdot)
\end{array}\right) = {\cal M} \left(\begin{array}{c} 
\varphi_h({\boldsymbol \theta}^{00},\cdot)\\ 
\varphi_h({\boldsymbol \theta}^{11},\cdot)\\ 
\varphi_h({\boldsymbol \theta}^{10},\cdot)\\ 
\varphi_h({\boldsymbol \theta}^{01},\cdot)
\end{array}\right) = \frac{1}{4}
\left(\begin{array}{cccc}
1 & 1 & 1 & 1 \\ 1 & 1 & -1 & -1 \\ 1 & -1 & -1 & 1 \\ 1 & -1 & 1 & -1
\end{array}\right)
\left(\begin{array}{c} 
\varphi_h({\boldsymbol \theta}^{00},\cdot)\\ 
\varphi_h({\boldsymbol \theta}^{11},\cdot)\\ 
\varphi_h({\boldsymbol \theta}^{10},\cdot)\\ 
\varphi_h({\boldsymbol \theta}^{01},\cdot)
\end{array}\right).
\end{equation}
%\end{example}
It is well-known that Fourier modes are eigenfunctions of any constant coefficient linear discrete operator $L_h$, that is,
$L_h\varphi_h({\boldsymbol \theta},{\mathbf x}) = \widetilde{L}_h({\boldsymbol \theta}) \varphi_h({\boldsymbol \theta},{\mathbf x})$.
Therefore, the representation of $L_h$ with respect to a basis of $n^2$ Fourier modes is a diagonal matrix with diagonal elements $\widetilde{L}_h({\boldsymbol \theta}^{kl})$ with $k,l=0,\ldots,n-1$. In general, the Fourier representation with respect to the basis of functions $\{\psi_h^{kl}\}_{k,l=0}^{n-1}$ is a dense matrix. We will denote it by $\widehat{L}_h({\boldsymbol \theta}^{00})$. 
%\begin{example}

If we consider the five-point standard discretization of the Laplace operator on a uniform grid of mesh size $h$, 
\begin{equation}\label{Laplace}
\frac{1}{h^2}\left[\begin{array}{ccc} & -1 & \\ -1 & 4 & -1 \\ & -1 &  \end{array}\right],
\end{equation}
its Fourier symbol with respect to the standard basis of Fourier modes is a diagonal matrix with diagonal elements equal to
 $$\frac{1}{h^2} (4 - 2 \cos(\theta_x^{kl}) - 2 \cos(\theta_y^{kl})),$$
% $$\displaystyle \frac{4}{h^2}\left(\sin^2(\theta_x^{kl}/2)+\sin^2(\theta_y^{kl}/2)\right),$$ 
(see \cite{MG1}, for instance) whereas the Fourier representation with respect to the new basis in the case $n=2$ is given by
\begin{equation}\label{symbol_new}
\widehat{L}_h({\boldsymbol \theta^{00}}) = 
\frac{2}{h^2}
\left( 
\begin{array}{cccc} 
2 & 0 & \cos(\theta_x^{00}) & \cos(\theta_y^{00})\\
0 & 2 & \cos(\theta_y^{00}) & \cos(\theta_x^{00}) \\
\cos(\theta_x^{00}) & \cos(\theta_y^{00}) & 2 & 0 \\
\cos(\theta_y^{00}) & \cos(\theta_x^{00}) & 0 & 2
\end{array}\right), \; \hbox{with} \; {\boldsymbol \theta^{00}} = (\theta_x^{00},\theta_y^{00}).
\end{equation}
Notice that, for example, the first row of the previous symbol is obtained by looking at the decomposition of the stencil \eqref{Laplace} into the connections among the unknowns located at the different subgrids $D_h^{kl}$ defined in \eqref{subgrids}. In particular, following the notations in Figure \ref{LFA_grid}, the $\bullet$-$\bullet$, $\bullet-\square$ and $\bullet-\circ$ connections are given by the following stencils $$\frac{1}{h^2}\left[4\right], \quad \frac{1}{h^2}\left[\begin{array}{ccc}-1 & \bullet & -1 \end{array}\right], \quad \frac{1}{h^2}\left[\begin{array}{c} -1 \\ \bullet \\ -1 \end{array}\right],$$ giving rise to the symbols $$\frac{4}{h^2}, \quad \frac{2}{h^2}\cos(\theta_x^{00}),\quad \frac{2}{h^2}\cos(\theta_y^{00}),$$ which appear in the first row of \eqref{symbol_new}, whereas there is no $\bullet-\times$ connection. The rest of the rows are analogously computed. \\
%\end{example}
The procedure to obtain the Fourier symbol of a smoothing operator $S_h$, which is based on a splitting of the discrete operator $L_h= L_h^+ + L_h^-$, is analogous with the new basis. The smoothing iteration is given by 
$$
L_h^+ \overline{w}_h + L^-_h w_h = f_h,
$$
with $w_h$ the approximation of the solution before the smoothing step and $\overline{w_h}$ the approximation after the smoothing step.
 By computing the symbols of $L_h^+$ and $L_h^-$ as before, the Fourier symbol of the smoothing operator is given by
\begin{equation} \label{sym_smoother}
\widehat{S}_h({\boldsymbol \theta^{00}}) = - (\widehat{L}_h^+)^{-1}({\boldsymbol \theta^{00}})\widehat{L}_h^-({\boldsymbol \theta^{00}}).
\end{equation}
%\begin{example}
The Fourier symbol corresponding to a lexicographic Gauss-Seidel iteration for the five-point standard discretization of the Laplace operator on a uniform grid of mesh size $h$, in the case $n=2$, is as in \eqref{sym_smoother}, where 
\begin{equation}\label{symbol_smth}
\widehat{L}_h^+({\boldsymbol \theta^{00}}) = 
\frac{2}{h^2}
\left( 
\begin{array}{cccc} 
2 & 0 & 0 & 0\\
0 & 2 & \cos(\theta_y^{00}) & \cos(\theta_x^{00}) \\
\cos(\theta_x^{00}) & 0 & 2 & 0 \\
\cos(\theta_y^{00}) & 0 & 0 & 2
\end{array}\right), \; \widehat{L}_h^-({\boldsymbol \theta^{00}}) = \widehat{L}_h({\boldsymbol \theta^{00}})-\widehat{L}_h^+({\boldsymbol \theta^{00}}).
\end{equation}
%\end{example}
Once the Fourier representation of the smoothing operator with respect to the new 
basis is obtained, we can define the smoothing factor by using the change of basis matrix. For example, for the case $n=2$, the smoothing factor is obtained by
$$
\mu(S_h) = \sup_{{\boldsymbol \theta}^{00}\in{\boldsymbol \Theta}_{2h}}  \rho( {\cal Q}_h 
{\cal M} \widehat{S}_h({\boldsymbol \theta}^{00}) {\cal M}^{-1}),
$$
where ${\cal Q}_h$ is the projection operator onto the space of high-frequency components and ${\cal M}$ is the change of basis matrix given in \eqref{rel_psi_phi}. \\

\subsection{LFA formulations for cell-centered grids}\label{sec:cellcentered}
In Section \ref{sec:2}, we didn't distinguish between cell- and vertex-centered grids. The generalized LFA indeed works well for both types of discretization. By introducing the coarse grids and their relation with the fine grids, we need to fix the approach of interest. Since here we will focus on cell-centered discretizations, from now on the description of the analysis will be given for this case, although it may be applied to the vertex-centered case in a similar way by defining appropriately the coarse meshes. \\
According to the location of the coarse-grid points in a regular cell-centered grid, we define for a fixed $n$ the following infinite coarse subgrids of $D_{2h}$:
\begin{equation}\label{subgrids_coarse}
D_{2h}^{kl} = \left\{(h/2,h/2) + (k,l)2h+(nk_1,nk_2)h \; | \; k_1,k_2\in {\mathbb Z}\right\},\; k,l=0,\ldots,n/2-1.
\end{equation}
Due to the relation between the grid-functions of the new Fourier basis given in \eqref{new_modes} and the standard Fourier modes, it can be shown that the coarse-grid correction operator $C_h = I_h - P_{2h}^h L_{2h}^{-1} R_h^{2h} L_h$, where $P_{2h}^h$ and $R_h^{2h}$ are the prolongation and restriction operators, $L_h$ and $L_{2h}$ are the fine- and coarse-grid operators and $I_h$ is the identity, satisfies the following invariance property: $$C_h: {\mathcal F}_h^{n^2}({\boldsymbol \theta}^{00}) \; \rightarrow {\mathcal F}_h^{n^2}({\boldsymbol \theta}^{00}).$$ More concretely, for ${\boldsymbol \theta}^{00}\in{\boldsymbol \Theta}_{2nh} = (-\pi/2nh, \pi/2nh]^2$, the following properties of the operators in $C_h$ are fulfilled:
\begin{enumerate}
\item $L_h,\, I_h: {\mathcal F}_h^{n^2}({\boldsymbol \theta}^{00}) \; \rightarrow {\mathcal F}_h^{n^2}({\boldsymbol \theta}^{00}),$
\item $L_{2h}: {\mathcal F}_{2h}^{n^2/4}(2{\boldsymbol \theta}^{00}) \; \rightarrow {\mathcal F}_{2h}^{n^2/4}(2{\boldsymbol \theta}^{00}),$
\item $R_h^{2h}: {\mathcal F}_h^{n^2}({\boldsymbol \theta}^{00}) \; \rightarrow {\mathcal F}_{2h}^{n^2/4}(2{\boldsymbol \theta}^{00}),$
\item $P_{2h}^h: {\mathcal F}_{2h}^{n^2/4}(2{\boldsymbol \theta}^{00}) \; \rightarrow {\mathcal F}_h^{n^2}({\boldsymbol \theta}^{00}).$
\end{enumerate}
From these invariance properties we can compute the Fourier representations of the corresponding operators. As an example, we will describe next the representation of $R_h^{2h}$ with respect to the grid-functions $\{\psi_{2h}^{kl}(2{\boldsymbol \theta}^{00})\}_{k,l=0}^{n/2-1}$ and $\{\psi_{h}^{kl}({\boldsymbol \theta}^{00})\}_{k,l=0}^{n-1}$, for the restriction operators considered in this work.
%\begin{example}

We first consider the basic restriction operator obtained as the adjoint of the piecewise constant prolongation operator with stencil form \eqref{cr}. Its Fourier representation with respect to the new Fourier basis is given by
$$
\widehat{R}_h^{2h} ({\boldsymbol \theta}^{00}) = 
\frac{1}{4}\left(
\begin{array}{cccc}
e^{-\imath (\theta_x^{00}+\theta_y^{00})/2} &
e^{\imath (\theta_x^{00}+ \theta_y^{00})/2} &
e^{\imath (\theta_x^{00}-\theta_y^{00})/2} &
e^{\imath (\theta_y^{00}-\theta_x^{00})/2} 
\end{array}
\right).
$$ 
In the case of the cell-centered restriction operator by Wesseling/Khalil~\cite{MG2}, that is,
\begin{equation}
R_h^{2h} = \frac{1}{16}
\left[
\begin{array}{ccccc}
1 & 1 & & 0 & 0 \\ 
1 & 3 & & 2 & 0 \\ 
& & \star & & \\ 
0 & 2 & & 3 & 1 \\ 
0 & 0 & & 1 & 1
\end{array}
\right]_h^{2h},
\label{WKstencil}
\end{equation}
the Fourier representation is given by
$$
\widehat{R}_h^{2h} ({\boldsymbol \theta}^{00}) = 
\frac{1}{16}
\left(
\begin{array}{c}
e^{-\imath (\theta_x^{00}+\theta_y^{00})/2} (2 + e^{2 \imath \theta_y^{00}} + 
e^{2 \imath \theta_x^{00}}) \\
e^{\imath (\theta_x^{00}+ \theta_y^{00})/2} (2 + e^{-2 \imath \theta_y^{00}} + 
e^{-2 \imath \theta_x^{00}}) \\
e^{\imath (\theta_x^{00}-\theta_y^{00})/2}(3 + e^{2 \imath (\theta_y^{00} - \theta_x^{00})} \\
e^{\imath (\theta_y^{00}-\theta_x^{00})/2}(3 + e^{2 \imath (\theta_x^{00} - \theta_y^{00})}).
\end{array}
\right)^T.
$$ 
%\end{example}
As an immediate consequence of these invariance properties and the invariance property of the smoothing operator, also the two-grid operator $K_h^{2h} = C_hS_h^{\nu}$, where $\nu$ denotes the number of smoothing steps, leaves the subspaces ${\mathcal F}_h^{n^2}({\boldsymbol \theta}^{00})$ invariant.  Its Fourier representation is given by 
$$
\widehat{K}_h^{2h}({\boldsymbol \theta}^{00}) = \widehat{C}_h({\boldsymbol \theta}^{00}) \widehat{S}_h^{\nu}({\boldsymbol \theta}^{00}) = (\widehat{I}_h({\boldsymbol \theta}^{00}) - \widehat{P}_{2h}^h({\boldsymbol \theta}^{00}) \widehat{L}_{2h}^{-1}({\boldsymbol \theta}^{00}) \widehat{R}_h^{2h}({\boldsymbol \theta}^{00}) \widehat{L}_h({\boldsymbol \theta}^{00})) \widehat{S}_h^{\nu}({\boldsymbol \theta}^{00}).
$$
Finally, we can compute the asymptotic two-grid convergence factor as the supremum of the spectral radii of $(n^2\times n^2)-$matrices, as follows
$$
\rho(K_h^{2h}) = \sup_{{\boldsymbol \theta}^{00}\in\widetilde{{\boldsymbol \Theta}}_{2nh}}\rho(\widehat{K}_h^{2h}({\boldsymbol \theta}^{00})),
$$
where $\widetilde{{\boldsymbol \Theta}}_{2nh}$ is the subset of ${\boldsymbol \Theta}_{2nh}$ in which we remove the frequencies ${\boldsymbol \theta}^{00}$ such that the determinant of the Fourier symbol of $L_h$ or $L_{2h}$ vanishes.

\section{LFA results for PDEs with jumping coefficients}\label{sec:4}

\setcounter{section}{4}

In this section, we apply the proposed LFA to predict the multigrid convergence factors for a collection of benchmark problems with jumping coefficients taken from the literature \cite{Alcouffe:1981:MGM,knapek,moulton1}. The test cases cover a variety of possible inhomogeneities including jumps that are not aligned with the coarse grid. In all these problems, equation \eqref{diffusion_problem} is numerically solved in the domain ${D}= (0,1)^2$, by using a mesh of grid-size $h = 1/128$. In particular, the following jumping coefficient benchmark problems, characterized by the distribution of the diffusion coefficient, are considered here: \\

\begin{figure}[H]
\begin{center}
\begin{tabular}{ccc}
&\textit{\footnotesize Multigrid} & \textit{\footnotesize LFA}\\
\rotatebox{90}{\hspace{0.7cm}\textit{\footnotesize (a) Vertical jump }}&\includegraphics[width = 0.3\textwidth]{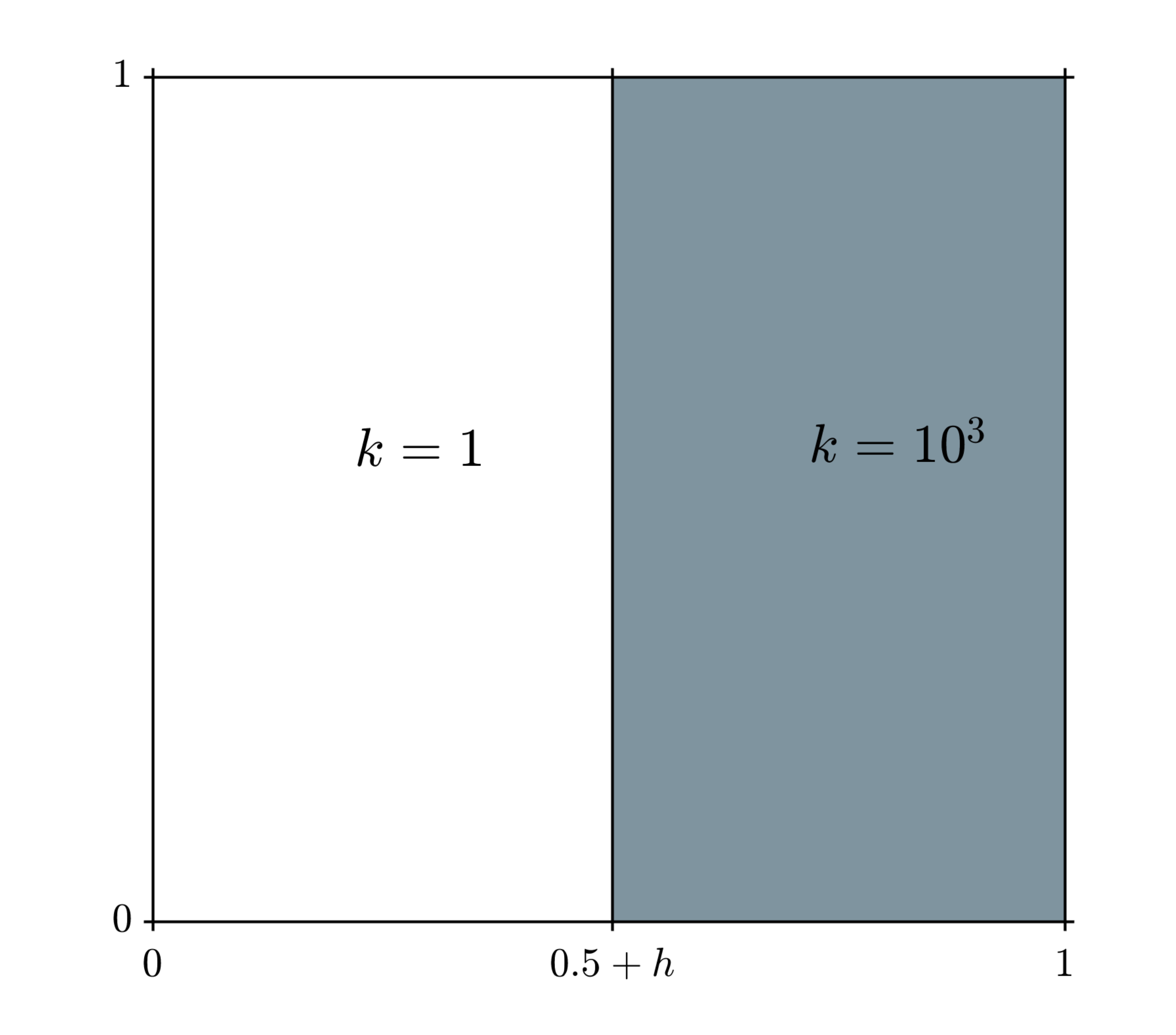} &\includegraphics[width = 0.32\textwidth]{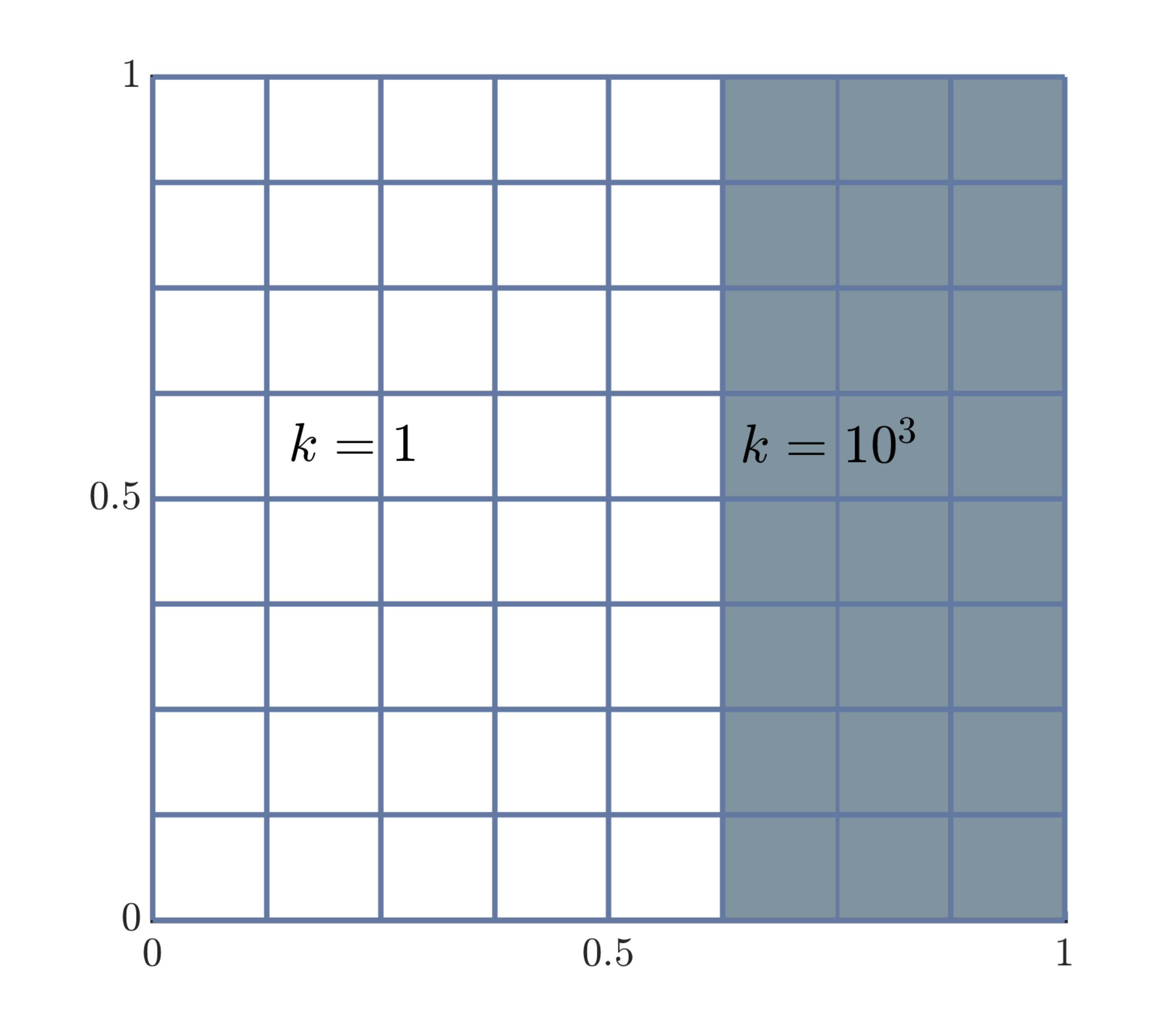} \\
\rotatebox{90}{\hspace{0.7cm}\textit{\footnotesize (b) Four corner problem}}&\includegraphics[width = 0.3\textwidth]{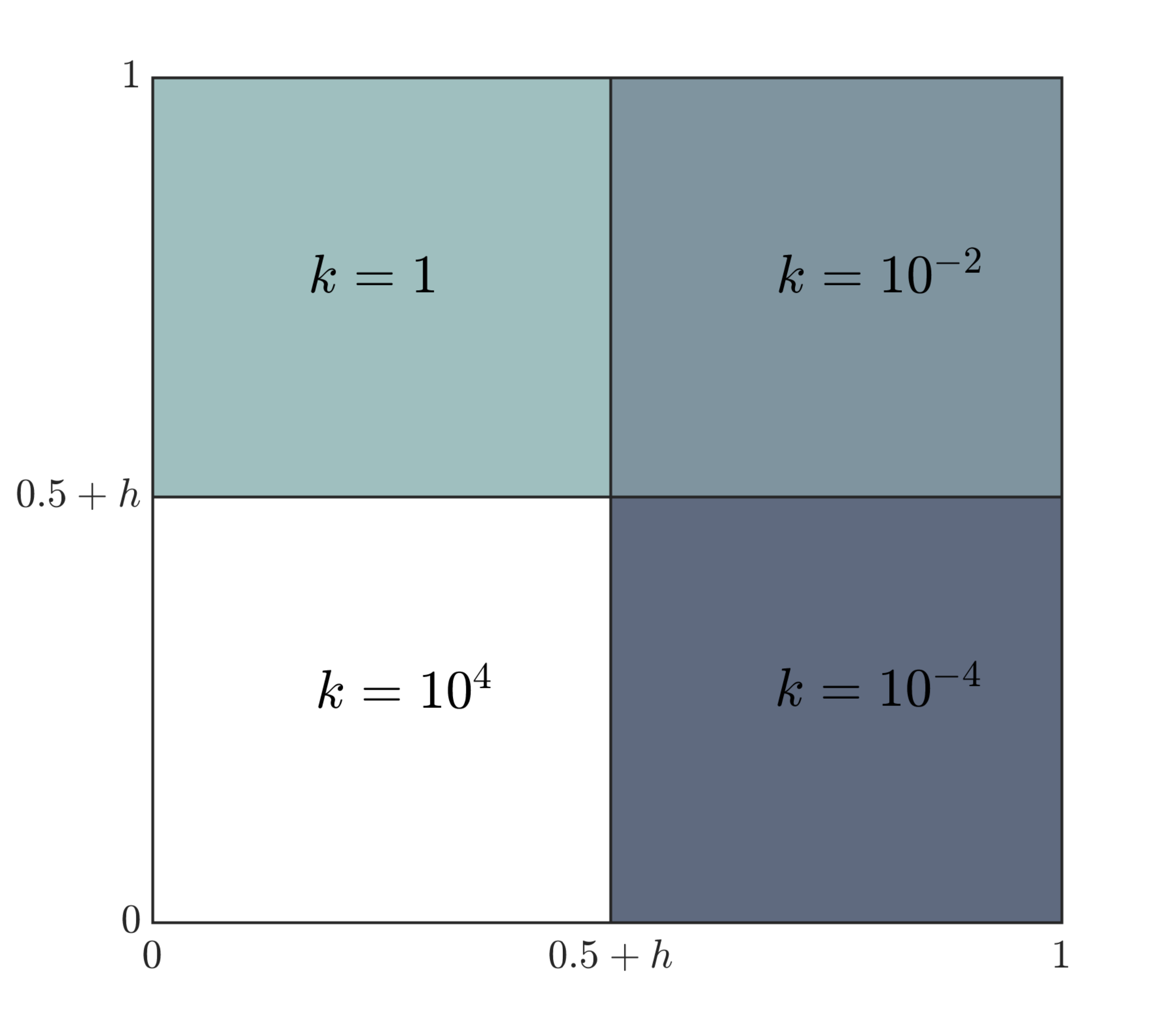}&\includegraphics[width = 0.32\textwidth]{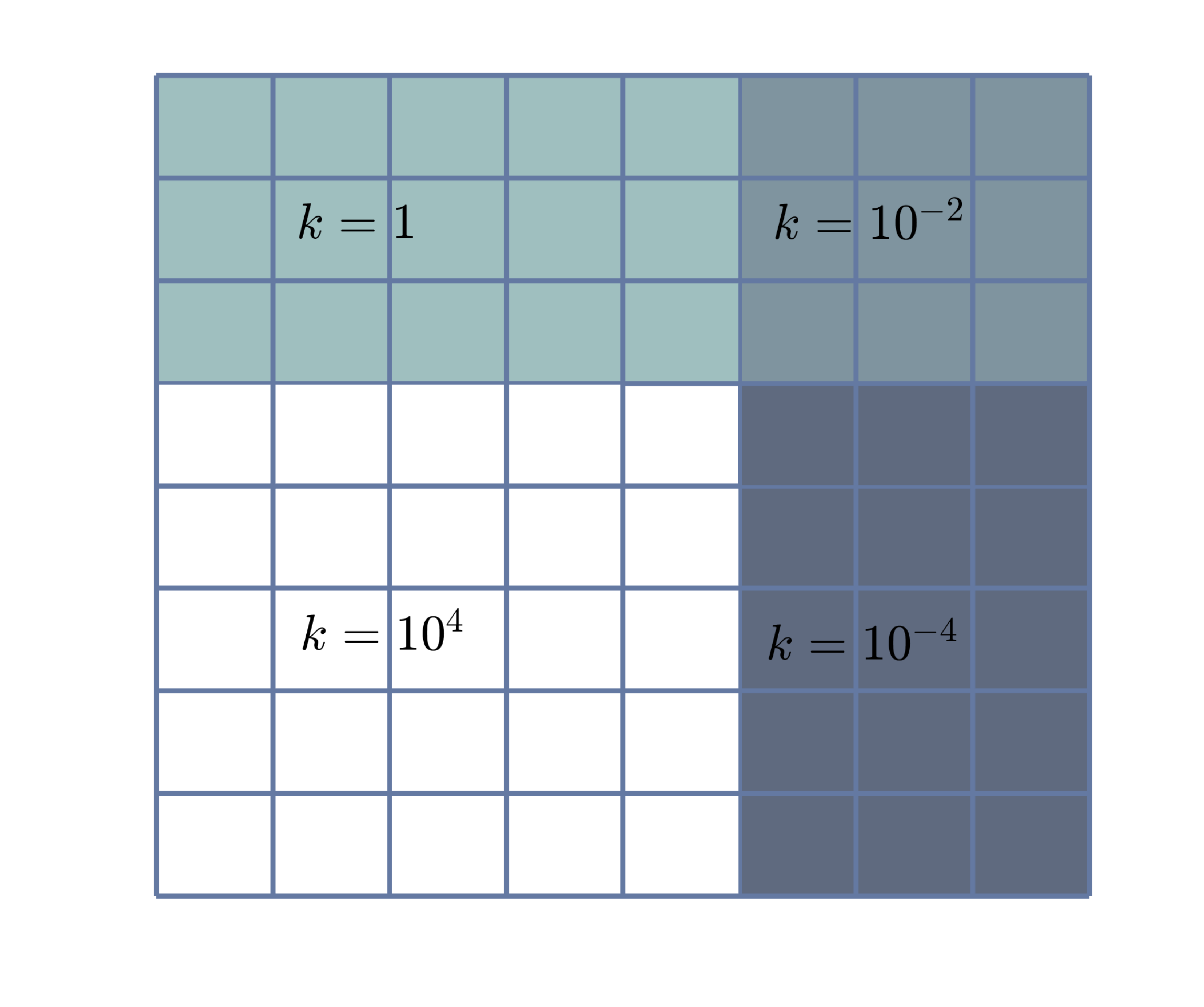}\\
%(a) Chessboard  & (b) Vertical jump \\
\rotatebox{90}{\hspace{0.3cm}\textit{\footnotesize (c) Square inclusions}}&\includegraphics[width = 0.3\textwidth]{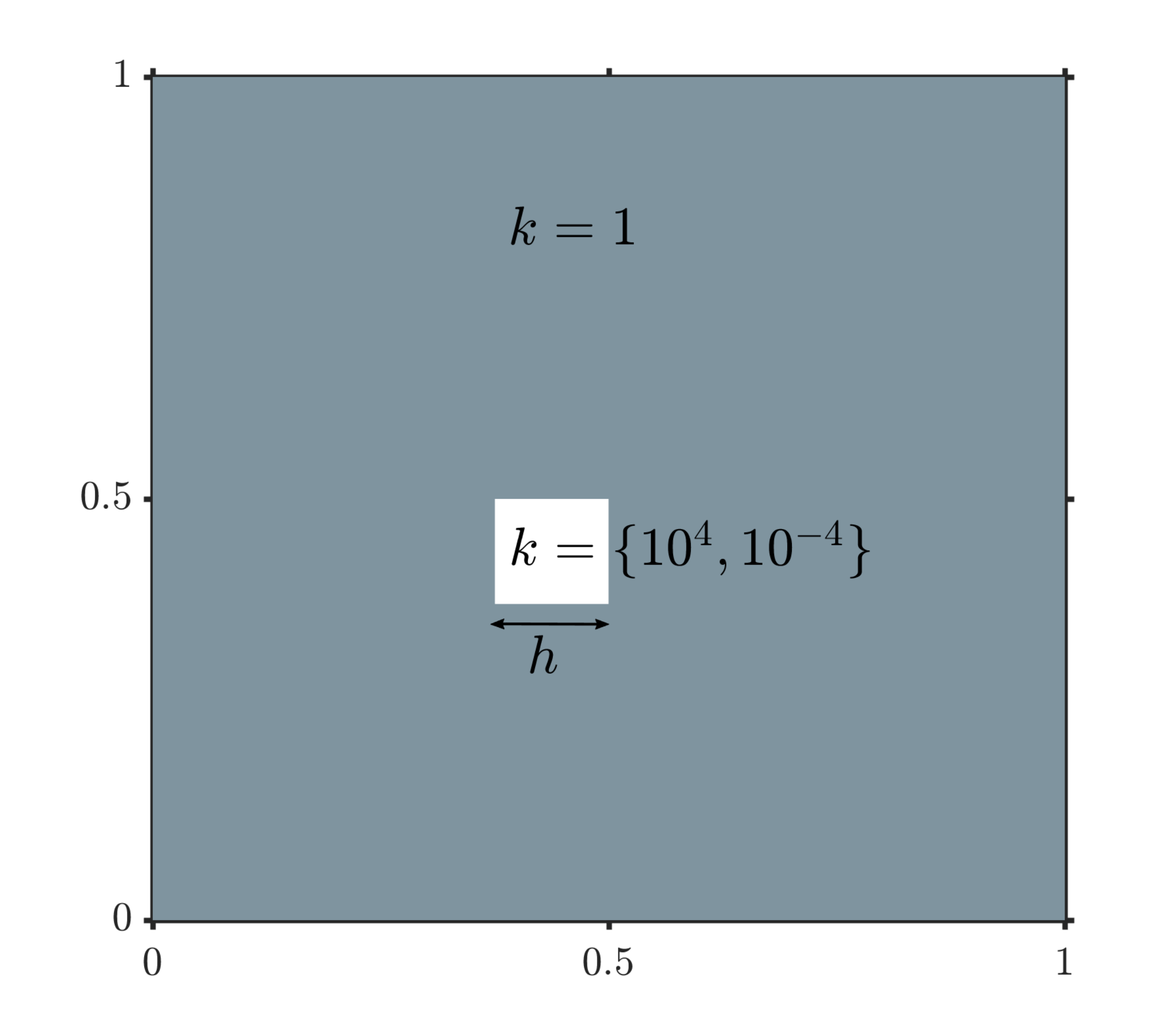}&\includegraphics[width = 0.32\textwidth]{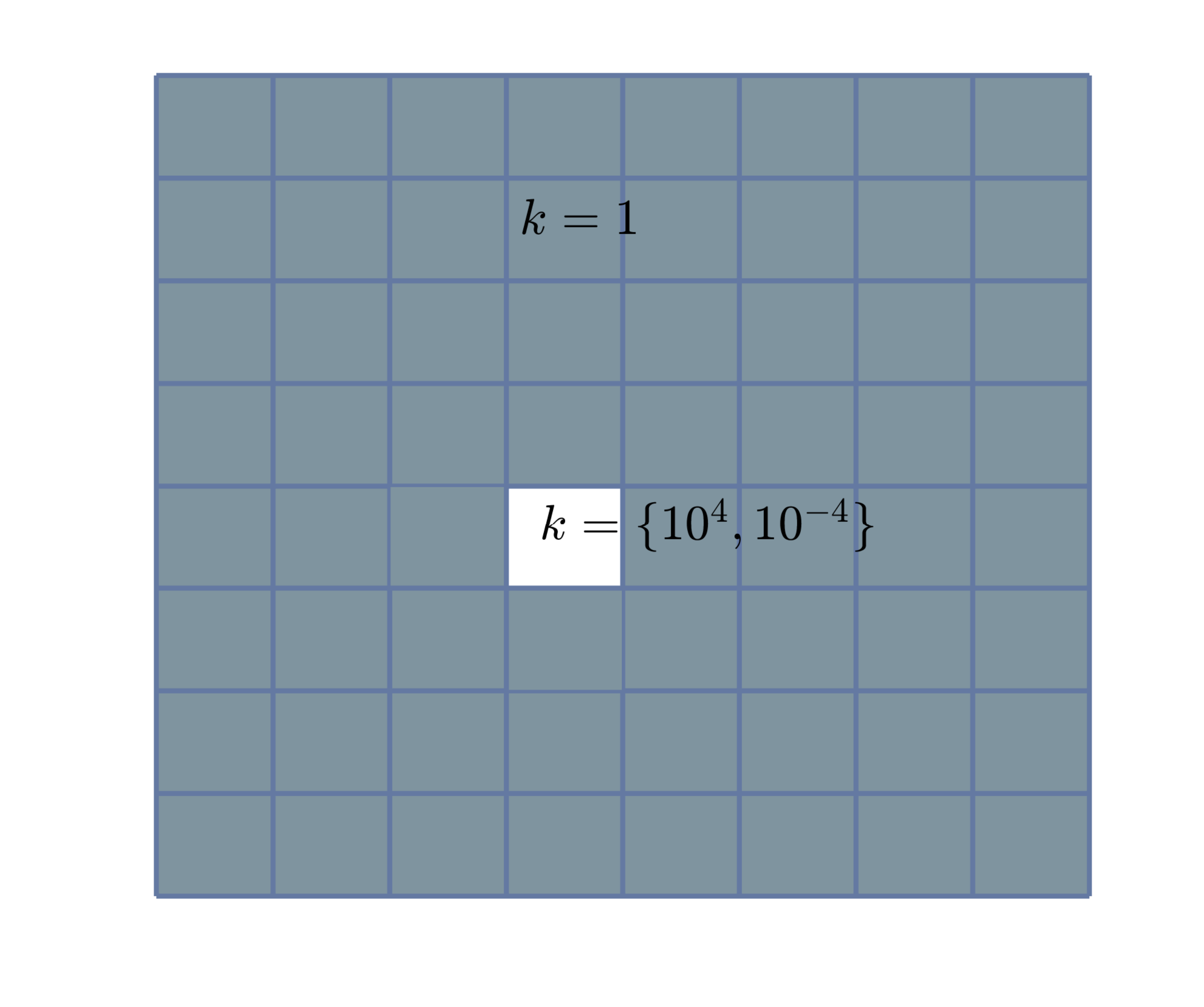}\\
\rotatebox{90}{\hspace{0.1cm}\textit{\footnotesize (d) Periodic Sq. inclusions}}&\includegraphics[width = 0.3\textwidth]{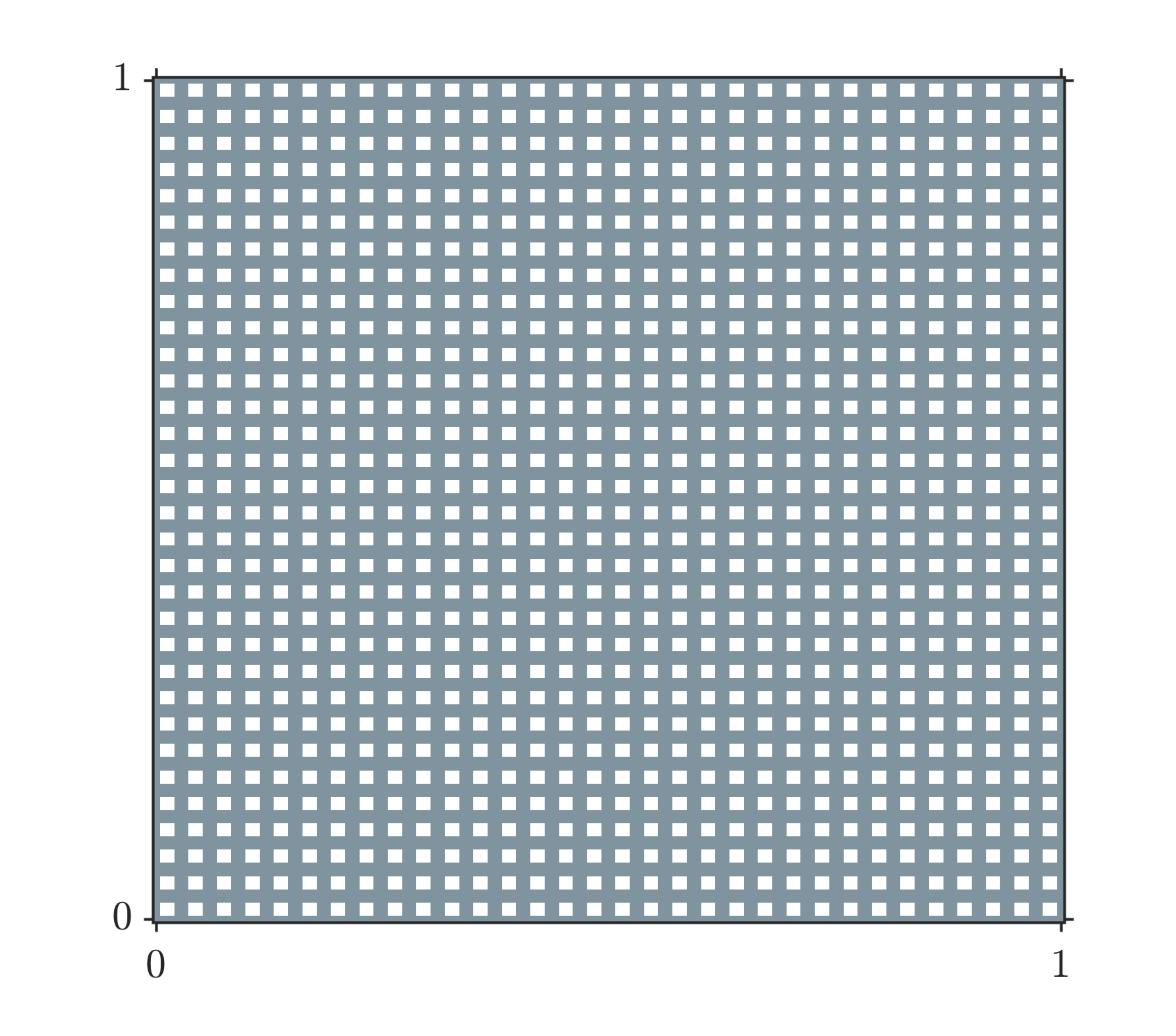}&\includegraphics[width = 0.32\textwidth]{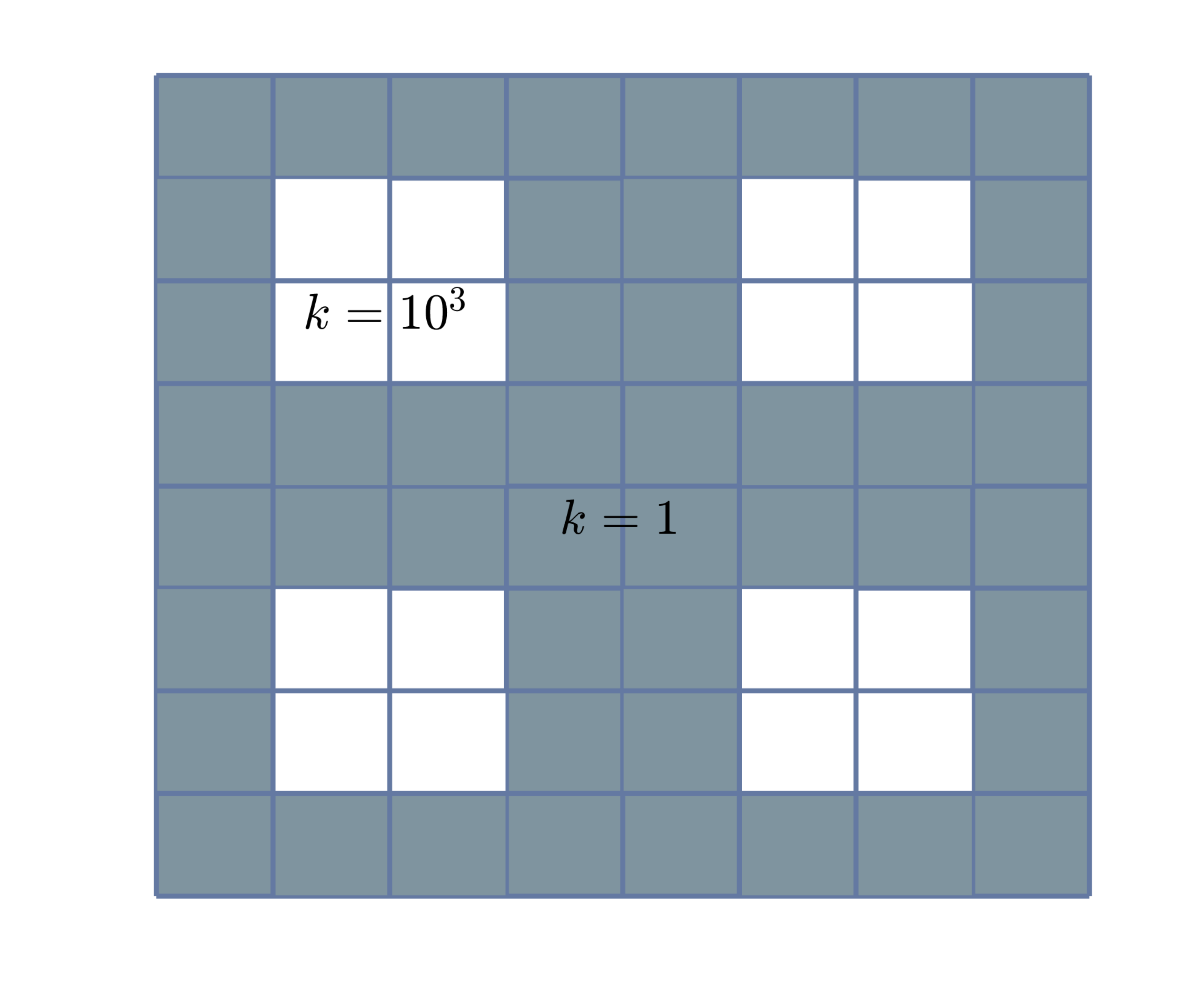}\\
%(c) Square inclusions & (d) Four corner problem\\
\rotatebox{90}{\hspace{0.2cm}\textit{\footnotesize (e) Periodic L-s. inclusions}}&\includegraphics[width = 0.3\textwidth]{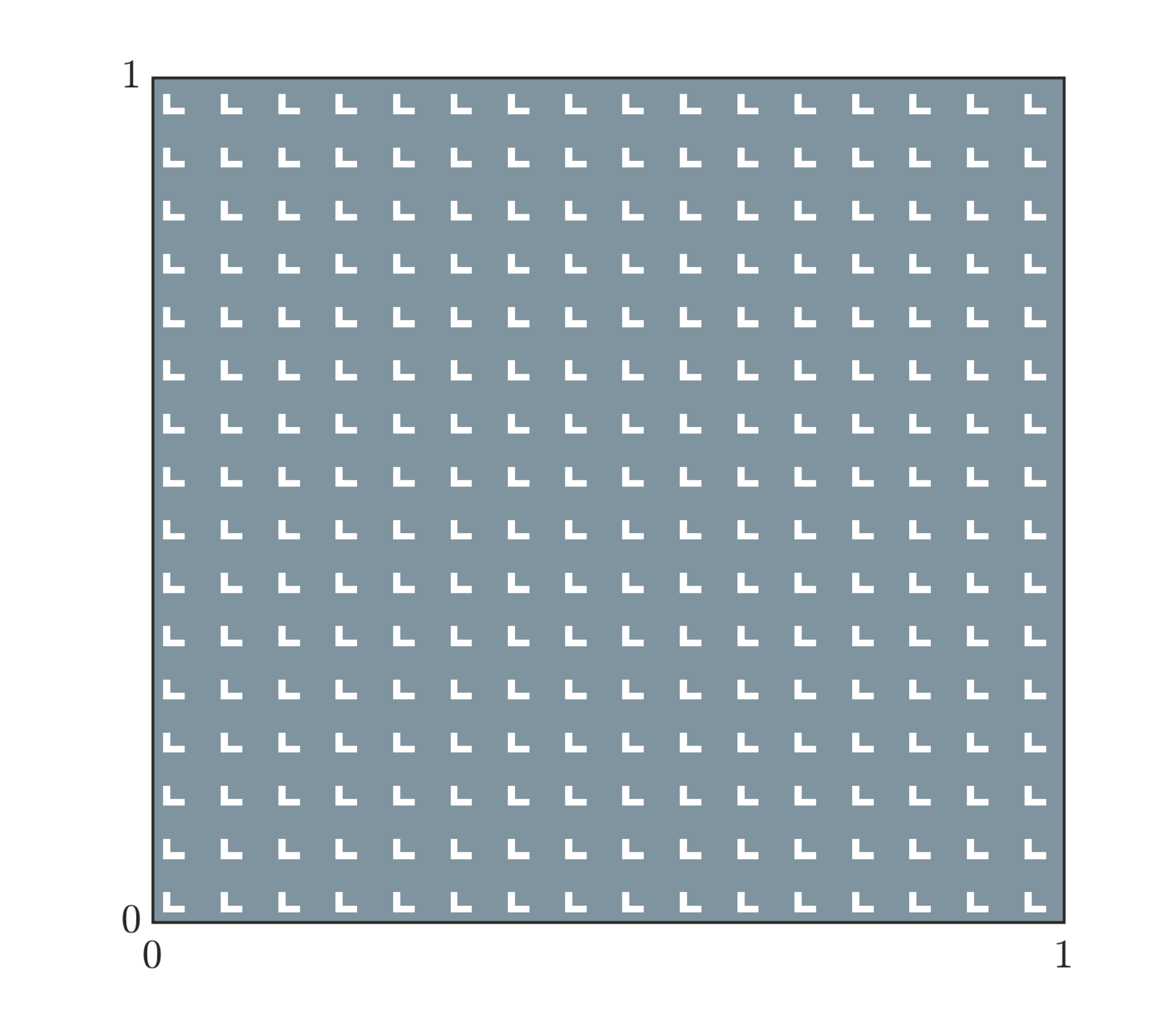}&\includegraphics[width = 0.32\textwidth]{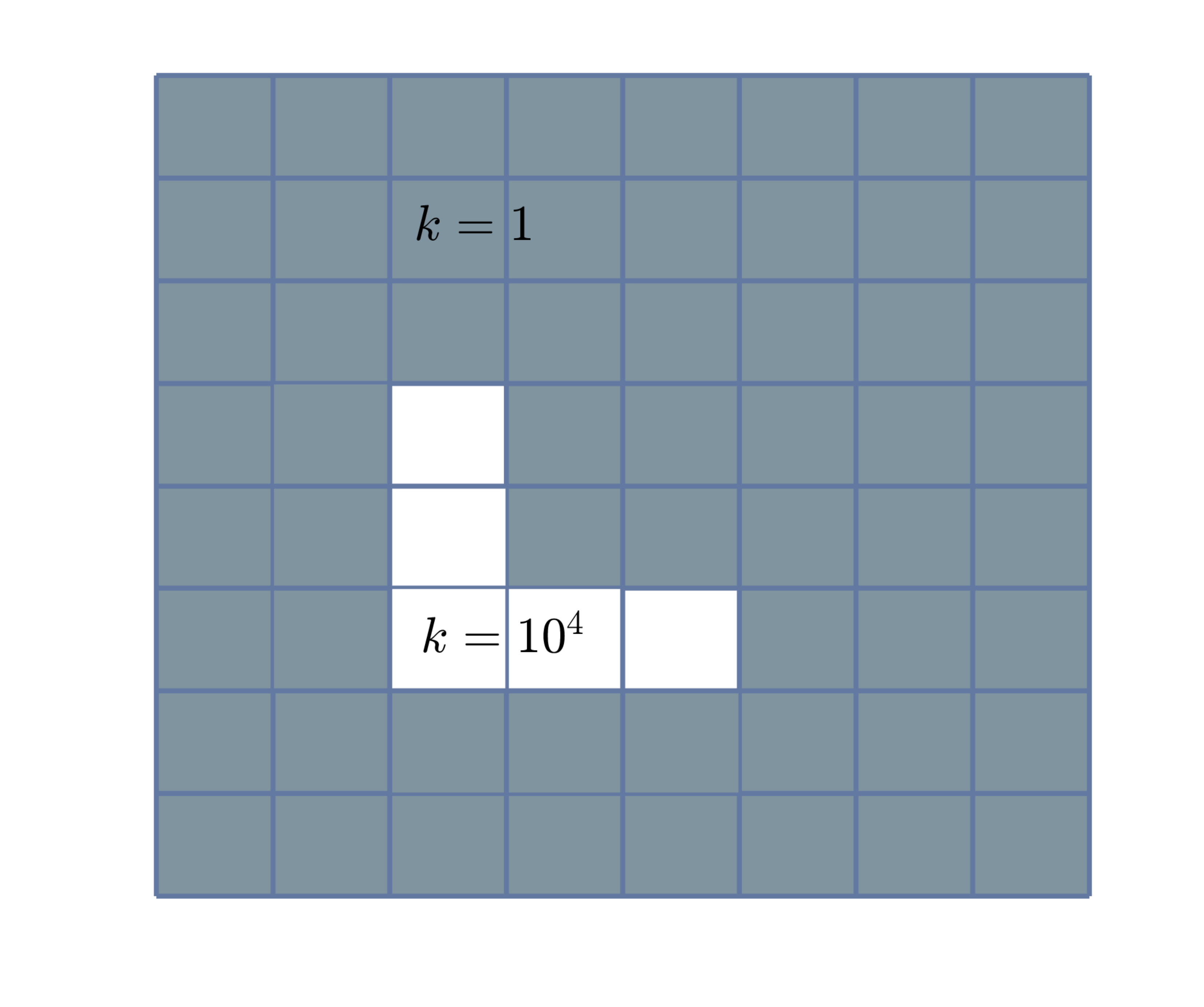}\\
%(e) L-shaped inclusion
\end{tabular}

\caption{Distribution of the diffusion coefficients for the five considered examples on a unit square domain and corresponding $8\times8$ window used in the LFA.}\label{figures_examples}
\end{center}
\end{figure}

\begin{enumerate}
\item {\bf Vertical jump.}
%This case is chosen such that the interface position ``moves'' after the first coarsening step and is subsequently eliminated on further coarsening. 
Function $k(x,y)$ is defined in the following way (see also Figure \ref{figures_examples} (a))
$$
k(x,y) = 
\left\{
\begin{array}{l}
1, \quad \quad{\mbox {\rm if} } \quad x < \frac12+h, \\ [1ex]
10^3, \quad {\mbox {\rm if} } \quad x \geq \frac12+h.
\end{array}
\right.
$$
\item {\bf Four corner problem.} 
The domain is divided into four regions in which the diffusion coefficient is varying, see Figure \ref{figures_examples} (b). In particular,
%This results in a cross-point formed at the intersection of the four interfaces. Furthermore, we move the cross-point by one cell from the center so that it disappears after coarsening. In detail,
$$
k(x,y) = 
\left\{
\begin{array}{l}
10^4, \quad {\mbox {\rm if} } \quad (x,y)\in\left(0,\frac12+h\right)^2, \\ [1ex]
1, \qquad {\mbox {\rm if} } \quad (x,y)\in\left(0,\frac12+h\right)\times\left(\frac12+h,1\right),\\ [1ex]
10^{-2} \quad {\mbox {\rm if} }\quad (x,y)\in\left(\frac12+h,1\right)\times\left(\frac12+h,1\right), \\ [1ex]
10^{-4} \quad \text{otherwise}.
\end{array}
\right.
$$
\item {\bf Square inclusion.}
In this example we assume a square inhomogeneity in one cell within the square domain, see Figure \ref{figures_examples} (c). 
%The position of the inclusion is chosen such that the interfaces disappear after subsequent coarsening steps. 
The diffusion coefficient is defined as
$$
k(x,y) = 
\left\{
\begin{array}{l}
k_0, \quad {\mbox {\rm if} } \quad (x,y)\in\left(\frac12-h,\frac12\right)^2, \\ [1ex]
1, \quad {\mbox {\rm otherwise}},
\end{array}
\right.
$$
where values $k_0=10^4$ and $k_0=10^{-4}$ are considered.
%Furthermore, as the effective permeability values depend on the magnitude of the inclusion relative to the background value of the coefficients \cite{moulton1}, we consider two cases with inclusions $k=10^4$ and $k=10^{-4}$. 
\item {\bf Periodic square inclusions.}
This test is taken from \cite{scott}. We consider a structured pattern of square inclusions of size $2h\times 2h$ as depicted in Figure \ref{figures_examples} (d). The diffusion parameter is $k(x,y)=1$ inside the dark region and $k(x,y) = 1000$ inside the white region.
\item {\bf Periodic L-shaped inclusions.}  In the last test case, we consider a structured pattern of L-shaped inclusions as in Figure \ref{figures_examples} (e). The diffusion parameter is $k(x,y)=10^4$ inside the white region and $k(x,y) = 1$ inside the dark region. 
\end{enumerate}
To perform the theoretical analysis, the periodic extension of a window of size $8\times8$ has been chosen, where the diffusion coefficient is prescribed in such window according to its definition, see Figure \ref{figures_examples} (right side). In all numerical tests a random initial guess is chosen, and the right-hand side and boundary conditions are set to zero to be able to determine asymptotic convergence factors. In this way, we avoid round-off errors permitting us to perform as many iterations as needed. In practice, we have seen that $50$ iterations are sufficient. \\
Next, we show the excellent correspondence between the theoretical analysis and the experimental results for these test cases. Two combinations of inter-grid transfer operators are considered. The first combination, denoted here by (CP,CR), is based on the use of piecewise constant prolongation \eqref{cp} and its adjoint as the restriction \eqref{cr}. In the second combination we change to a higher polynomial order interpolation operator which is the adjoint to the Wesseling/Khalil restriction \eqref{WKstencil}. This choice is denoted by (WP,CR). Moreover, a standard damped Jacobi smoother (damping with $\omega = 0.8$) is considered as well as the proposed lexicographic Gauss-Seidel smoother.
\begin{table}[hbt]
\begin{center}
\begin{tabular}{|c|c|c|ccc|}
\cline{4-6}
\multicolumn{3}{c|}{}& \multicolumn{3}{c|}{Pre-, post-smoothing steps} \\
\cline{4-6}
\multicolumn{3}{c|}{}& $(1,0)$ & $(1,1)$ & $(2,2)$ \\
\hline
\multirow{4}{*}{Vertical jump} & \multirow{2}{*}{GS} &  (CP,CR) &  0.42(0.42)&0.18(0.19)&0.04(0.03) \\
& & (WP,CR) &  0.41(0.37)&0.19(0.16)&0.07(0.11) \\
\cline{2-6}
&  \multirow{2}{*}{Jac.} & (CP,CR)  & 0.65(0.65)&0.43(0.42)&0.19(0.19) \\
& & (WP,CR)& 0.63(0.59)&0.40(0.35)&0.19(0.19) \\
\hline
\multirow{4}{*}{Four corner problem} & \multirow{2}{*}{GS} & (CP,CR) &  0.42(0.37)&0.15(0.12)&0.04(0.03) \\
& & (WP,CR) &  0.40(0.39)&0.16(0.16)&0.09(0.09)  \\
\cline{2-6}
& \multirow{2}{*}{Jac.} & (CP,CR) & 0.63(0.61)&0.40(0.39)&0.16(0.16)  \\
& & (WP,CR) &  0.62(0.62)&0.40(0.40)&0.19(0.18)  \\
\hline
\multirow{4}{*}{Sq. inclusion ($k=10^4$)} & \multirow{2}{*}{GS} & (CP,CR) & 0.45(0.44)&0.21(0.19)&0.04(0.04) \\
& & (WP,CR) & 0.41(0.40)&0.18(0.17)&0.11(0.11) \\
\cline{2-6}
& \multirow{2}{*}{Jac.} & (CP,CR) & 0.60(0.65)&0.36(0.42)&0.13(0.19)  \\
& & (WP,CR) & 0.60(0.62)&0.38(0.37)&0.22(0.21) \\
\hline
\multirow{4}{*}{Sq. inclusion ($k=10^{-4}$)} & \multirow{2}{*}{GS} &  (CP,CR) & 0.46(0.45)&0.21(0.20)&0.05(0.05)  \\
& & (WP,CR) &  0.41(0.40)&0.19(0.19)&0.12(0.12) \\
\cline{2-6}
& \multirow{2}{*}{Jac.} & (CP,CR) &  0.61(0.65)&0.38(0.42)&0.15(0.19)  \\
& & (WP,CR) &  0.61(0.59)&0.39(0.39)&0.23(0.23)  \\
\hline
\multirow{4}{*}{Periodic Sq. inclusions} & \multirow{2}{*}{GS} & (CP,CR)  & 0.64(0.61)&0.43(0.42)&0.41(0.41)  \\
& & (WP,CR) &  0.62(0.61)&0.43(0.41)&0.41(0.40) \\
\cline{2-6}
& \multirow{2}{*}{Jac.} & (CP,CR) & 0.81(0.78)&0.66(0.65)&0.44(0.46) \\
& & (WP,CR) & 0.81(0.80)&0.66(0.65)&0.44(0.43) \\
\hline
\multirow{4}{*}{Periodic L-S. inclusions} & \multirow{2}{*}{GS} & (CP,CR) & 0.50(0.50)&0.32(0.26)&0.21(0.21)  \\
& & (WP,CR) &  0.54(0.53)&0.40(0.40)&0.30(0.30) \\
\cline{2-6}
& \multirow{2}{*}{Jac.} & (CP,CR) & 0.71(0.63)&0.54(0.48)&0.36(0.35)  \\
& & (WP,CR) & 0.71(0.71)&0.56(0.56)&0.42(0.42)  \\
\hline
\end{tabular}
\caption{Asymptotic two-grid convergence factors predicted by LFA  and the corresponding computed average multigrid convergence factors (in parenthesis) using two-grid cycles using different pre- and post-smoothing steps for the five examples.}
\label{table_1d}
\end{center}
\end{table}

In Table \ref{table_1d}, for different numbers of smoothing steps, for two different smoothers and for the two combinations of restriction and prolongation operators, we provide the two-grid convergence factors predicted by the novel LFA for each of the proposed numerical experiments.  We also display in parenthesis the average after $50$ iterations of the experimentally computed multigrid convergence factors by using two-grid cycles. For all these cases, we observe a very accurate match between the analysis results and the rates experimentally obtained. Regarding the size of the window 
to perform the LFA, we have observed that a window of size $8\times8$ is enough to achieve 
excellent predictions in all considered benchmark problems. For example, for the vertical jump test, the two-grid analysis considering four smoothing steps of Gauss-Seidel smoother and the combination (CP,CR) of inter-grid transfer operators provides a factor of $0.11$ when a $2\times2$ window is used, a factor of $0.06$ for a $4 \times4$ window, and a factor of $0.04$, which matches perfectly the real convergence, when the $8\times 8$ window is considered.

We remark that with the current multigrid approach, the quality of the coarse grid discretization may not be satisfactory, for example, for chessboard or L-shaped inclusion. For such cases, we either recommend using more \emph{powerful} smoothers such as the ILU smoother or adapt the coarsening. Furthermore, for PDEs with strong local variations in the coefficient fields, \emph{homogenization techniques} \cite{Alcouffe:1981:MGM,knapek, moulton1} may also be used to obtain coarser representation of the fine grid problem.
Comparing the results of the two combinations of inter-grid transfer operators, we observe a very similar performance for all five test cases studied here. In the rest of the paper, we therefore choose the strategy (CP,CR) because of its simplicity and low computational cost.
%%%%%%%%%%%%%%%%%%%%%%%%%%%%%%%%%%%%%%%%%%%%%%%%%%
\section{LFA results for PDEs with random coefficients}\label{sec:5}

\setcounter{section}{5}
Here, we consider the SPDE \eqref{1a} defined on a unit square domain ${D} = (0,1)^2$ with homogeneous Dirichlet boundary conditions. Two different types of diffusion coefficients based on random jumps and lognormal random fields are studied. The randomly jumping coefficient problem can be seen as a transition from the deterministic to a stochastic setting. 
\subsection{Randomly jumping coefficients} 
\begin{figure}[hbt]
\begin{center}
\includegraphics[width = 0.55\textwidth]{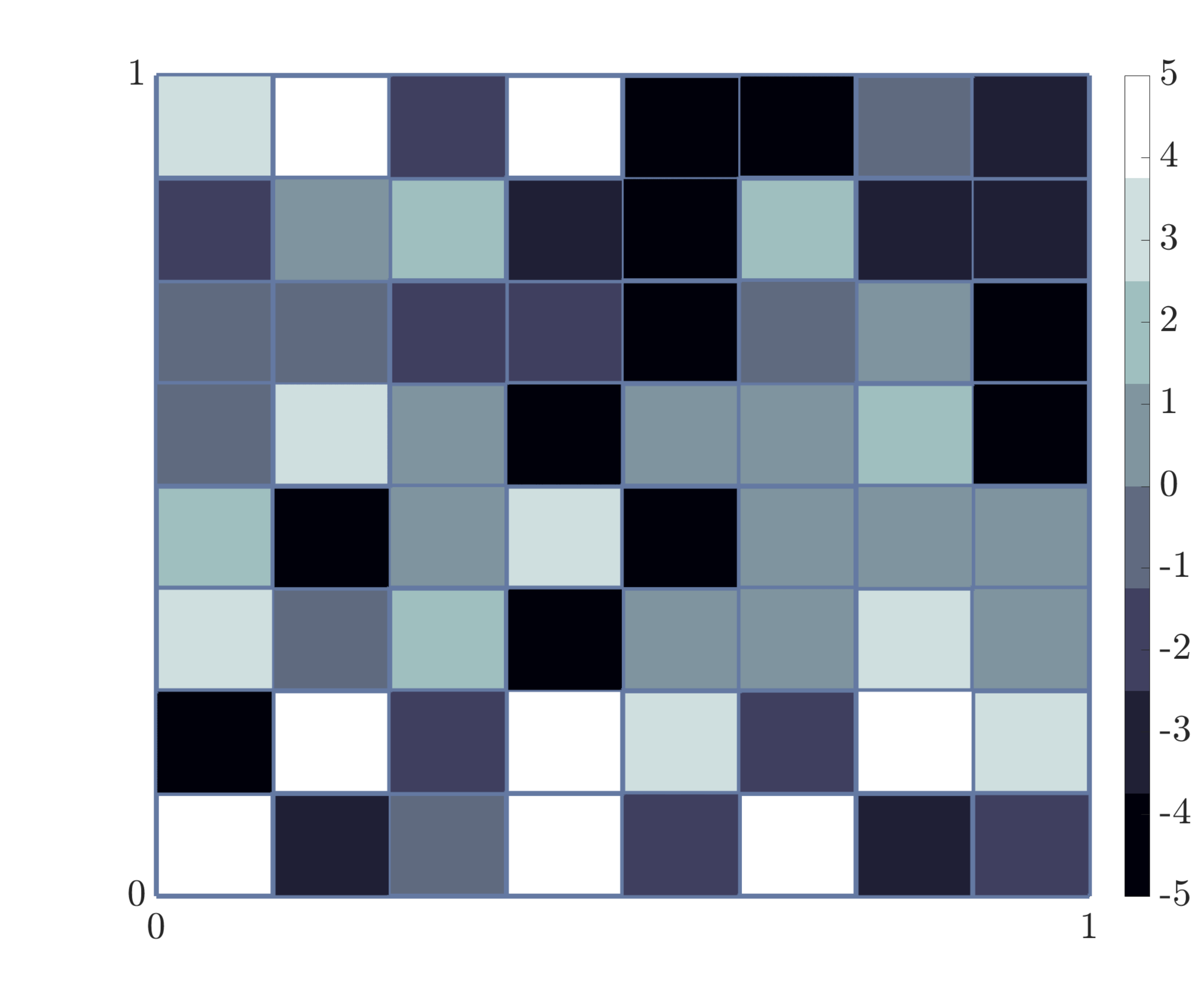}
\end{center}
\caption{An example of random realization of $U = \ln k$ with $m=5$ on a unit square domain.}\label{RJ_sample}
\end{figure}

To simulate random jumps, the domain $D$ is subdivided into square-blocks of size $[\frac18\times\frac18]$ and the value of the coefficients on each of the blocks is sampled as
\begin{equation}\label{RJ}
 k = e^{U}\quad \text{with} \quad U\sim\mathcal{U}\{-m,m\}\quad\text{and}\quad m\in\mathbb{Z}.
\end{equation}
In other words, ${U}$ is an independent identically distributed (i.i.d.) integer sampled from a discrete uniform distribution $\mathcal{U}\{-m,m\}$. Here, the integer $m$ defines the order of magnitude of the jumps, an example for $m=5$ is shown in Figure \ref{RJ_sample}.  Notice that for this choice of $m$, we may encounter interfaces with maximum jumps of magnitude equal to $e^{10}$.
For each random realization of the jumping coefficient field, we compare the LFA two-grid convergence factors with the computed asymptotic convergence factors of the multigrid method by using W-cycles. To perform the LFA, we again use a window of size $8 \times 8$. Furthermore, the LFA equivalent of the randomly jumping coefficient problem is similar to the four corner problem in Figure \ref{figures_examples} (d) with the magnitude of each block given by \eqref{RJ} and the cross-point located exactly at the center of the LFA block. Regarding the multigrid components, a lexicographic Gauss-Seidel iteration is employed as the smoother, and the simplest combination (CP,CR) of inter-grid transfer operators is chosen. Also, we use a $4\times4$ grid as the coarsest in the multigrid hierarchy. To determine the asymptotic convergence factors of the multigrid method, the right-hand side is again set to zero.  

The experimental convergence factor of the multigrid method for the $i$-th realization of the random field is then computed, as follows:
\begin{equation}\label{convFac}
\rho_{i} = \bigg\{\frac{||res^{k_i}||_\infty}{||res^0||_\infty}\bigg\}^{1/k_i},\qquad \text{for}\quad  i=1,2,...,N_{MG},
\end{equation} 
where $||\text{res}^0||_{\infty}$ is the infinity norm of the residual obtained from an initial solution and $||\text{res}^{k_i}||_{\infty}$ is the residual after $k_i$ iterations of the multigrid cycle. We use these quantities to calculate the average and the standard deviation of the asymptotic convergence factors:
\begin{equation}
\langle\rho\rangle_{MG} =\frac{1}{N_{MG}}\sum^{N_{MG}}_{i=1}\rho_{i},\qquad \sigma_{MG} = \sqrt{\frac{1}{(N_{MG}-1)} \sum_{i=1}^{N_{MG}} \left(\rho_{i} - \langle\rho\rangle_{MG}\right)^2},
\end{equation}
respectively. These averaged quantities are defined similarly for the LFA results (based on LFA two-grid factors), and are denoted, respectively, by $\langle\rho\rangle_{LFA}$ and $\sigma_{LFA}$.

\begin{figure}[H]
\begin{center}
\includegraphics[clip, trim=1.2cm 0cm 1.5cm 1cm,scale=.28]{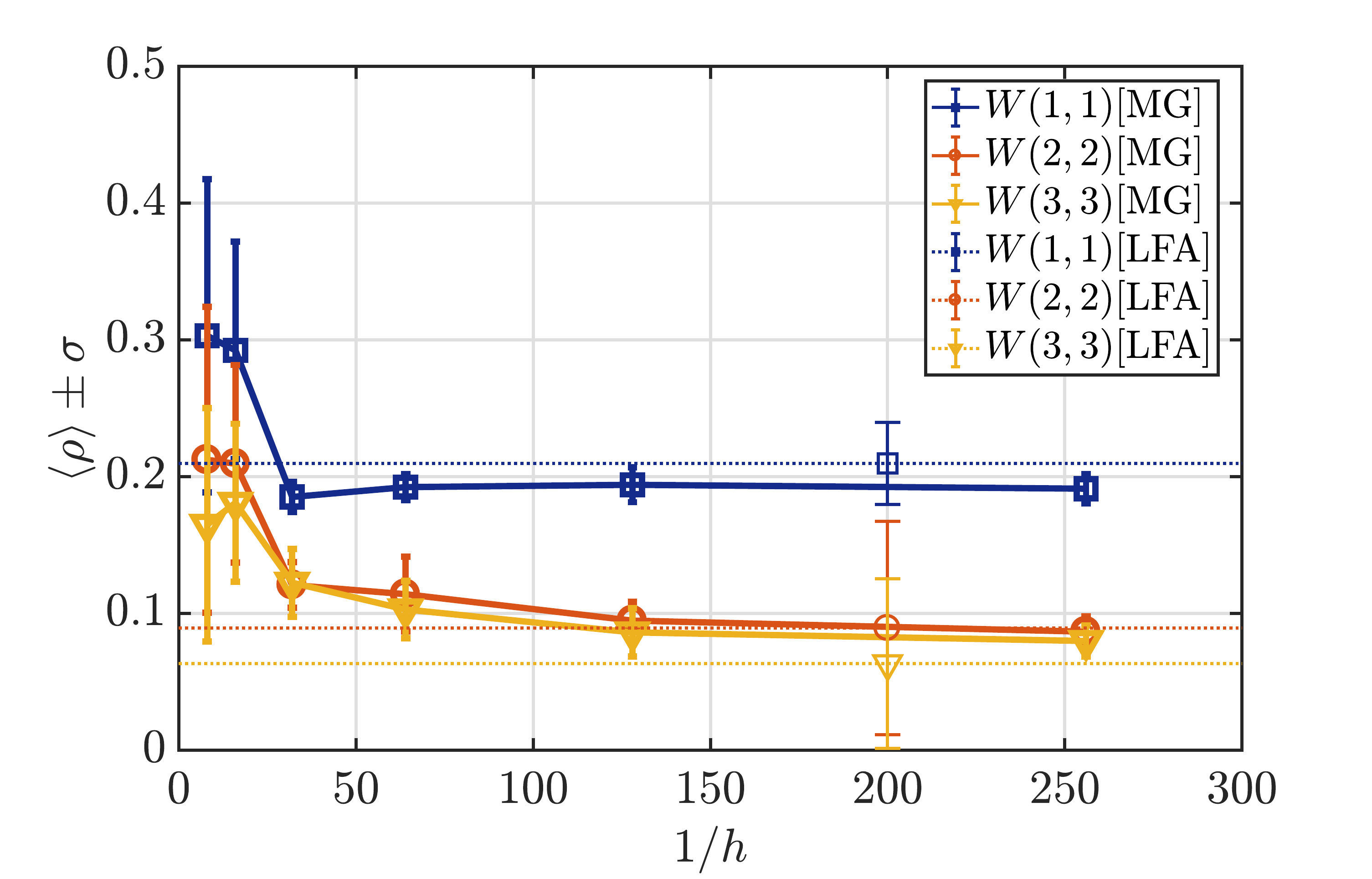} 
\includegraphics[clip, trim=1.2cm 0cm 1.5cm 1cm,scale=.28]{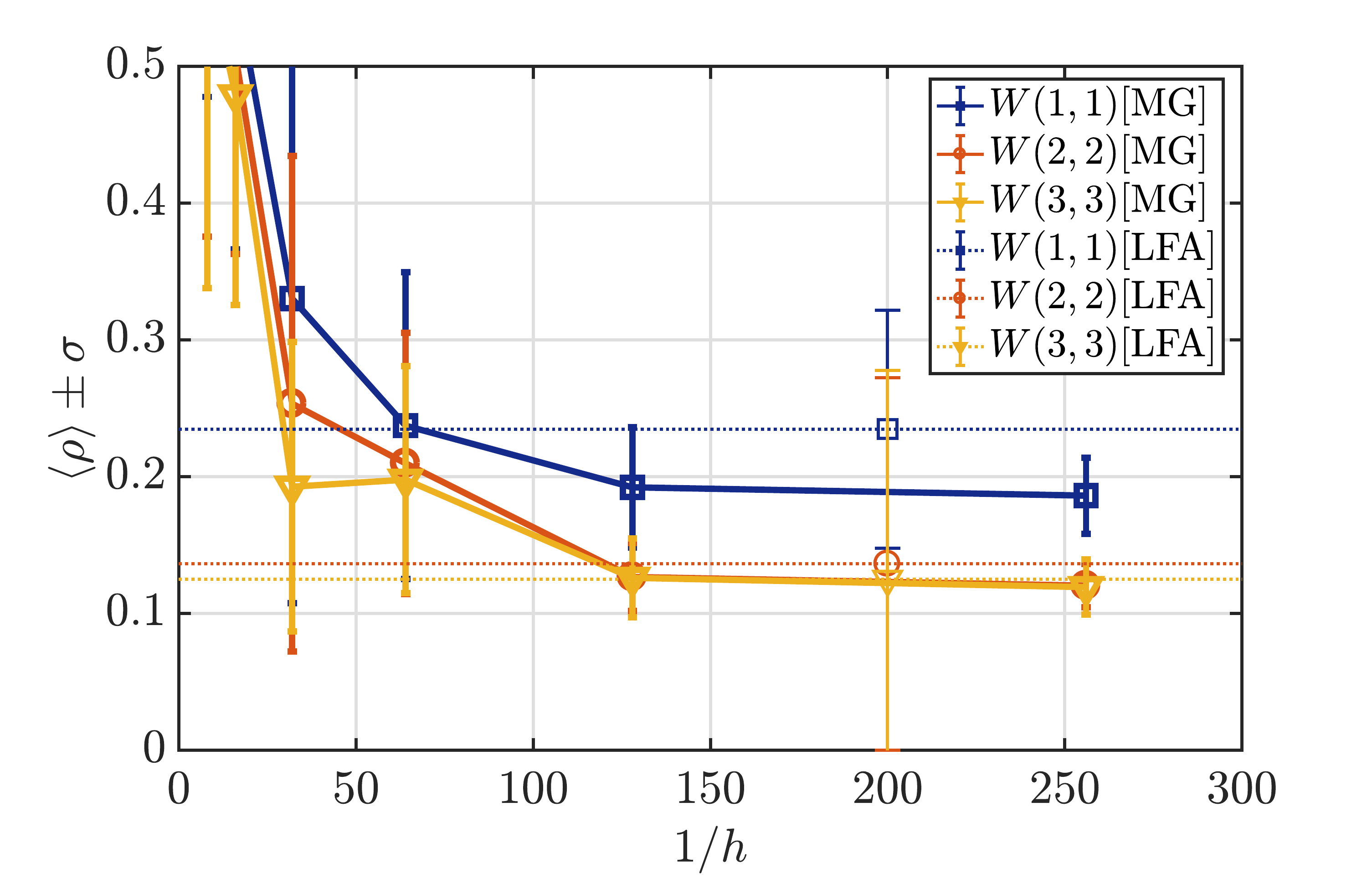} 
\end{center}
\caption{Comparison of the mean and the standard deviation of the LFA (dashed line) and MG (solid line) convergence factors for different $W-$cycling strategy for randomly jumping coefficients with $m=2$ (left) and $m=5$ (right).}
\label{randomjump}
\end{figure}

In Figure \ref{randomjump}, we show the comparison, the mean $\pm$ standard deviation, of the LFA prediction and the multigrid convergence for jump parameter $m=2$ (left) and $m=5$ (right) computed using $N_{LFA}=N_{MG} =100$. Overall, a good match between the LFA and MG convergence is seen up to one decimal place. We also observe that for this specific jumping coefficient problem, there is no further improvement with an increase in the number of smoothing steps after the $W(2,2)-$cycle.

\subsection{Lognormal random fields} Next, we test the LFA prediction capability for a more realistic lognormal diffusion problem. Log-normality leads to positive permeability throughout the domain. The logarithm of the permeability field, $Z=\log k$, is modeled by a zero-mean Gaussian random field $Z: \Omega \times \overline{{ D}} \rightarrow \mathbb{R}$. 
%so that
%\begin{equation}
%\E[Z({\mathbf x}, \cdot) = 0, \mbox{ and } \cov(Z({\mathbf x}, \cdot) Z({\mathbf y}, \cdot)) = \E[Z({\mathbf x}, \cdot)Z({\mathbf y}, \cdot)],\;  {\mathbf x, y} \in \D.
%\end{equation}
A simplification is the use of a homogeneous covariance function $C_\Phi : \mathbb{R} \rightarrow \mathbb{R}$, so that
\begin{equation}
\text{Cov}(Z({\mathbf x}, \cdot), Z({\mathbf y}, \cdot)) = C_\Phi(r), \mbox{  with }r = ||{\mathbf x}-{\mathbf y}||_2.
\end{equation}
By the so-called {\em Mat\'ern family of covariance functions}, random coefficient fields with different degrees of smoothness can be generated.
The Mat\'ern covariance function \cite{MATERN} is characterized by a parameter set $\Phi = (\nu_c, \lambda_c, \sigma_c^2)$, as follows,
\begin{equation}
{C}_\Phi(r) = \sigma_c^2 \frac{2^{1-\nu_c}}{\Gamma(\nu_c)}\left(2\sqrt{\nu_c} \frac{r}{\lambda_c}\right)^{\nu_c} K_\nu\left(2\sqrt{\nu_c} \frac{r}{\lambda_c}\right).
\label{mat}
\end{equation}
Here, $\Gamma$ is the gamma function and $K_\nu$ the modified Bessel function of the second kind. The different parameters have different roles regarding the field's randomness. Parameter $\nu_c$ defines the field's smoothness, $\sigma_c^2$ represents its variance and $\lambda_c$ is the correlation length of the covariance function. Moreover, parameters $\lambda_c$ and $\sigma_c^2$ prescribe the number of peaks and the amplitude of the random field, respectively. When the smoothness parameter is $\nu_c = 1/2$, the Mat\'ern function corresponds to an exponential model, whereas when $\nu_c \rightarrow \infty$ it represents a Gaussian model. 

Realizations of the random field are ``almost surely'' H\"older continuous, $Z,k \in C^{\eta}({D})$, with $0 < \eta < \nu_c$ (see e.g. \cite{Adler_book,Nobile2015}).
For Eq.(\ref{1a}) regularity results have been obtained by taking into account the regularity of the lognormal coefficient field, see \cite{doi:10.1137/110853054, doi:10.1137/080717924, Nobile2015}. 
%\subsection{LFA results} 
\begin{remark}
The random fields are generated using the circulant embedding method \cite{RF4,RF2}. In \emph{Appendix A}, we briefly describe this sampling algorithm. This method employs the fast Fourier transform (FFT) for the covariance matrix decomposition, and, therefore, requires $\mathcal{O}(h^{-2} \log h^{-2})$ operations to generate one sample of the random field on a uniform grid mesh size $h$. Other techniques, e.g., the Karhunen-Lo\'eve (KL) expansion can also be utilized to sample these fields. The benefit of the circulant embedding technique is that it yields an exact representation of the random field on the sampling mesh for all $\Phi$, whereas the KL expansion gives a low-dimensional representation of the field introducing a bias (due to truncation after a finite number of eigenmodes). The number of terms in the KL expansion is dependent on $\Phi$ \cite{RF1,RF3} and can be large for stochastic processes with small correlation lengths.
\end{remark}

\begin{figure}[htb]
\begin{tabular}{cc}
\includegraphics[width =0.45\textwidth]{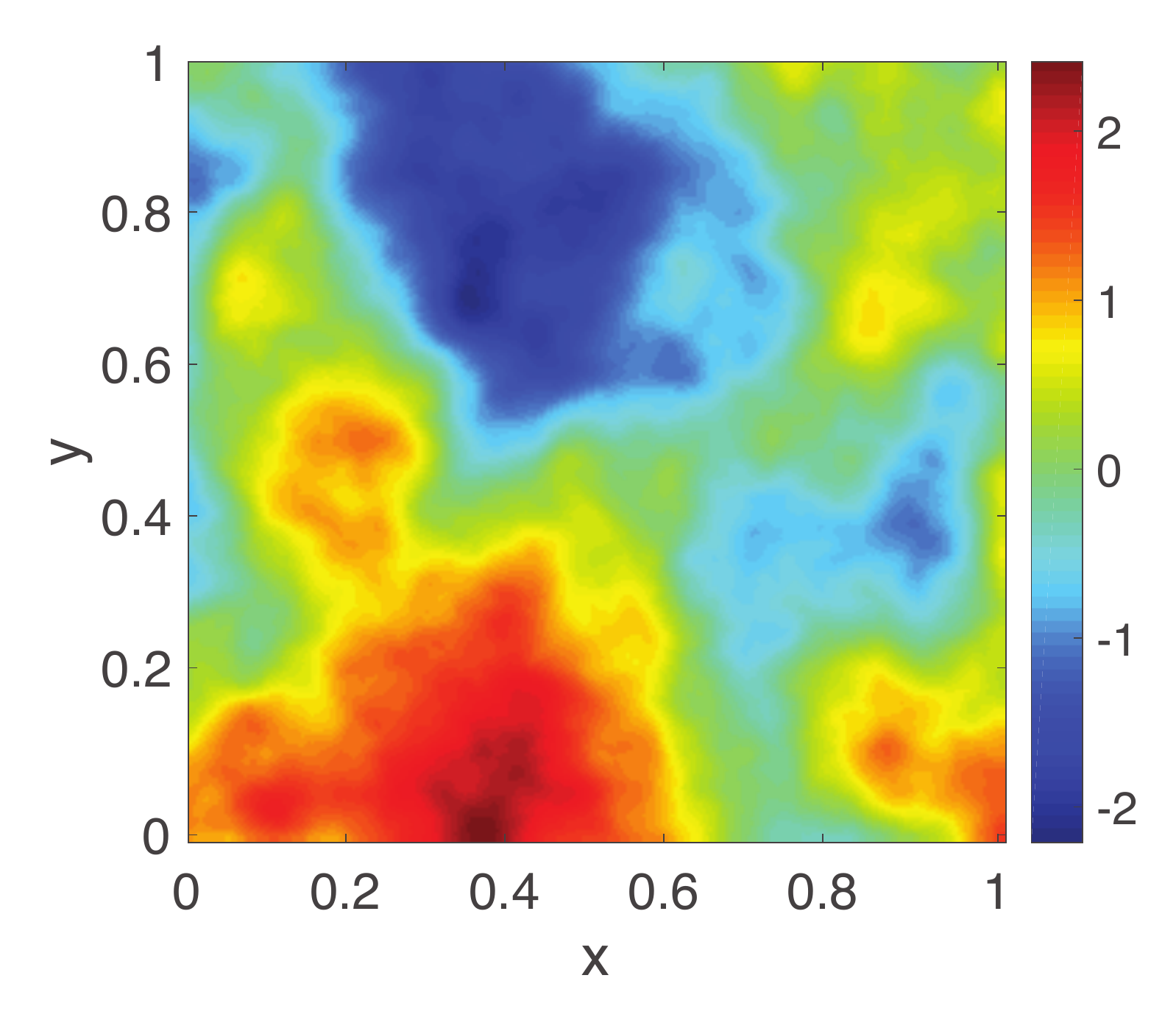} & 
\includegraphics[width =0.45\textwidth]{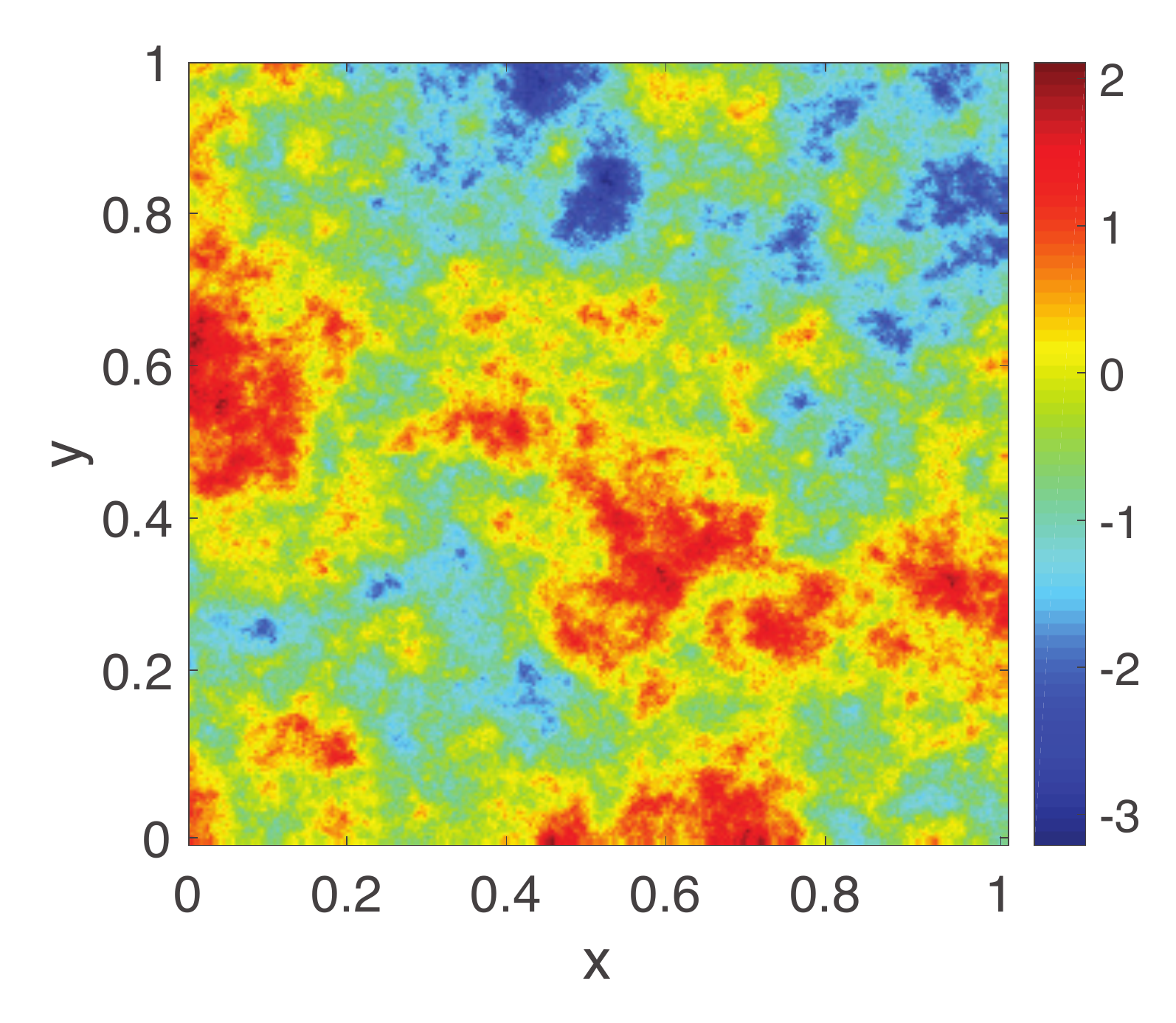} \\
(a) $\Phi_1 = (1.5,0.3,1)$ & (b) $\Phi_2 = (0.5,0.3,1)$ \\
\includegraphics[width =0.45\textwidth]{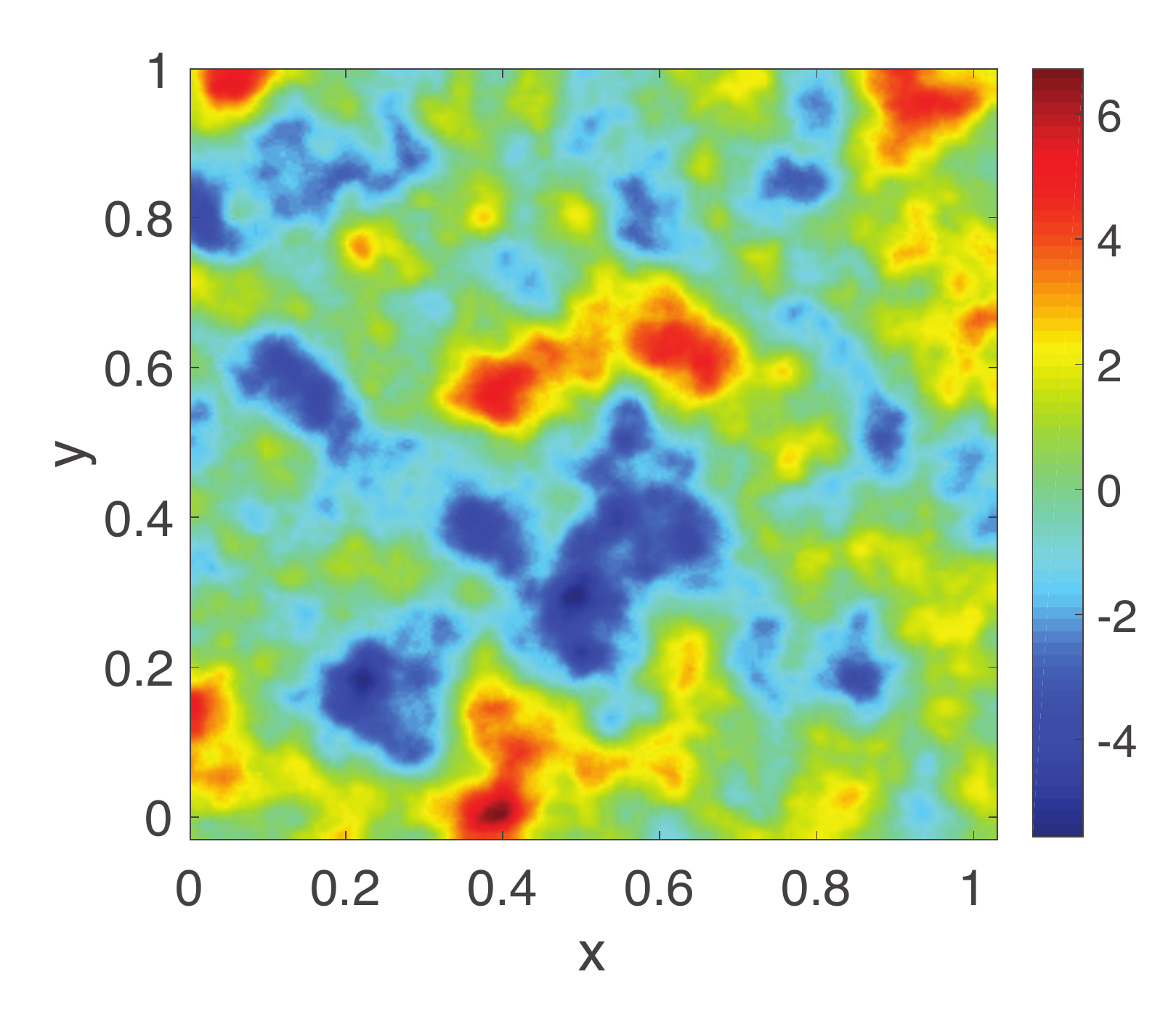} & 
\includegraphics[width =0.45\textwidth]{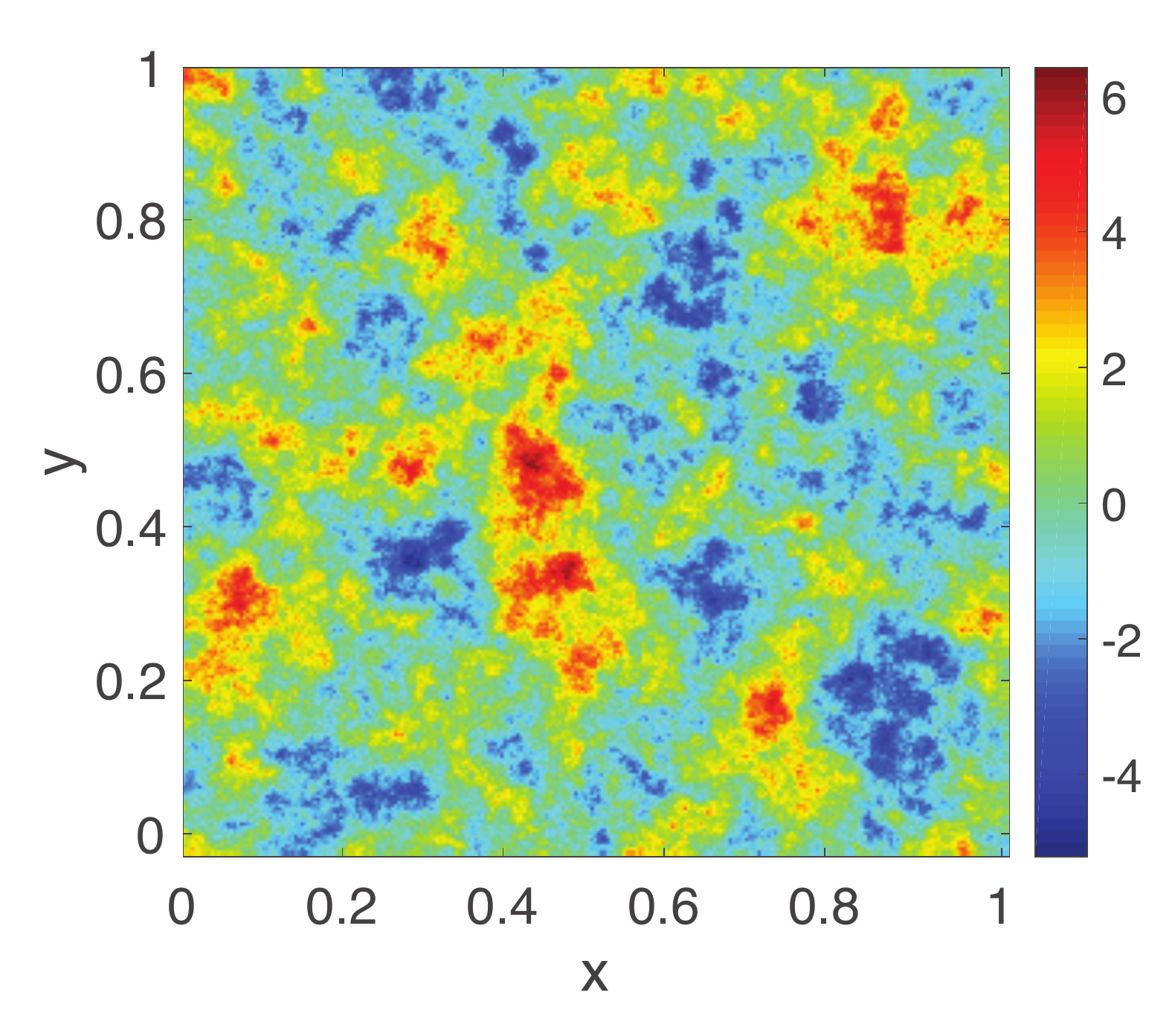} \\
(c) $\Phi_3 = (1.5,0.1,3)$ & (d) $\Phi_4 = (0.5,0.1,3)$ \\
\end{tabular}
\caption{Logarithm of the permeability field, $\log_{10} k$, generated using four reference parameter sets 
$\Phi = (\nu_c, \lambda_c, \sigma_c^2)$.}
\label{logperm_figs}
\end{figure}

To illustrate the performance of LFA for PDEs with random parameters, we consider four Mat\'ern reference parameter sets $\Phi$ with increasing order of complexity, listed in Table \ref{fourMat}. Random fields generated with these Mat\'ern parameter sets, sampled on a uniform mesh, are presented in Figure \ref{logperm_figs}.
\begin{table}[H]
\caption{Different combinations of the Mat\'ern reference parameters ${\Phi}=(\nu_c,\lambda_c,\sigma_c^2)$ with increasing complexity from left to right.}\label{fourMat}
\begin{center} 
\begin{tabular}{cccc}
\hline
${\Phi_1}$&${\Phi_2}$ & ${\Phi_3}$ & ${\Phi_4}$\\
\hline
(1.5,0.3,1)&(0.5,0.3,1)&(1.5,0.1,3)&(0.5,0.1,3)\\
\hline
\end{tabular}
\end{center}
\end{table}

For a fixed $\Phi$ and mesh size $h$, we generate $N_{MG}$ realizations of the permeability field. For each of these random fields, we can compare the LFA two-grid convergence factors with the computed asymptotic convergence factors of the multigrid method by using the same components that were utilized for randomly jumping coefficient fields. The comparison is shown in Table \ref{results}, where different numbers of smoothing steps are considered and the comparison is done for the four reference parameter sets, described in Table \ref{fourMat}. For all experiments, we set $h=1/64$ and $N_{MG} =N_{LFA}=100$. In general, we observe a good agreement between the experimental and the LFA quantities. For the first three sets of parameters excellent convergence factors are obtained already with two pre- and two post-smoothing steps. When the more difficult set of parameters is considered, however, more smoothing steps may be necessary to obtain a good convergence factor. It is also pointed out that the performance of W(1,1)-cycle is well predicted by LFA for all $\Phi_i, \; i=1,4$. A slight discrepancy is observed for the W(2,2)- and W(3,3)-cycles in the case of $\Phi_4$.
\begin{table}[htb]
\begin{center}
\begin{tabular}{|c|c|c|c|c|c|c|c|c|c|}
\cline{3-10}
\multicolumn{2}{c|}{} & \multicolumn{2}{|c|}{$\Phi_1$} & \multicolumn{2}{|c|}{$\Phi_2$} & \multicolumn{2}{|c|}{$\Phi_3$} & \multicolumn{2}{|c|}{$\Phi_4$}  \\
\cline{3-10}
\multicolumn{2}{c|}{} & MG & LFA & MG & LFA & MG & LFA & MG & LFA \\
\hline
\multirow{2}{*}{$W(1,1)$} & $\langle\rho\rangle$  &0.20  &0.20   & $0.20$ & $0.20$ &0.19 &   0.20& $0.23$ & $0.21$ \\
& $\sigma$ &0.004&0.002  & $0.004$ & $0.005$ &0.005   & 0.01 & $0.04$ & $0.02$ \\
\hline 
\multirow{2}{*}{$W(2,2)$} & $\langle\rho\rangle$ &0.04  &0.04   & $0.04$ & $0.04$ &0.07   & 0.06   & $0.17$ & $0.08$\\
& $\sigma$ &0.001  &0.001  & $0.002$ & $0.003$ &0.03 &0.03  & $0.04$ & $0.03$ \\
\hline
\multirow{2}{*}{$W(3,3)$} & $\langle\rho\rangle$ &0.01 &0.02   & $0.02$ & $0.02$ &0.05  & 0.03 & $0.13$ & $0.06$ \\
& $\sigma$ &0.001 & 0.002 & $0.002$ & $0.004$ &0.02   & 0.02  & $0.03$ & $0.03$\\
\hline
\end{tabular}
\end{center}
\caption{Comparison of the average and the standard deviation of the LFA and MG convergence factors. Different numbers of smoothing steps, $h = 1/64$, and the four reference parameter sets in Table  \ref{fourMat} are considered.}
\label{results}
\end{table}
\subsection{Mesh dependency}
\begin{figure}[hbt]
\begin{center}
{\includegraphics[clip, trim=2cm 12cm 2cm 2.5cm,scale=.9]{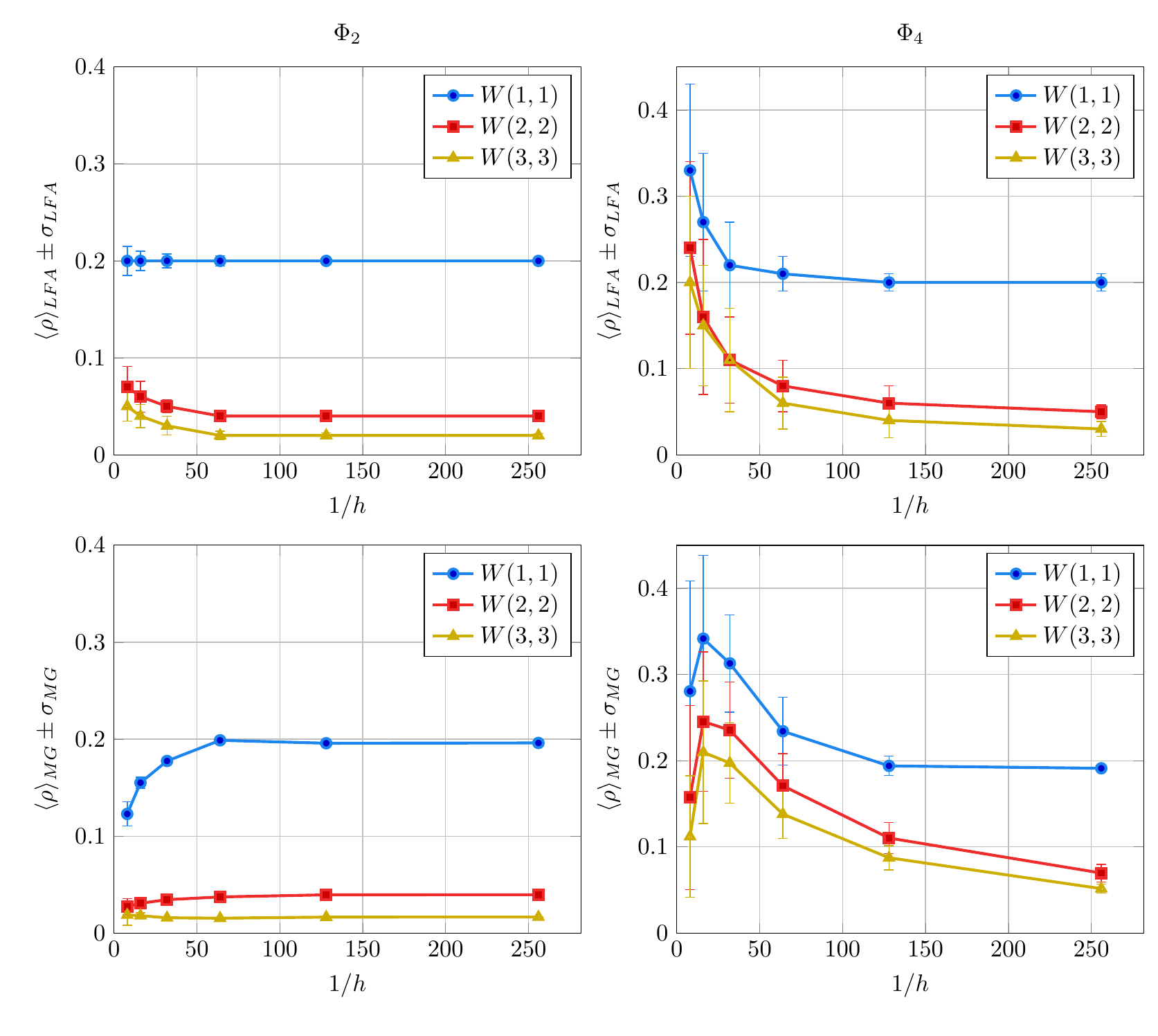}}
\end{center}
\caption{Average reduction LFA two-grid factors (top)  and asymptotic MG convergence (bottom) along with the standard deviation for two reference parameter sets  $\Phi_2$ (left column) and $\Phi_4$ (right column); $N_{MG}=N_{LFA}=100$. } 
\label{second_fourth_figs}
\end{figure}

Next, we wish to study the influence of the size of the sampling mesh on the multigrid convergence. For this purpose, we choose two representative parameter sets describing a smooth and a highly oscillating random fields, that are generated by the parameter sets $\Phi_2$ and $\Phi_4$, respectively. Figure \ref{second_fourth_figs} shows the average convergence factors for $\Phi_2$ (left side) and $\Phi_4$ (right side), predicted by LFA (top) and experimentally observed multigrid convergence (bottom), for different mesh sizes and different numbers of smoothing steps. For $\Phi_2$, the average reduction factor is roughly the same, independently of the size of the mesh. These predictions coincide well with the experimentally observed factors. For the results corresponding to parameter set $\Phi_4$, we observe robustness of the method when the mesh is sufficiently fine, which is also confirmed by the multigrid experiments. In the same figure, the standard deviation is presented, which decreases when $h\rightarrow 0$. The LFA predictions on the coarser grids are less reliable when compared to the multigrid convergence rates. We would like to mention that a three-grid LFA~\cite{three_grid} is not helpful here. 
%This can be attributed to the fact that while performing the LFA only a subset of the random field is utilized to derive the two-grid LFA operators which might not be fully representative of the properties of the random field.   
On coarse grids, boundary conditions have an impact on the method's convergence, but they are not taken into the account in the analysis. %Therefore, if an additional boundary relaxation is included in the multigrid algorithm, we can expect even better match.  
\begin{figure}[hbt]
\begin{center}
{\includegraphics[clip, trim=0cm 6cm 0cm 6cm,scale=0.38]{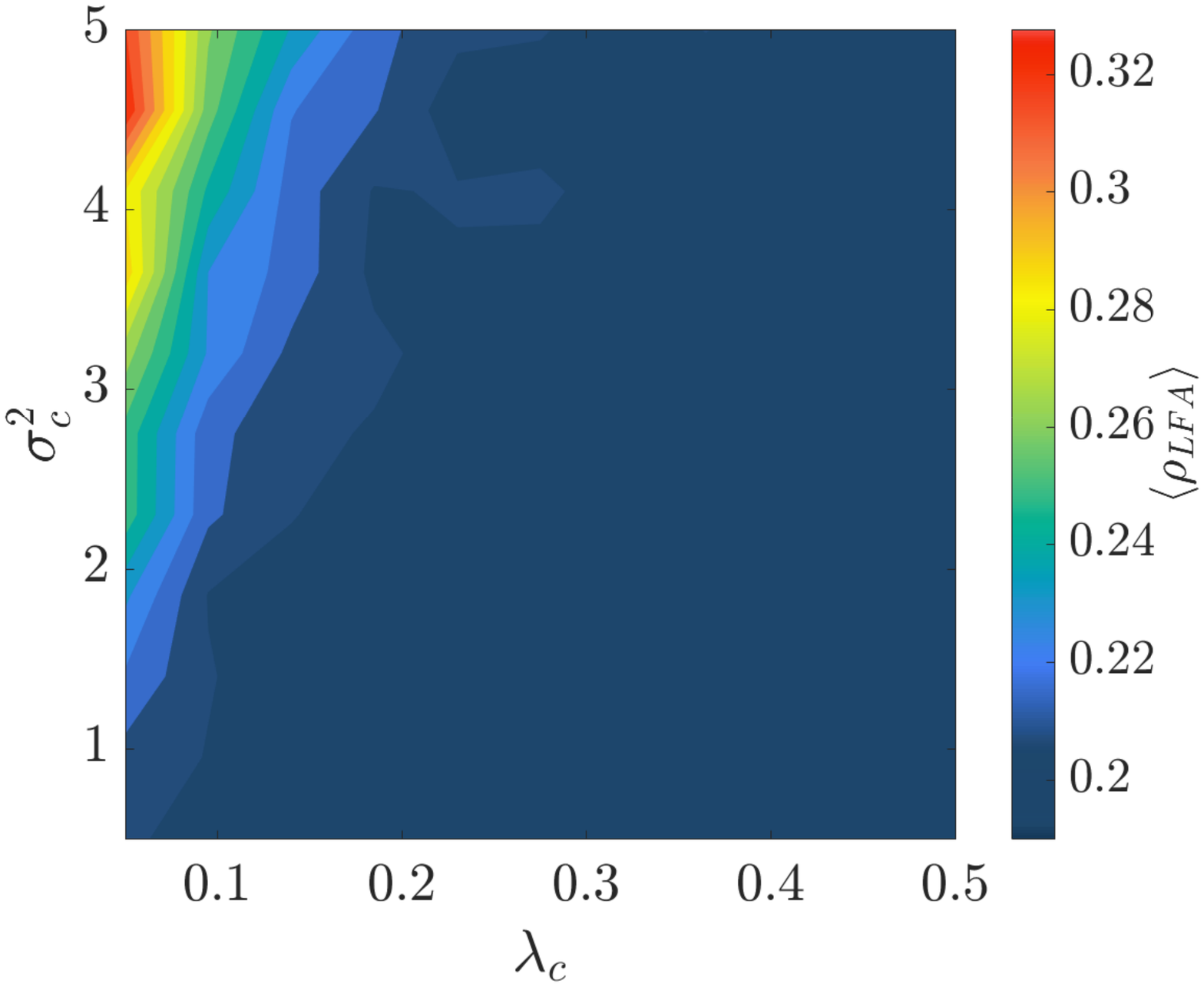}}
{\includegraphics[clip, trim=0cm 6cm 0cm 6cm,scale=0.38]{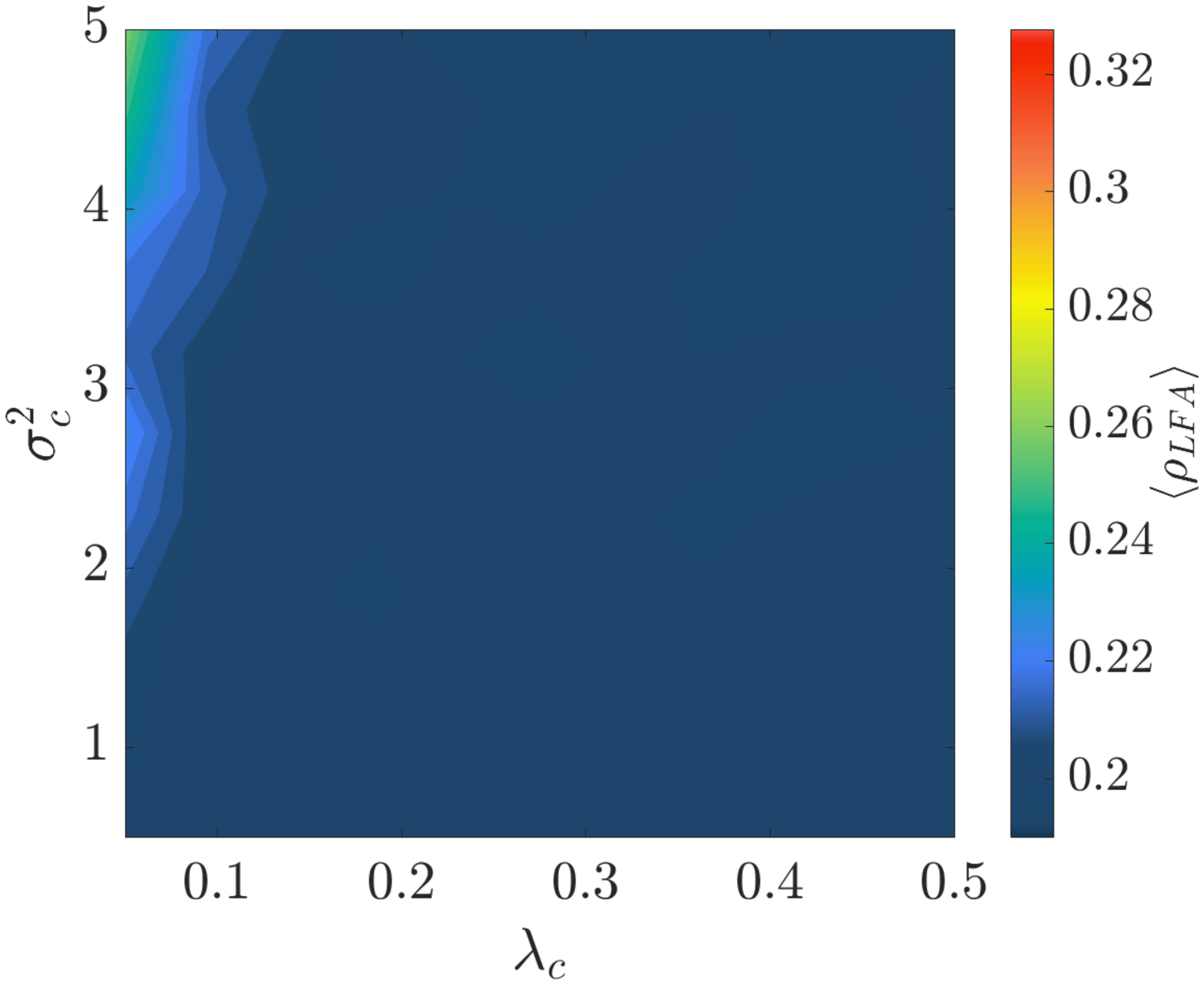}}
\end{center}
\caption{Contour of average LFA two-grid convergence factors, $\langle\rho\rangle_{LFA}$, with two smoothing steps for different $\lambda_c$ and $\sigma_c^2$, and for fixed $\nu_c=0.5$ with $h=1/32$ (left) and $h=1/64$ (right).}
\label{Heatmaps}
\end{figure}
Moreover, to analyze the robustness of the proposed multigrid method, in the next experiment we fix the field's smoothness parameter as $\nu_c=0.5$, and vary the other parameters $\sigma_c^2 \in[0.5,5]$ and $\lambda_c \in [0.05,0.5]$. In Figure~\ref{Heatmaps} we show the average LFA two-grid convergence factors when two smoothing steps are considered for $h=1/32$ (left) and $h=1/64$ (right).  The multigrid convergence is very satisfactory for all combinations of the parameters and increases slightly when $\lambda_c$ tends to be small and $\sigma_c^2$ becomes large. This case, however, represents a rather extreme situation in which the jumps in the permeability field are of more than $15$ orders of magnitude. Again, we see an improvement in the convergence rate with grid refinement. Note that with the considered range of parameters we cover all realistic cases. 

Here we have demonstrated the multigrid convergence for isotropic random field. In \emph{Appendix C} we also perform tests on anisotropic and non-grid aligned random fields. There we show that by using the same multigrid algorithm, convergence rates close to $0.2$ are achievable.
\section{Multilevel Monte Carlo computations}\label{sec:6}
In the preceding section, we have demonstrated that LFA performs well in predicting the mean convergence rate of the multigrid method for a given Mat\'ern parameter set. In this section, we wish to utilize LFA based information to select the multigrid cycle for MLMC simulations. %We begin by describing the single-level and multilevel variants of the MC method. When solving PDEs with high-dimensional input uncertainties, MC type methods are often preferred due to their dimension independent convergence properties and the convenience of implementation. For any quantity of interest $Q$, the expected value $\E[Q]$ is approximated by means of $N$ independent identically distributed (i.i.d.) realizations
Within MLMC methods, a {\em hierarchy of grids}, $\{D_{h_{\ell}}\}_{\ell=0}^L$, characterized by $h_0 > \ldots > h_L$ is defined to accurately estimate the quantity of interest $Q$, based on the linearity of the expectation operator, i.e.,
\begin{equation}
\mathbb{E}[Q_{h_L}] = \mathbb{E}[Q_{h_0}] + \sum_{\ell=1}^{L } \mathbb{E}[Q_{h_\ell}-Q_{h_{\ell-1}}].
\label{ml}
\end{equation}
The advantage of using the above telescopic decomposition is that on the coarsest grid, for $\ell = 0$, large numbers of samples to accurately determine $\mathbb{E}[Q_{h_0}]$ can be computed at a small computational effort. For larger values of $\ell$ the numerical solution is comparatively expensive, however, fewer samples are required on these finer levels, as the variance of the correction terms, $\mathbb{V}[Q_{h_\ell}-Q_{h_{\ell-1}}]$, is significantly smaller compared to the variance $\mathbb{V}[Q_{h_\ell}]$. A multilevel estimator for $\mathbb{E}[Q_{h_L}] $ can be defined in terms of the standard MC estimator as:
%\begin{align}
%\mathpzc{E}^{ML}_{L}[Q_{h_L}]:&=\mathpzc{E}^{MC}_{N _0}[Q_{h_0}] + \sum_{\ell=1}^L\mathpzc{E}^{MC}_{N _\ell}[Q_{h_\ell} - Q_{h_{\ell-1}}]\nonumber\\ 
%&=\frac{1}{N_0} \sum_{i=1}^{N_0} Q_{h_0}(\omega_i)+ \sum_{\ell=1}^L \left(\frac{1}{N_\ell} \sum_{i=1}^{N_\ell} (Q_{h_\ell}(\omega_i) -Q_{h_{\ell-1}}(\omega_i)) \right),
%\label{mlest}
%\end{align}
\begin{equation}
\mathpzc{E}^{ML}_{L}[Q_{h_L}]:=\frac{1}{N_0} \sum_{i=1}^{N_0} Q_{h_0}(\omega_i)+ \sum_{\ell=1}^L \left(\frac{1}{N_\ell} \sum_{i=1}^{N_\ell} (Q_{h_\ell}(\omega_i) -Q_{h_{\ell-1}}(\omega_i)) \right),
\label{mlest}
\end{equation}
where the level dependent samples $N_\ell \in \mathbb{N}$ form a decreasing sequence for increasing $\ell$-values.  We will use the same geometric sequence of meshes in the MLMC estimator as employed for the multigrid method described in Subsection \ref{MGalgo}, based on uniform coarsening such that $h_\ell = h_{\ell-1}/2 =2^{-\ell}h_0$, with $h_0$ the coarsest mesh size. As each of the estimators in \eqref{mlest} is independent, the variance of the multilevel estimator is the sum of the variances of individual estimators, i.e.
\begin{equation}
\mathbb{V}[\mathpzc{E}^{ML}_{L}[Q_{h_L}]] = \sum_{\ell=0}^L \frac{\mathcal{V}_\ell}{N_\ell},
\label{mlvar}
\end{equation}
where the level dependent sample variance, $\mathcal{V}_\ell=\mathbb{V}[Q_{h_\ell}-Q_{h_{\ell-1}}]= \mathcal{O}(2^{-\beta\ell}),\beta>0$ and $Q_{h_{-1}}=0$. Assuming that the mean FV error for mesh width $h_\ell$ satisfies
\begin{equation}
|\mathbb{E}[Q_{h_\ell}-Q]| = \mathcal{O}(h_\ell^\alpha),
\end{equation}
with $\alpha>0$. The mean square error (MSE) of the MLMC estimator can be expressed as a sum of bias and sampling errors, 
\begin{align}
\mathbb{E}[(\mathpzc{E}^{ML}_{L}[Q_{h_L}]-\mathbb{E}[Q])^2] &=(\mathbb{E}[Q_{h_L}-Q])^2+ \mathbb{V}[\mathpzc{E}^{ML}_{L}[Q_{h_L}]],\nonumber\\
&=c_0h_L^{2\alpha}+\sum_{\ell=0}^L \frac{\mathcal{V}_\ell}{N_\ell},
\label{mlvar2}
\end{align}
with $c_0$ a constant independent of $h_L$ and $\alpha$.  %As a direct measure of the FV bias is not available, the following identity can be utilized to determine the FV bias with respect to mesh size $h_L$,
%\begin{align}\label{Dbias}
%\left|\mathbb{E}[Q - Q_{h_L}]\right|  &\leq \frac{1}{2^\alpha-1}\left|\mathbb{E}[Q_{h_L} - Q_{h_{L-1}}]\right|,
%\end{align} 
%where $\left|\mathbb{E}[Q_{h_L} - Q_{h_{L-1}}]\right|$ can be approximated numerically. Moreover, for  a given user-specified tolerance $\varepsilon$, the above criteria can be used to determine the finest resolution in the MLMC hierarchy. 
%Similar to the standard MC method the sampling error is balanced with the discretization bias. For this, 
The level dependent sample  size $N_\ell$ can be chosen such that the sampling errors of each of the MC estimators in \eqref{mlest} reduces to the size of the discretization error. This is achieved by taking $N_\ell$ such that  $\mathcal{V}_\ell/N_\ell = \mathcal{O}(h_L^{2\alpha}) = \mathcal{O}(2^{-2\alpha L})$. Moreover, by fixing the number of the finest level samples $N_L$, the number of coarse grid samples can be simply expressed as
\begin{equation}\label{sampNL}
N_\ell = \lceil N_L2^{\beta(L - \ell)}\rceil.
\end{equation}
In practice the value $N_L$ is small $\sim \mathcal{O}(1)$ and can be chosen heuristically \cite{mishra2012sparse,Mishra20123365}. The sampling error on the coarsest level  does not depend on $\beta$, and therefore, by using \eqref{sampNL}, the balance ${\mathcal{V}_0}/N_0 = \mathcal{O}(h_L^{2\alpha})$ holds only asymptotically for large $L$.  Alternatively, one can also solve an optimization problem, as in \cite{MLMC2,MLMC1,ANU:9672986}, which minimizes the total cost of the MLMC estimator for a given tolerance $\varepsilon$. %The considered sampling approach requires a-priori values for the convergence rates $\alpha$ and $\beta$. The numerically observed convergence rates for $|\E[Q_{h_\ell} - Q_{h_{\ell-1}}]|$ and $\V[Q_{h_\ell} - Q_{h_{\ell-1}}]$ with respect to $h_\ell$ can be utilized for $\alpha$ and $\beta$, respectively.

%One of the advantages of this sampling strategy \eqref{sampNL} is that the number of samples on all levels is known and can be used to optimally distribute the cost on parallel computing systems. Following this, t
The computational cost of the MLMC estimator is given by
\begin{equation}\label{sum_cost}
\mathcal{W}^{ML}_L = \sum_{\ell=0}^{L}{N_\ell}\cdot({\mathcal{W}}_{\ell}+{\mathcal{W}}_{\ell-1}),
\end{equation}
where $\mathcal{W}_{\ell}$ represents total number of arithmetic operations to obtain a sample of the quantity of interest $Q$ on grid level $\ell$. This typically includes the cost of generating the random permeability field, assembling and solving the linear system of equations and post-processing. Out of these, the cost of solving the linear system is usually dominant, and, therefore, $\mathcal{W}_{\ell}$ will only represent the cost of the multigrid method. For the current combination of multigrid components, it is straightforward to show that the total amount of computational work for one MG cycle is proportional to the number of grid points on the finest MG level \cite{MG1}. The cost of the multigrid solver can thus  be expressed as:
\begin{equation}\label{MGcost}
\mathcal{W}_{\ell} =  \bar{c}\nu k_{\ell} h_{\ell}^{-2},
\end{equation}
where $\bar{c}$ is a constant, $k_{\ell}\in \mathbb{N}$ is the number of multigrid iterations, $\nu$ is the number of smoothing steps and $h_{\ell}^{-2}$ represents the number of unknowns on grid level $\ell$.  The cost associated with the other components, such as computations of the defect, transfer to coarser grids, interpolation of corrections, is accounted by the constant $\bar{c}$, to facilitate the comparison of W-cycles with different numbers of smoothing steps.

The product ${N}_\ell\cdot({\mathcal{W}}_{\ell}+{\mathcal{W}}_{\ell-1})$ in \eqref{sum_cost} gives the individual cost contributions at any level $\ell$. For the proposed multigrid solver, we have $( {\mathcal{W}}_{\ell}+ {\mathcal{W}}_{\ell-1}) =\mathcal{O}(h_\ell^{-2})$. For simplicity, let's assume that the cost of one sample on level  $\ell$ is given by $\mathcal{O}(h_\ell^{-\gamma})$, $\gamma\geq d$, where $d$ is the number of spatial dimensions. The computational cost \eqref{sum_cost} can be conveniently expressed as $\mathcal{W}^{ML}_L =\mathcal{O}\left(\sum^L_{\ell=0}2^{(\gamma-\beta)\ell}\right)$, leading to the following three cases. When the rate of decay of $N_\ell$ is faster than the growth in computational cost ($\beta> \gamma$), the dominating cost comes from the coarsest level. When the sample decay rate is similar to the growth in computational cost ($\beta = \gamma$), all levels equally contribute to the total cost. If the sample decay rate is slower than the growth in computational cost ($\beta<\gamma$), the dominant cost is on the finest level. This scenario is typically encountered in many practical settings involving 3D problems with low regularity permeability models, resulting in a slow variance decay (or a small $\beta$).

%Alternatively, one can also obtain $N_\ell$ by solving an optimization problem, as in \cite{MLMC2,MLMC1,ANU:9672986}, which minimizes the total computational cost $\mathcal{W}_L^{ML}$ such that $\V[\mathpzc{E}^{ML}_{L}[Q_{h_L}]] < \varepsilon^2$. This strategy yields the number of samples $N_\ell \propto \sqrt{V_\ell/\mathcal{W}_\ell}$. However, for both sampling strategies, it can be shown that asymptotically the cost of the MLMC estimator to achieve a RMSE of $\varepsilon$ scales like $${\mathcal{W}}^{ML}_L(\varepsilon) =\mathcal{O}\Big(\varepsilon^{-2-max\big(0,\tfrac{\gamma-\beta}{\alpha}\big)}\Big),$$ see \cite{MLMC2} for the proof.

Note that in many cases the rates $\alpha,\beta,\gamma$ are not available. A common practice is to use a few ``warm-up samples and levels'' to obtain an estimate for these quantities. Having a reliable LFA tool can be further helpful in efficiently optimizing the MLMC simulation, especially to assess the performance of a multigrid solver for different combinations of transfer operators and smoothers. This way, one can select the most suitable combination for the considered random input data. We can also use the LFA to obtain solver statistics on fine grids where computing warm-up samples is not feasible due to high computational cost.

\subsection{MLMC numerical experiments}
\begin{figure}[H]
\begin{center}
\includegraphics[clip, trim=0cm 0cm 1cm 0cm,scale=0.24]{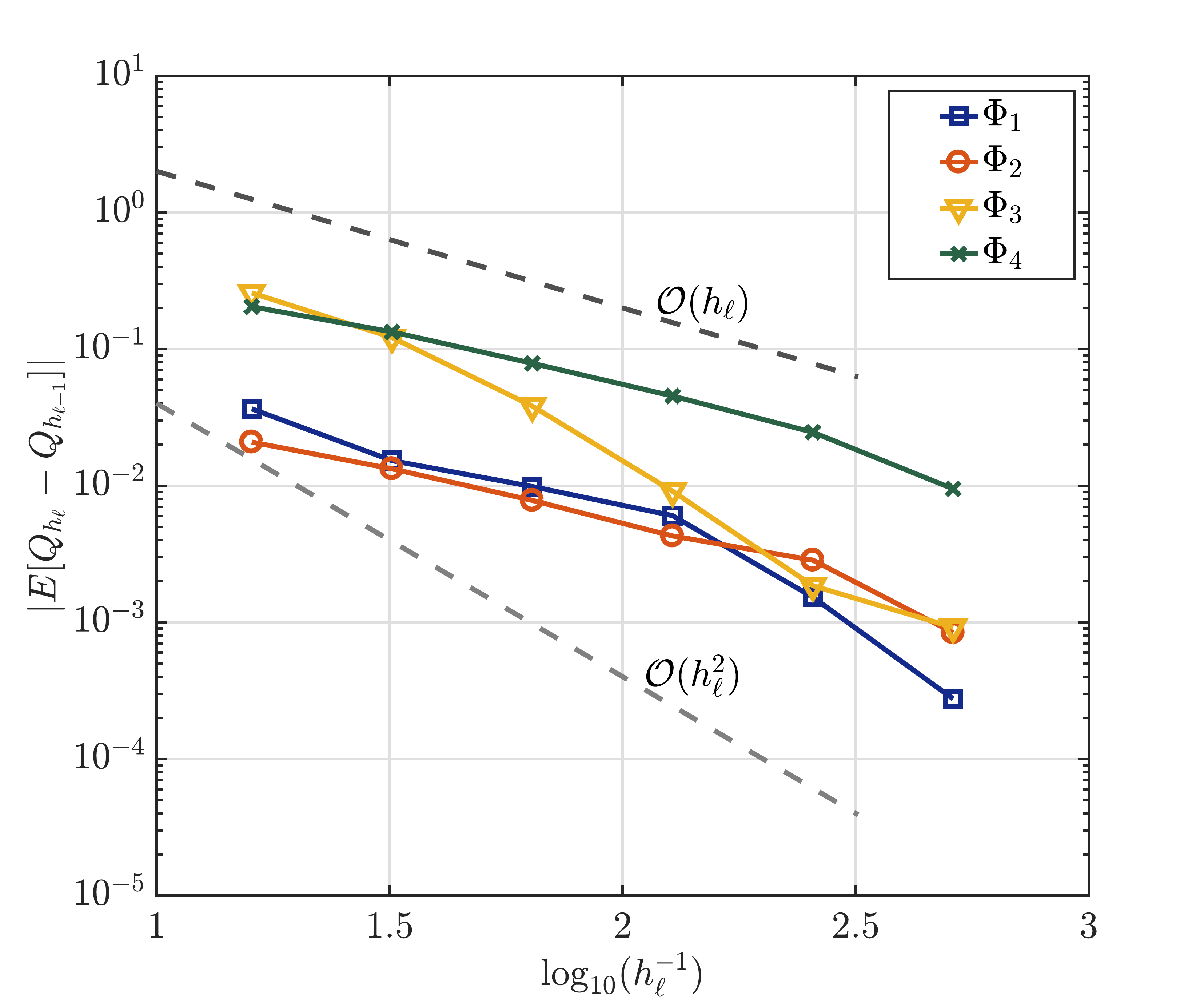}
\includegraphics[clip, trim=0cm 0cm 1cm 0cm,scale=0.24]{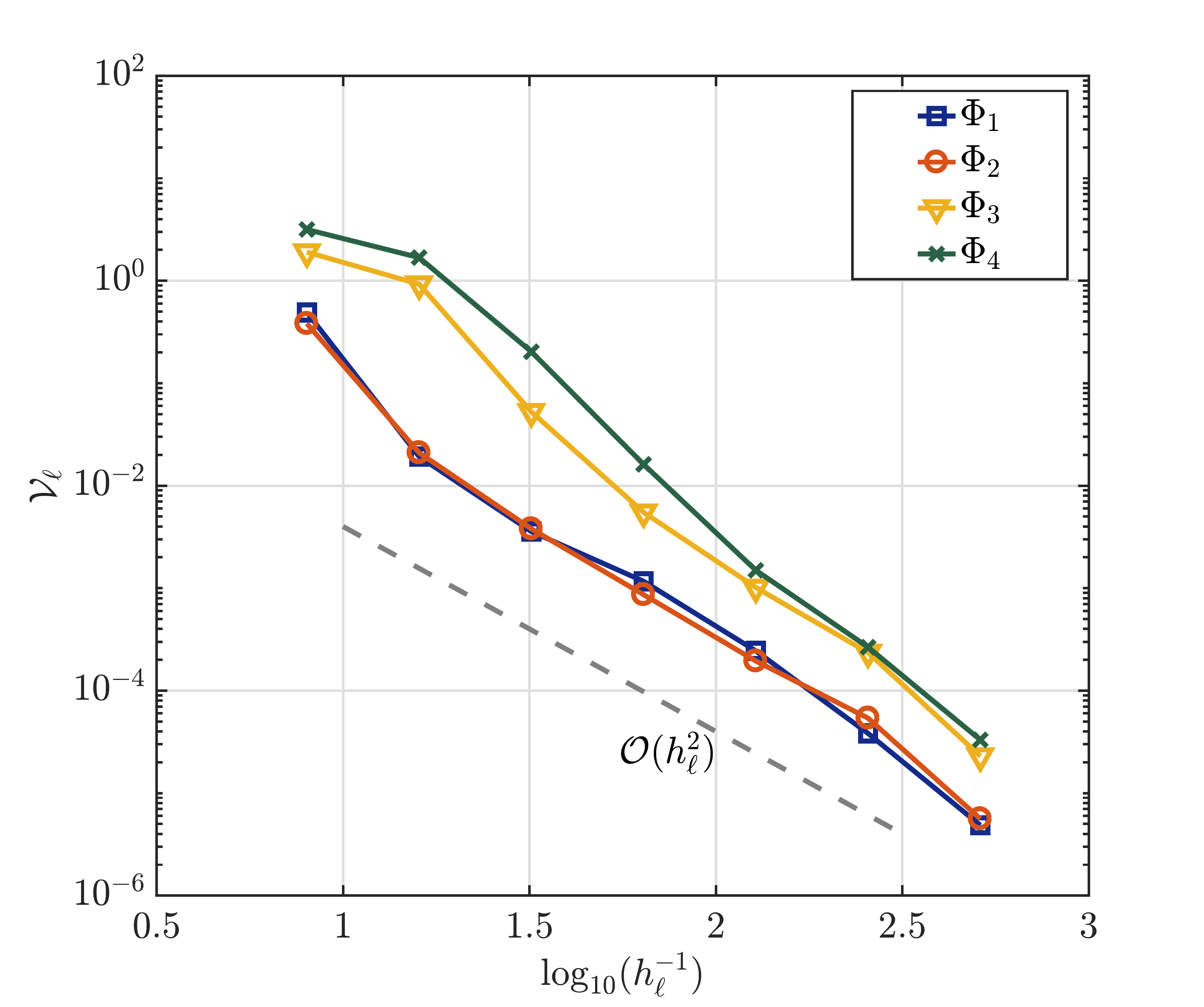}
\end{center}
\caption{(Left) Convergence of finite volume error, (Right) decay of level dependent variance for the outflow, $Q_{h_\ell}$, with mesh refinement for different Mat\'ern parameters.}
\label{alpha_beta}
\end{figure}

In this section we analyze the performance of the proposed multigrid MLMC method. 
%We will primarily focus on bench-marking different cycling strategies and reliability of the LFA prediction of MLMC cost for a given Mat\'ern parameter. 
We consider PDE \eqref{1a} on domain ${D}\in(0,1)^2$ with mixed Dirichlet-Neumann boundary conditions,
\begin{eqnarray}\label{BCdarcy}
u(0,y,\omega) =1, & \qquad u(1,y,\omega) = 0, \quad \text{and}\\ \nonumber
k(\mathbf{x},\omega)\frac{\partial u}{\partial x}\bigg|_{y=0} = 0, & \qquad k(\mathbf{x},\omega)\frac{\partial u}{\partial x}\bigg|_{y=1}  = 0, 
\end{eqnarray}
respectively. For all tests the source term is set to zero, i.e. $f = 0$. As the quantity of interest, the outflow through the boundary  $x=1,y=(0,1)$, also referred to as {\em effective permeability}, is considered:
\begin{equation}\label{outflow}
Q(u)\:= -\int_0^1 k(\mathbf{x},\omega)\frac{\partial u}{\partial x}(\mathbf{x},\omega)\bigg|_{x=1} dy.
\end{equation}
In Figure \ref{alpha_beta}, we show the convergence of the FV bias, $|\mathbb{E}[Q_{h_\ell} - Q_{h_{\ell-1}}]|$, and the level dependent variance $\mathcal{V}_\ell$ for the considered quantity of interest \eqref{outflow}. We observe that the rate of decay of the FV bias depends on the smoothness parameter $\nu_c$ of the random fields. In the case of parameter sets $\Phi_1$ and $\Phi_3$, we see a second-order convergence whereas a first-order convergence is observed for $\Phi_2$ and $\Phi_4$. The correlation length $\lambda_c$ and variance $\sigma^2_c$ only affect the proportionality constant. Also, for all four cases, the level dependent variance shows a quadratic decay.
To compute samples of $Q_{h_\ell}(\omega_i) - Q_{h_{\ell-1}}(\omega_i)$, we invoke the multigrid solver twice on level $\ell$ and $\ell-1$ with the permeability field ${k}(\mathbf{x}_{\ell},\omega_i)$ and an upscaled version of permeability ${k}(\mathbf{x}_{\ell-1},\omega_i)$, respectively. We define ${k}(\mathbf{x}_{\ell-1},\omega_i)$ using the \emph{covariance upscaling} technique from \cite{mishra2016multi} (see \emph{Appendix B}). It is important to emphasize that homogenization techniques, e.g. \cite{moulton1,scott} for obtaining the coarse-scale solution, cannot be used to compute $Q_{h_{\ell-1}}(\omega_i)$. This is due to the fact that the homogenized permeabilities are typically obtained by averaging procedures that can modify the spatial covariance as well as the marginal distribution defined for the permeability field. This will lead to a violation of the telescopic identity \eqref{ml}. However, homogenization methods can still be used within a multigrid solver to obtain coarse grid discretization operators.
\begin{remark}
The experiments are performed on a parallel machine with a homogeneous set of cores. The MLMC algorithm allows for three degrees of parallelization: over samples, over Monte Carlo levels and within the multigrid solver. Recent advancements on parallelization can be found in \cite{Mishra20123365,doi:10.1137/16M1083591,10.1007/978-3-642-31464-3_25}. We only consider parallelization over the MLMC samples and levels. For the considered MLMC sampling strategy with known number of samples and levels, a greedy algorithm based approach can be followed where cores are distributed between levels such that the sum of idle times for all cores is minimized. 
\end{remark}
Next, we conduct MLMC experiments for Mat\'ern parameter sets $\Phi_2$ and $\Phi_4$, using W(1,1), W(2,2) and W(3,3)-cycles. The multigrid components described in Section \ref{MGalgo} are used with stopping criterion $\varepsilon_{MG} = 10^{-10}$.  We consider $\alpha = 1$ and $\beta = 2$ for both parameter sets $\Phi_2$ and $\Phi_4$. Recall that the MLMC sampling strategy \eqref{sampNL} requires the total number of MLMC levels, the number of finest level samples $N_L$ along with the rate $\beta$. We set $N_L=64$ within MLMC for both parameter sets resulting in a sample sequence $N_\ell = N_L4^{(L-\ell)}$.  This sampling approach is simple but is sensitive to the choice of $N_L$, and a small value of $N_L$ may result in unconverged mean and variance on coarser levels.
\subsubsection{MLMC results for $\Phi_2$} 
\begin{figure}[hbt]
\begin{center}
{\includegraphics[clip, trim=1.8cm 19.5cm 9.8cm 2.4cm,scale=.9]{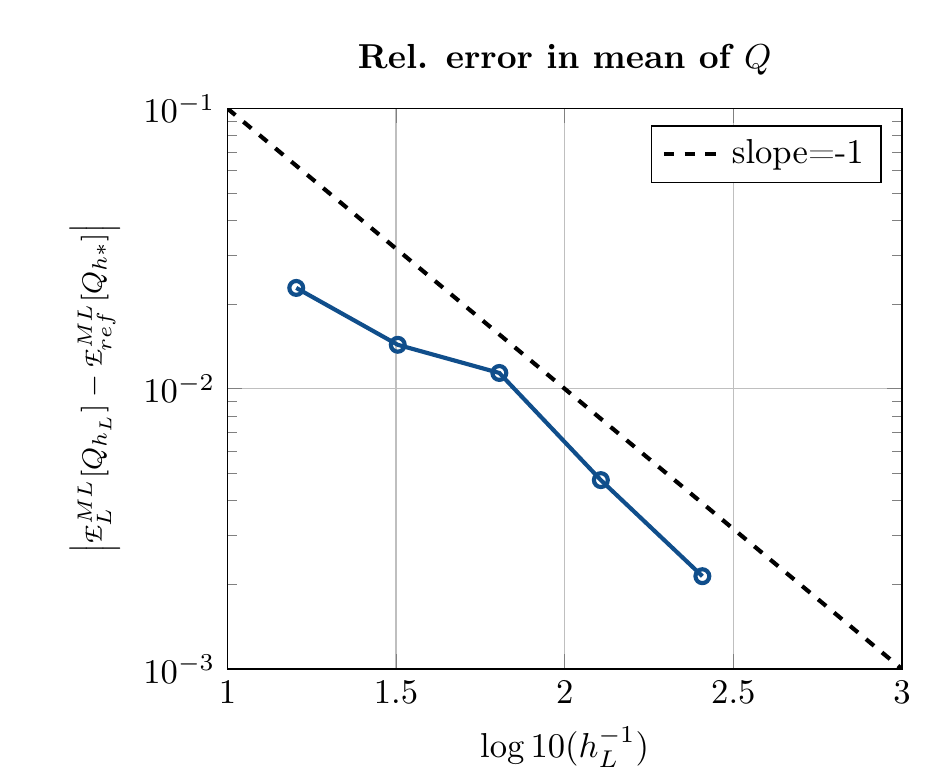}}
{\includegraphics[clip,trim=6.5cm 19.5cm 6.2cm 2.4cm,scale=.9]{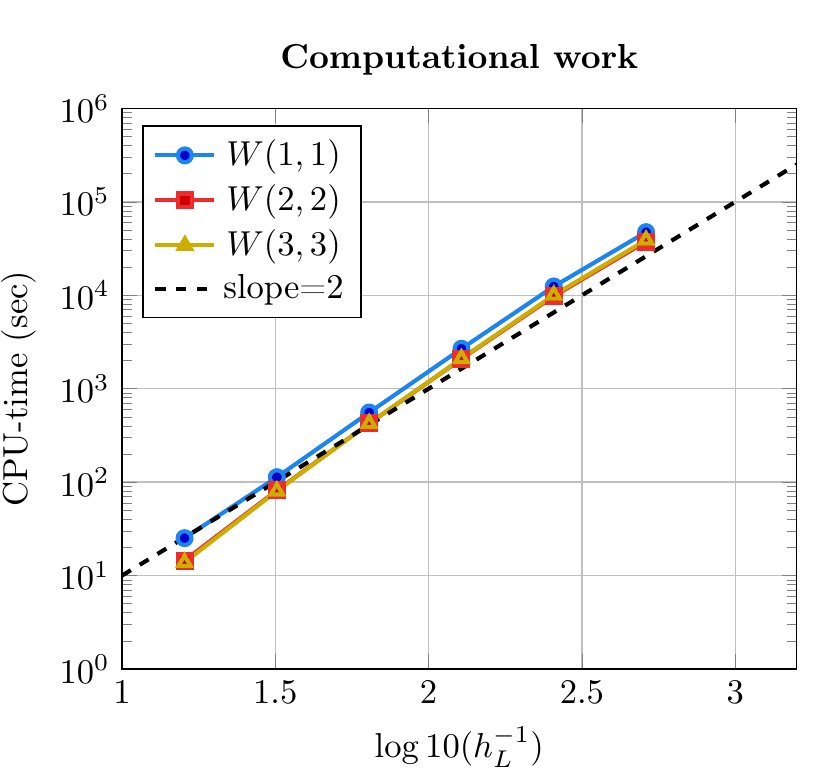}}
\end{center}
\caption{(Left) Convergence of error in $\mathpzc{E}_L^{ML}[Q_{h_L}]$ with increasing number of MLMC levels for parameter $\Phi_2$. The reference solution is based on mesh size $h^* = 1/512$. (Right) Mean CPU-times versus accuracy for different W-cycles. Computational cost proportional to $\mathcal{O}(h_L^{-2})$ is observed for all cycling strategies.}
\label{Rel_err}
\end{figure}

We first analyze the error convergence of the expected value and the variance of the quantity of interest. In Figure \ref{Rel_err} (left), relative errors with respect to a reference solution  computed using the MLMC estimator with the finest level $h^* =1/512$ are shown. We see the error $|\mathpzc{E}^{ML}_{L}[Q_{h_L}] -\mathpzc{E}^{ML}_{ref}[Q_{h^*}]|$ converges with $\mathcal{O}(h_L)$. Next, we analyze the performance of the MLMC estimator for different multigrid cycling strategies. In Figure \ref{Rel_err} (right), we present the CPU-times for the MLMC simulations with an increasing number of MLMC levels. These CPU-times are derived by summing up the run-times from the multigrid solves for all samples over all levels. Here, we do not include the cost of generation of the random field and the post-processing costs as they are the same for different cycles. The cost scales as $\mathcal{O}(h_L^{-2})$ coinciding with the theoretical MLMC complexity when $\beta = \gamma$. No substantial difference in the CPU-times from the three cycles is observed. We also show the level-wise CPU-times in Table \ref{MLMC_phi2} for the three cycles with the same number of MLMC samples per level. In general, the computational cost improves with levels for all cycles. We observe that the W(2,2)-cycle slightly outperforms the other two variants. 

\begin{figure}[hbt]
\begin{center}
{\includegraphics[clip, trim=1.5cm 22.2cm 1cm 2.7cm,scale=.9]{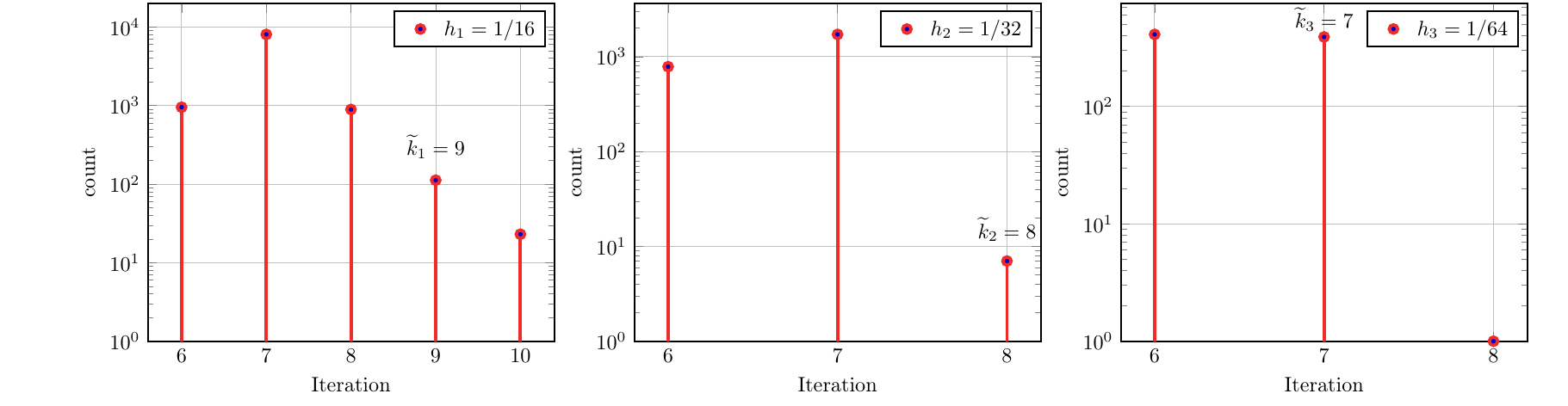}}
\end{center}
\caption{Distribution of number of multigrid iterations of $W(2,2)$-cycles for different mesh sizes for parameter $\Phi_2$ to reach $\varepsilon_{MG}<10^{-10}$ and the LFA predicted number of iteration $\widetilde{k}_\ell$.} 
\label{histo_phi2}
\end{figure}

We also show the distribution of the number of multigrid iterations for the W(2,2)-cycle on different grid levels along with the number of iterations predicted by the LFA, $\tilde{k}_\ell := \lceil \log(\varepsilon_{MG})/\log (\langle \rho \rangle_{LFA,\ell}) \rceil$ in Figure \ref{histo_phi2}. The cost is more heterogeneous on coarser levels with larger variance in  the number of iterations and is homogeneous for $h_\ell=1/32$ and onwards. A high variance in the multigrid convergence rate was also predicted by the LFA experiments in the previous section. Thus, a sparse direct solver can replace the multigrid solver on these coarser levels. In this work, the considered multigrid stopping criteria of $10^{-10}$ is quite conservative. For many engineering applications a residual reduction of $10^{-6}$ may already be sufficient to reach a converged solution, therefore, reducing the computational cost roughly by a factor of two. 

\begin{table}[hbt]
\begin{center}
\begin{tabular}{|c|c|c|c|c|}
\cline{3-5}
\multicolumn{1}{c}{} &\multicolumn{1}{c}{} & \multicolumn{3}{|c|}{level-wise CPU-time (sec.)}\\
\cline{3-5}
\hline
$h_\ell$ &$N_\ell$& W(1,1) &W(2,2) & W(3,3) \\
\hline
1/8     &262144   &489  &331 &323  \\
1/16   &65536     &528  &391 &390  \\
1/32   &16384     &454  &362 &363  \\
1/64   &4096       &407  &336 &348  \\
1/128 &1024       &376  &320 &342  \\
1/256 &256         &357  &308 &338  \\
1/512 &64           &350  &300 &333  \\
\hline
\end{tabular}
\end{center}
\vspace{0.2cm}
\caption{Comparison of the three W-cycles in terms of the level-wise CPU-times for parameter $\Phi_2$.}\label{MLMC_phi2}
\end{table}     
\subsubsection{MLMC results for $\Phi_4$} 
\begin{figure}[hbt]
\begin{center}
{\includegraphics[clip, trim=2cm 19.4cm 10.3cm 2.5cm,scale=.9]{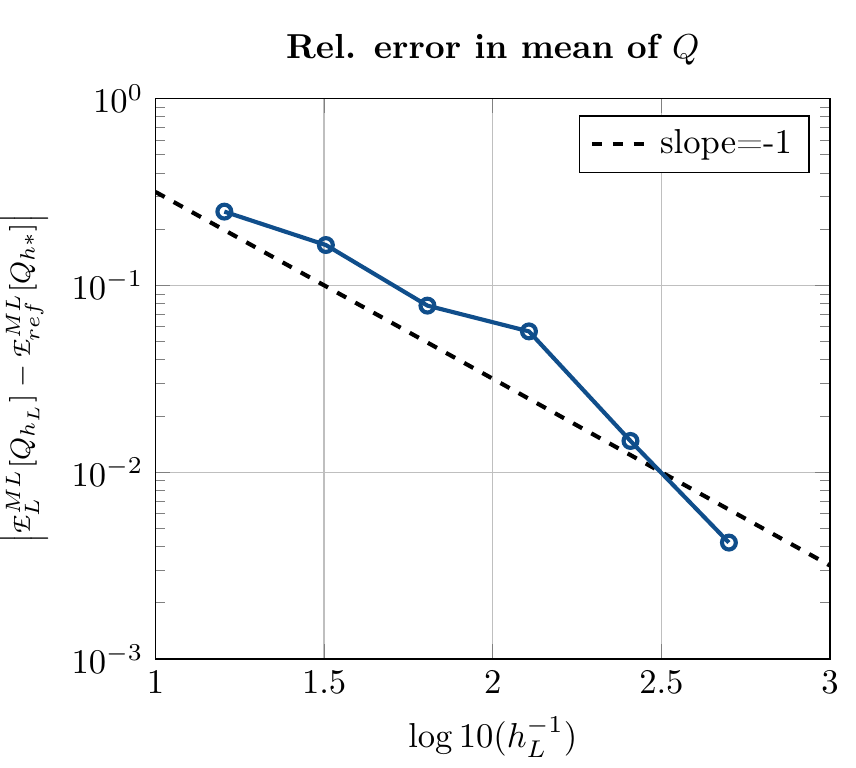}}
{\includegraphics[clip,trim=6.5cm 19.4cm 6.2cm 2.5cm,scale=.9]{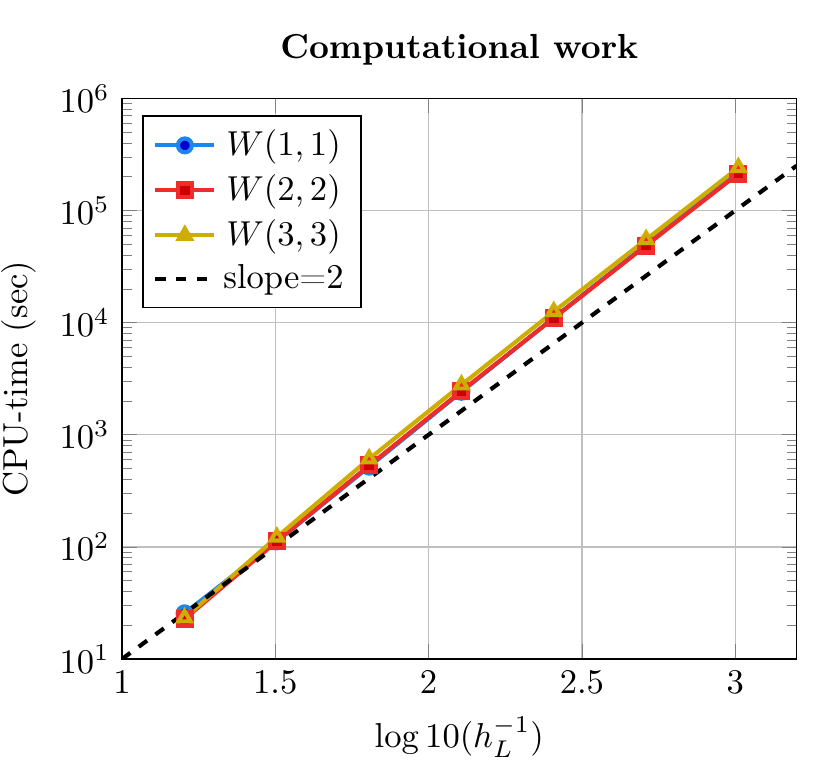}}
\end{center}
\caption{ (Left) Convergence of error in $\mathpzc{E}_L^{ML}[Q_{h_L}]$ with increasing number of MLMC levels for parameter $\Phi_4$. The reference solution is based on mesh size $h^* = 1/1024$. (Right) Mean CPU-times versus accuracy for different W-cycles. Computational cost proportional to $\mathcal{O}(h_L^{-2})$ is observed for all cycling strategies.}
\label{Rel_err_phi4}
\end{figure}

Now we describe MLMC results for the challenging parameter set $\Phi_4$ for the same quantity of interest. In Figure \ref{Rel_err_phi4} (left), we display the convergence of relative error in the expected value of the quantity of interest using a reference solution obtained with finest resolution $h^*=1/1024$. Due to a large variance $\sigma_c$ and a small correlation length $\lambda_c$ of the random field, the error is larger, as compared to the results of parameter set $\Phi_2$ for same mesh size $h_L$ and the \emph{expected} convergence rates are visible on relatively fine grids. A similar trend is observed for the computational cost with scaling as $\mathcal{O}(h_L^{-2})$ in Figure \ref{Rel_err_phi4} (right). Again, in Table \ref{MLMC_phi4}, the level-wise CPU-times for the three cycles are listed. Similar to the results for set $\Phi_2$, the cost improves with levels but the gain from W(2,2)-cycles is more prominent. 

\begin{figure}[hbt]
\begin{center}
{\includegraphics[clip, trim=1.5cm 17cm 1cm 2.7cm,scale=.87]{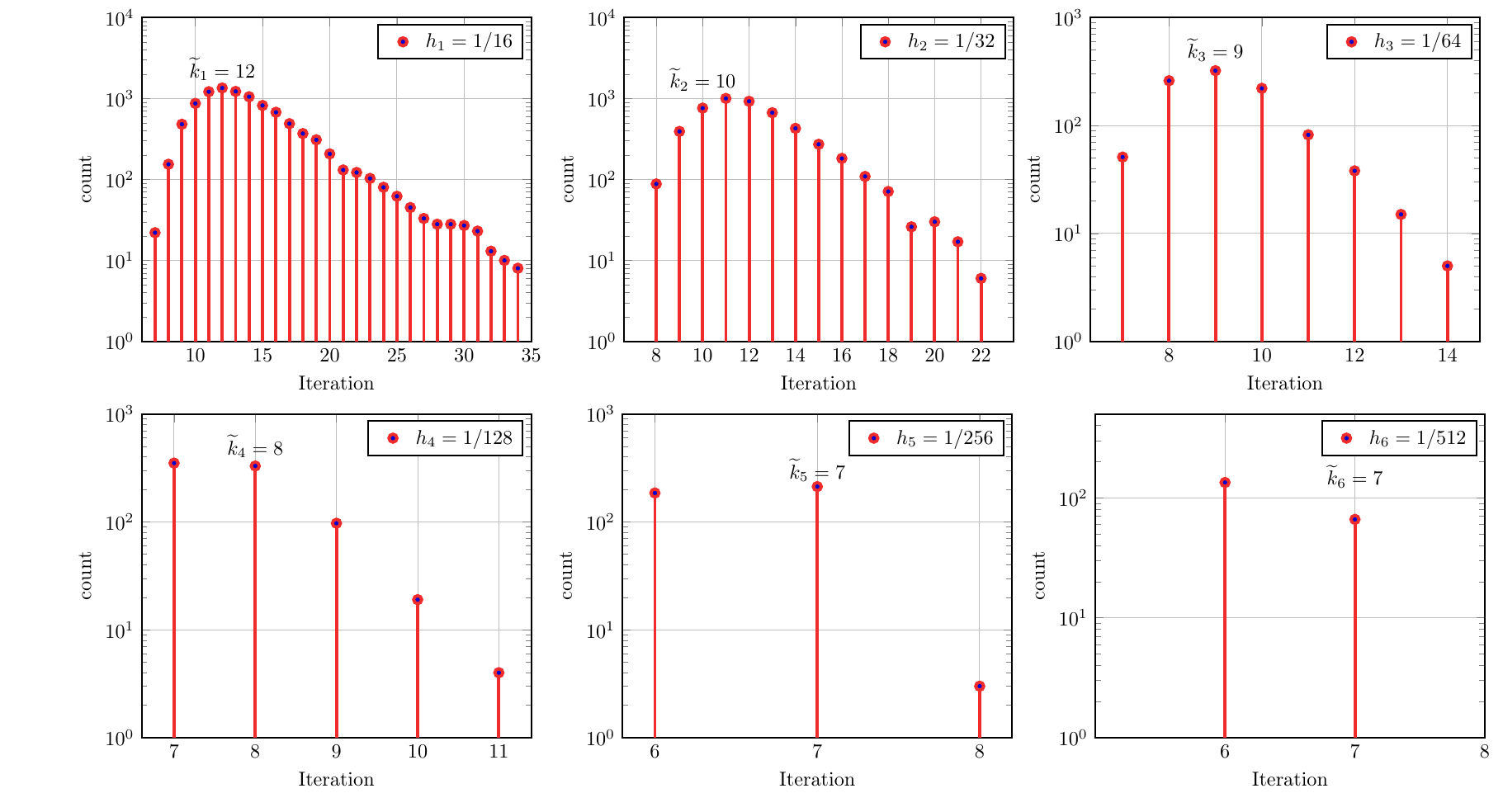}}
\end{center}
\caption{Distribution of number of multigrid iterations of $W(2,2)$-cycles for different mesh sizes for parameter $\Phi_4$ to reach $\varepsilon_{MG}<10^{-10}$ and the LFA predicted number of iteration $\widetilde{k}_\ell$.} 
\label{histo_phi4}
\end{figure}

Lastly, the distribution of the number of iterations for W(2,2)-cycles is depicted in Figure \ref{histo_phi4}.  High variability in the number of multigrid iterations persists until a relatively fine grid ($h=1/128$).  However, the average number of iterations predicted by LFA coincides very well with the mode of these distributions. Again, a direct solver can be applied on these coarser levels to get a more reliable estimate of the cost. Here, we point out that the W(2,2)-cycle takes roughly 7 iterations similar to parameter $\Phi_2$ to reduce the residual by $10^{-10}$, indicating robustness with respect to the Mat\'ern parameter.
\begin{table}[hbt]
\begin{center}
\begin{tabular}{|c|c|c|c|c|}
\cline{3-5}
\multicolumn{1}{c}{} &\multicolumn{1}{c}{} & \multicolumn{3}{|c|}{level-wise CPU-time (sec.)}\\
\cline{3-5}
\hline
$h_\ell$ &$N_\ell$& W(1,1) &W(2,2) & W(3,3) \\
\hline
1/8       &1048576 &2.23[+3]  & 2.10[+3] &2.18[+2]  \\
1/16     &262144   &2.42[+3]  & 2.52[+3] &2.80[+3]  \\
1/32     &65536     &1.85[+3]  & 2.10[+3] &2.45[+3]  \\
1/64     &16384     &1.56[+3]  & 1.59[+3] &1.98[+3]  \\
1/128   &4096       &1.45[+3]  & 1.33[+3] &1.60[+3]  \\
1/256   &1024       &1.40[+3]  & 1.23[+3] &1.40[+3]  \\
1/512   &256         &1.38[+3]  & 1.23[+3] &1.34[+3]  \\
1/1024 &64           &1.36[+3]  & 1.21[+3] &1.31[+3]  \\
\hline
\end{tabular}
\end{center}
\caption{Comparison of the three W-cycles in terms of the level-wise CPU-times for parameter $\Phi_4$.}\label{MLMC_phi4}
\end{table}                  
 \section{Conclusions}\label{sec:7}
A novel, generalized Local Fourier analysis which can be employed for quantitative assessment of multigrid methods for PDEs involving jumping and random coefficients has been presented. In particular, a cell-centered multigrid algorithm for solving a model problem based on a single-phase flow problem in a random porous medium has been used to show the accuracy of the proposed analysis. This technique, however, is appropriate for the prediction of the performance of a wider class of multigrid methods for solving PDEs with random fields. The effectiveness of the proposed LFA method was also confirmed numerically using a number of challenging test cases with jumping and random coefficients. 
%However, we would like to mention that the LFA will also help when other vertex-centered multigrid methods, like black-box MG, is utilized.  
%We also employed a simple multigrid algorithm based on a cell-centered discretization for solving problems elliptic PDEs with random fields. The simplicity of the approach is achieved by using a standard point-wise Gauss-Seidel smoother together with a piecewise constant prolongation operator and its adjoint as the restriction, plus a direct discretization to construct the operators on coarse grids. The fine convergence of CCMG was confirmed using the aforementioned generalized LFA technique. Further, we demonstrated that the prediction capability of the LFA method improves with grid refinement. 
Further, the novel local Fourier analysis helps us to estimate a-priori the time needed for solving certain uncertainty quantification  problems by using a multigrid multilevel Monte Carlo method. 
%The proposed multigrid solver was utilized for MLMC simulations for two Mat\'ern parameter sets. The optimal MLMC complexity of $\mathcal{O}(h_L^{-2})$ was achieved by all three variants of the W-cycle, the W(2,2)-cycle being the most efficient one. 
%%%%%%%%%%%%%%%%%%%%%%%%%%%%%%%%%%%%%%%%%%%%%%%%%%

\appendix
\section{Sampling Gaussian random fields}
A fast algorithm for generating Gaussian random fields is critical for obtaining an efficient multigrid MLMC estimator. When dealing with stationary covariance functions, one can utilize spectral generators, e.g. in \cite{RF4,RF2,Ravalec2000} that exploit the efficiency of the FFT algorithm to achieve fast sampling of these random fields on a uniform mesh. In this work, we use the Fast Fourier Transform moving average (FFT-MA) technique from \cite{Ravalec2000} to decompose the covariance matrix $C_{\Phi}(r)$ (recall \eqref{mat}). Although, this sampling method is similar to the Cholesky factorization technique, the key idea is to make the computational domain periodic. Thus, the resulting covariance operator is also periodic which can now be decomposed as a convolutional product. The samples of random fields are computed using cheaper vector-vector product compared to the expensive matrix-vector operation required when using Cholesky factorization. As a periodic covariance function sampled on a uniform grid results in a \emph{circulant} covariance matrix, these methods are sometimes also referred to as the circulant embedding technique.  In the following, we provide a brief description of FFT-MA method from \cite{Ravalec2000}.

Samples of correlated Gaussian random vectors $\mathbf{z}_\ell(\omega)$ can be obtained using a Cholesky decomposition of the covariance matrix $\mathbf{C}_{\Phi}^\ell$ on mesh ${D}_\ell$ as:
\begin{equation}
\mathbf{C}_{\Phi}^\ell  = \mathbf{L}_\ell\mathbf{L}_\ell^T\quad \text{ and use }\quad\mathbf{z}_\ell = \mathbf{L}_\ell\mathbf{y}_\ell,
\end{equation}
where $\mathbf{y}_{\ell}$ is a vector of i.i.d. samples from the standard normal distribution. The FFT-MA method is based on a decomposition of the covariance function $C_\Phi(r)$ as a convolutional product of some function $S_\Phi(r)$ and its transpose $S'_\Phi(r)$ ($S'_\Phi(r) = S_\Phi(-r)$). In a discrete setting, we can express this decomposition as
\begin{equation}\label{convol1}
\mathbf{c}_\ell = \mathbf{s}_\ell* \mathbf{s}'_\ell,
\end{equation}
where $\mathbf{c}_\ell,\mathbf{s}_\ell$ are vectors obtained by evaluating $C_\Phi(r)$ and $S_\Phi(r)$, respectively, at grid points of the mesh ${D}_\ell$. A correlated random vector $\mathbf{z}_\ell$ can now be synthesized by using the convolution product
\begin{equation}\label{convol2}
\mathbf{z}_\ell= \mathbf{s}_\ell*\mathbf{y}_\ell.
\end{equation}
The FFT-MA method performs the above computations in the frequency domain. First the vector $\mathbf{c}_\ell$ is transformed into a periodic signal, which is also real, positive and symmetric. More details on the practical aspects of this transformation, see \cite{RF4,RF2,QMC,mishra2016multi}. Note that the resulting vector $\mathbf{s}_\ell$ is also real, positive and symmetric and $\mathbf{s}_\ell = \mathbf{s}'_\ell$. As the convolution product in spatial domain is equivalent to the component-wise product in the frequency domain, we can rewrite \eqref{convol1} as
\begin{equation}\label{product1}
\mathcal{F}(\mathbf{c}_\ell) = \mathcal{F}(\mathbf{s}_\ell)\cdot \mathcal{F}(\mathbf{s}_\ell) \implies \mathcal{F}(\mathbf{s}_\ell) = \sqrt{\mathcal{F}(\mathbf{c}_\ell)},
\end{equation}
where $\mathcal{F}$ denotes the discrete FFT and $\cdot$ denotes component-wise multiplication. It is pointed out that the component-wise square-root operation does not pose any problems as the power spectrum $\mathcal{F}(\mathbf{c}_\ell)$ is real, positive. Next, we express the convolution product in \eqref{convol2} as a vector-vector product in the frequency domain as
\begin{equation}\label{product2}
\mathcal{F}(\mathbf{z}_\ell) = \mathcal{F}(\mathbf{s}_\ell*\mathbf{y}_\ell)=\mathcal{F}(\mathbf{s}_\ell)\cdot\mathcal{F}(\mathbf{y}_\ell).
\end{equation}
In the final step, an inverse fast Fourier transform is applied to obtain the samples for Gaussian random fields
\begin{equation}
\mathbf{z}_\ell = \mathcal{F}^{-1}(\mathcal{F}(\mathbf{s}_\ell)\cdot\mathcal{F}(\mathbf{y}_\ell)).
\end{equation} 
Due to the periodicity in the covariance vector $\mathbf{c}_\ell$, the resulting random field $\mathbf{z}_\ell$ is also periodic. Thus, we only retain the part of the vector that corresponds to the physical domain. We also remark that it takes two FFT evaluations to obtain one sample of $\mathbf{z}_\ell$ (ignoring the FFT operation in \eqref{product1} that is performed just once). Therefore, in terms of the number of floating point operations, the sampling cost is significantly smaller compared to one multigrid solve for the mesh sizes considered.
\section{Upscaling Gaussian random fields}
While estimating the correction term $\mathbb{E}[Q_{h_{\ell}} - Q_{h_{\ell-1}}]$ in the telescopic sum \eqref{ml}, the approximations $Q_{h_\ell}(\omega_i)$ and $Q_{h_{\ell-1}}(\omega_i)$ need to be positively correlated such that the variance $\mathbb{V}[Q_{h_{\ell}}- Q_{h_{\ell-1}}]$ is small. This is typically achieved by first sampling the fine grid permeability vector, ${k}(\mathbf{x}_\ell,\omega_i)$ to compute $Q_{h_\ell}(\omega_i)$ and using an upscaled version, ${k}(\mathbf{x}_{\ell-1},\omega_i)$ for $Q_{h_{\ell-1}}(\omega_i)$. While performing such upscaling of random fields, it is important to ensure that the telescopic sum \eqref{ml} is not voilated. In other words, the expectation of the random variable $Q_{h_\ell}$ when estimating $\mathbb{E}[Q_{h_{\ell}} - Q_{h_{\ell-1}}]$ and $\mathbb{E}[Q_{h_{\ell}} - Q_{h_{\ell+1}}]$ should be same, i.e.
\begin{equation}\label{fine_coarse}
\mathbb{E}[Q_{h_\ell}]^{(coarse)} = \mathbb{E}[Q_{h_\ell}]^{(fine)}, \quad\text{for} \quad \ell = \{0,1,...,L-1\}.
\end{equation}

\begin{figure}[H]
\begin{center}
\includegraphics[clip, trim=0.5cm 0cm 1cm 0cm,scale=0.26]{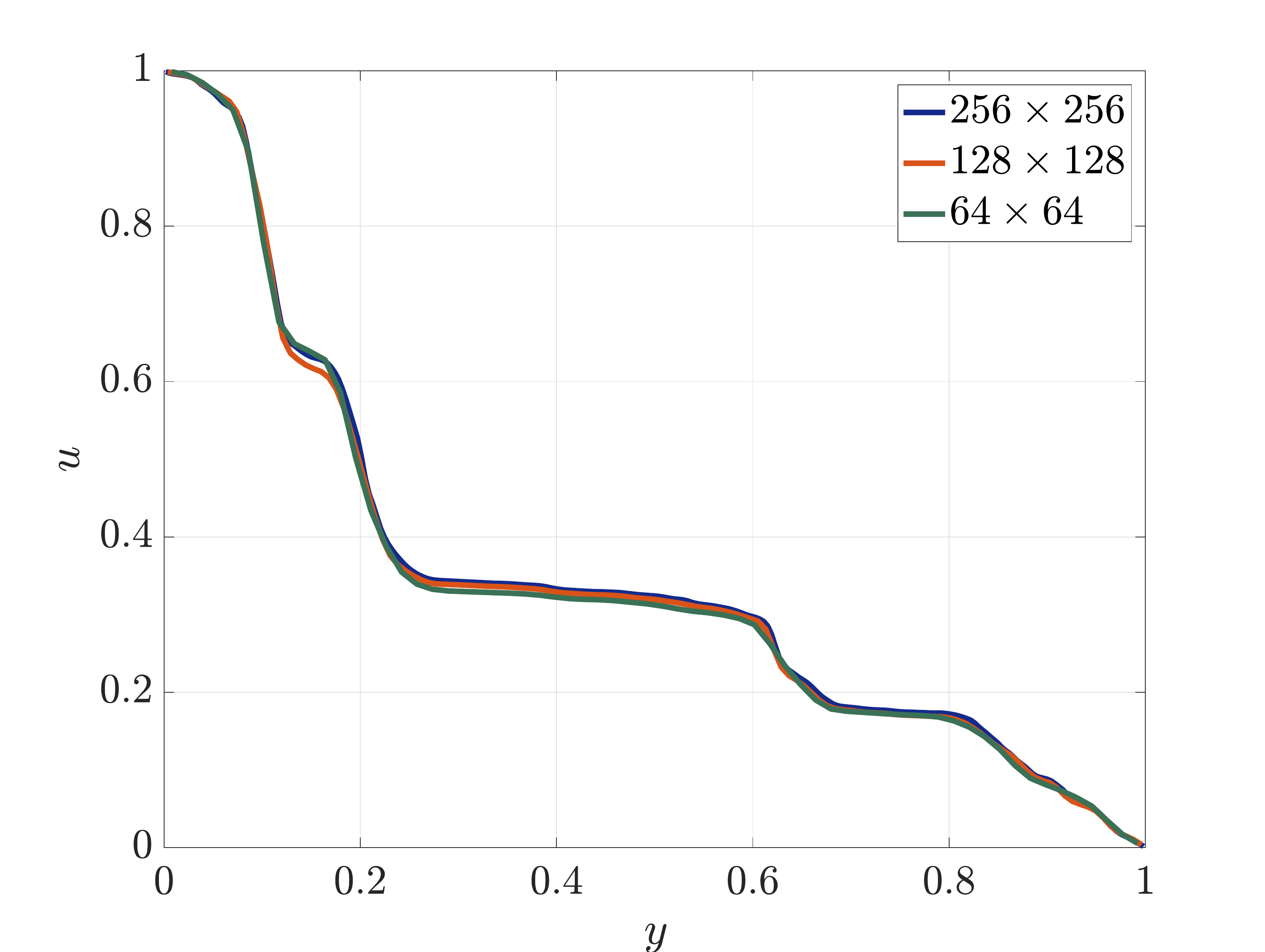}
\includegraphics[clip, trim=0.5cm 0cm 1cm 0cm,scale=0.26]{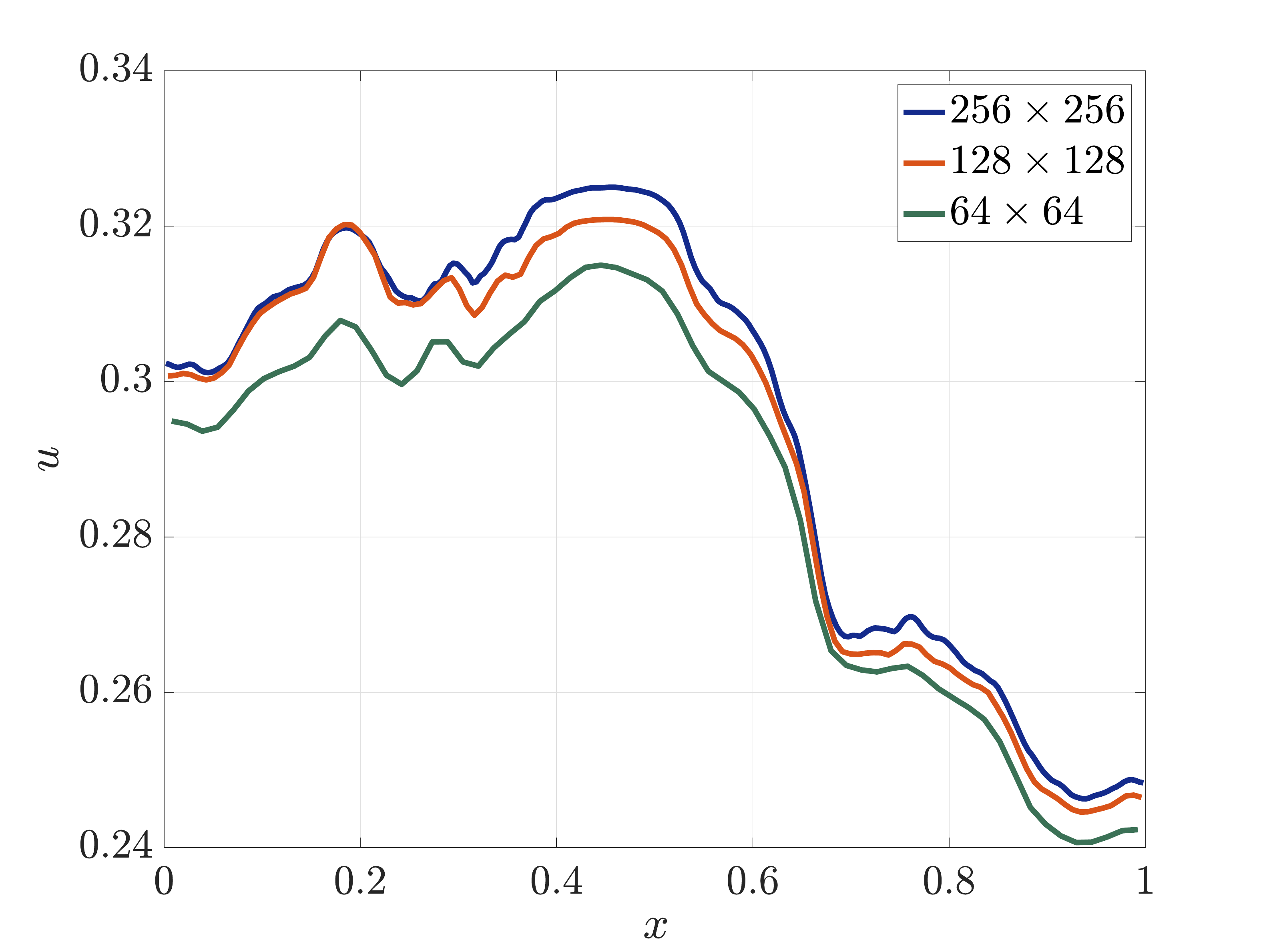}
\end{center}
\caption{Comparison of cross section of the pressure fields $u$ at $x =0.5$(Left) and $y = 0.5$(Right) obtained from a random permeability field (based on $\Phi_4$) on $256\times256$ grid, and from the corresponding upscaled permeabilities at $128\times128$ and $64\times64$ grids. }
\label{upscaled}
\end{figure}

Many of the upscaling algorithms based on homogenization techniques in the context of deterministic PDEs, such as \cite{moulton1,scott}, cannot be directly applied to obtain the upscaled permeability. This is because these homogenization procedures may result in a modified covariance structure on the coarser levels, violating \eqref{fine_coarse}. This issue can be avoided by using the \emph{covariance upscaling} \cite{mishra2016multi} that employs the spectral generator on two consecutive grids using the same normally distributed vector $\mathbf{y}_\ell$. Furthermore, in the case of the FFT-MA algorithm, the vector $\mathbf{y}_\ell$ is associated with respective grid points, coarser realizations of the fine grid Gaussian random field $\mathbf{z}_\ell$ ($\log k(\mathbf{x_\ell})$) can be obtained by using multi-dimensional averaging of the vector $\mathbf{y}^\ell$. For instance, in two dimensions for the cell-centred mesh,
\begin{equation}\label{averaging}
\mathbf{y}^{i,j}_{\ell-1} = \frac{1}{2}(\mathbf{y}^{2i-1,2j-1}_{\ell}+\mathbf{y}^{2i-1,2j}_{\ell}+\mathbf{y}^{2i,2j-1}_{\ell}+\mathbf{y}^{2i,2j}_{\ell}),
\end{equation}
where $(i,j)$ is the cell index for the mesh $\mathcal{D}_{\ell-1}$. The scaling by a factor 2 is needed to obtain a standard normal distribution for the averaged quantity $\mathbf{y}^{i,j}_{\ell-1}$. The coarse random field can now be simply assembled as
\begin{equation}\label{cov_upscale}
\mathbf{z}_{\ell-1} = \mathcal{F}^{-1}(\mathcal{F}(\mathbf{s}_{\ell-1})\cdot\mathcal{F}(\mathbf{y}_{\ell-1})).
\end{equation}
This process can be recursively applied to generate upscaled random fields on next coarser scales. As the averaging in \eqref{averaging} smooths out high frequencies, the upscaled version $\mathbf{z}_{\ell-1}$ will also be  smoother compared to ${z}_{\ell}$. In Figure \ref{upscaled}, we compare the cross sections of pressure fields obtained from a permeability field sampled (using parameter set $\Phi_4$) on a $256\times256$ grid and from the corresponding upscaled versions on $128\times128$ and $64\times64$ grids. The fine scale properties are well preserved on the coarser levels.

\section{Multigrid performance for anisotropic random fields} 
\begin{figure}[htb]
\begin{center}
\begin{tabular}{cc}
\includegraphics[clip, trim=2cm 1cm 0cm 0cm,scale=.32]{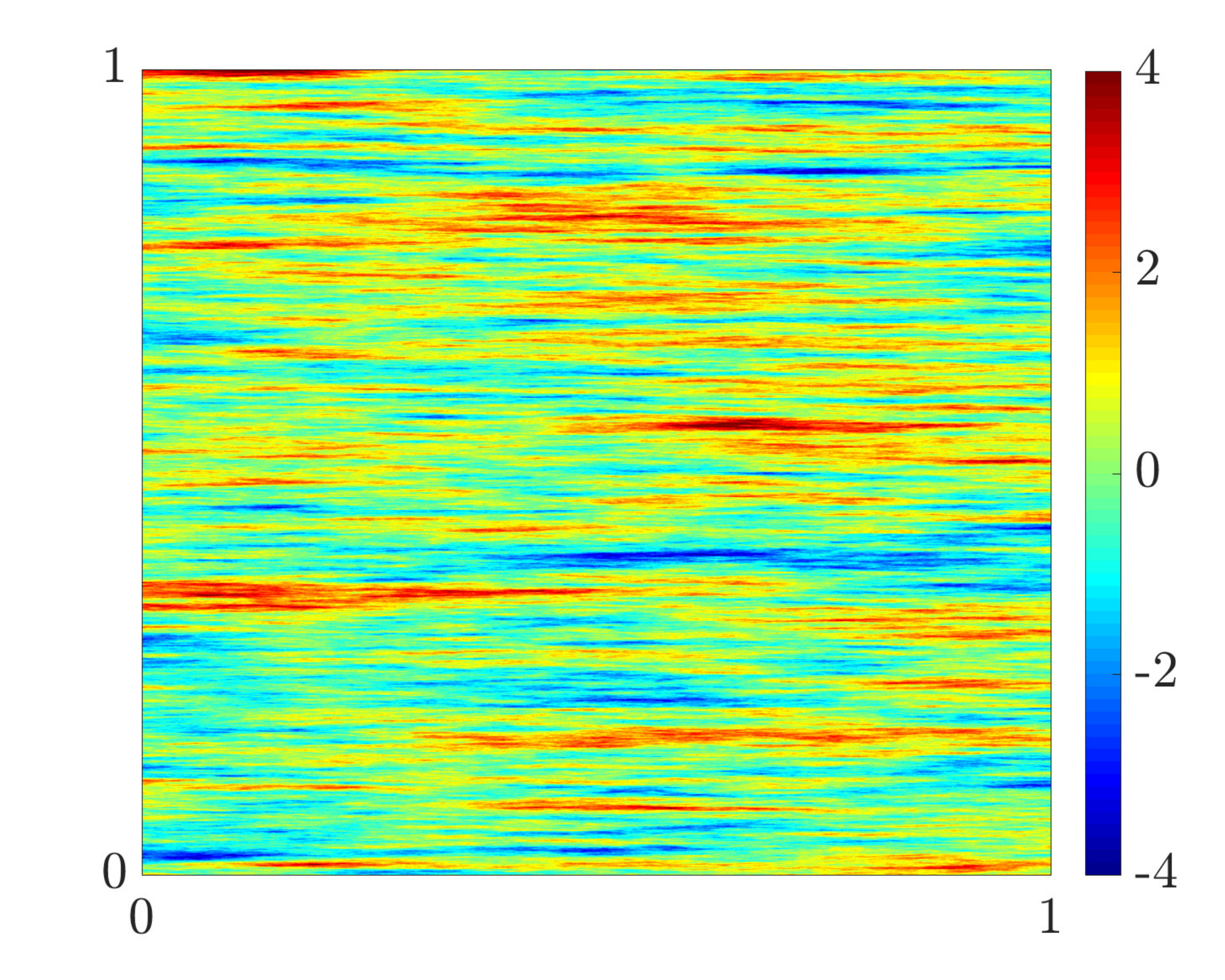} & 
\includegraphics[clip, trim=1cm 1cm 0cm 0cm,scale=.32]{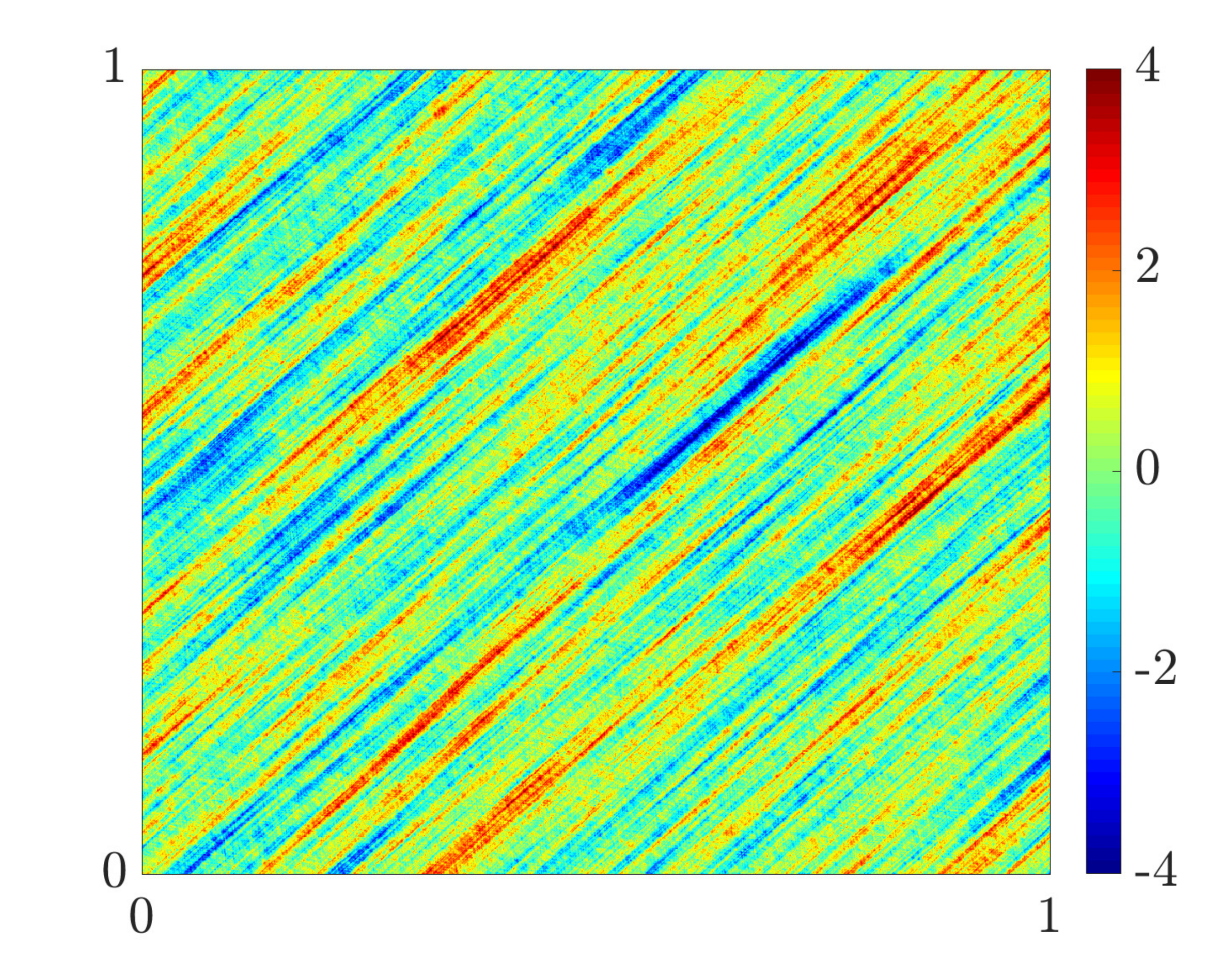} \\
(a) $\widetilde{\Phi}_1 = (0.5,0.3,0.01,1,0)$ & (b) $\widetilde{\Phi}_2 = (0.5,0.3,0.01,1,\pi/4)$ \\
\end{tabular}
\end{center}
\caption{Logarithm of the permeability field, $\log_{10} k$, generated using two reference parameter sets on  a unit square domain.}
\label{anisotropic}
\end{figure}

We study the performance of the cell-centered multigrid for anisotropic random fields. To generate these random fields, we use a modified version of the Mat\'ern covariance function, from \cite{kumarTransport},
\begin{equation}\label{aniso}
\begin{cases}
C_{\widetilde{\Phi}}(\mathbf{x_1},\mathbf{x_2}) = \sigma_c^2\dfrac{2^{1-\nu_c}}{\Gamma(\nu_c)}\left( 2\sqrt{\nu_c}\tilde{r}\right)^{\nu_c} K_{\nu_c}\left( 2\sqrt{\nu_c}\tilde{r}\right)\quad \mathbf{x_1},\mathbf{x_2}\in\mathcal{D},\\ 
\tilde{r} = \sqrt{\dfrac{(x_1'-x_2')^2}{\lambda^2_{cx}}+ \dfrac{(y_1'-y_2')^2}{\lambda^2_{cy}}},
\end{cases}
\end{equation}
where $C_{\widetilde{\Phi}}$ is a stationary covariance function depending on the parameter set ${\widetilde{\Phi}} = (\nu_c, \lambda_{cx}, \lambda_{cy}, \sigma_c^2,\theta)$ and $(x',y')$ corresponds to rotated coordinates by some angle $\theta$  in counterclockwise direction with respect to the horizontal axis, for e.g.
\begin{align}
x_1'  &= x_1\cos\theta - y_1\sin\theta,\nonumber\\
y_1'  &= x_1\sin\theta + y_1\cos\theta,\quad\text{with}\quad\mathbf{x_1} = (x_1,y_1)\nonumber.
\end{align}
The quantities $\lambda_{cx}$ and $\lambda_{cy}$ are correlation lengths along the x- and y-coordinates, respectively. The covariance function $C_{\widetilde{\Phi}}$ only differs from the isotropic covariance $C_\Phi$ defined in \eqref{mat} in terms of the distance function $\tilde{r}$. In Figure \ref{anisotropic}, we present realizations of the anisotropic random field (generated using circulant embedding) with two different $\widetilde{\Phi}$ values. As the two parameter sets only differ in terms of the rotation parameter $\theta$, they exhibit a similar magnitude of the jumps. Note that the random fields generated from $\widetilde{\Phi}_2$ are more challenging for the cell-centered multigrid as the long-range correlations are not aligned with the grid.

\begin{figure}[htb]
\begin{center}
\includegraphics[clip, trim=.9cm 0cm 2cm 0cm,scale=.29]{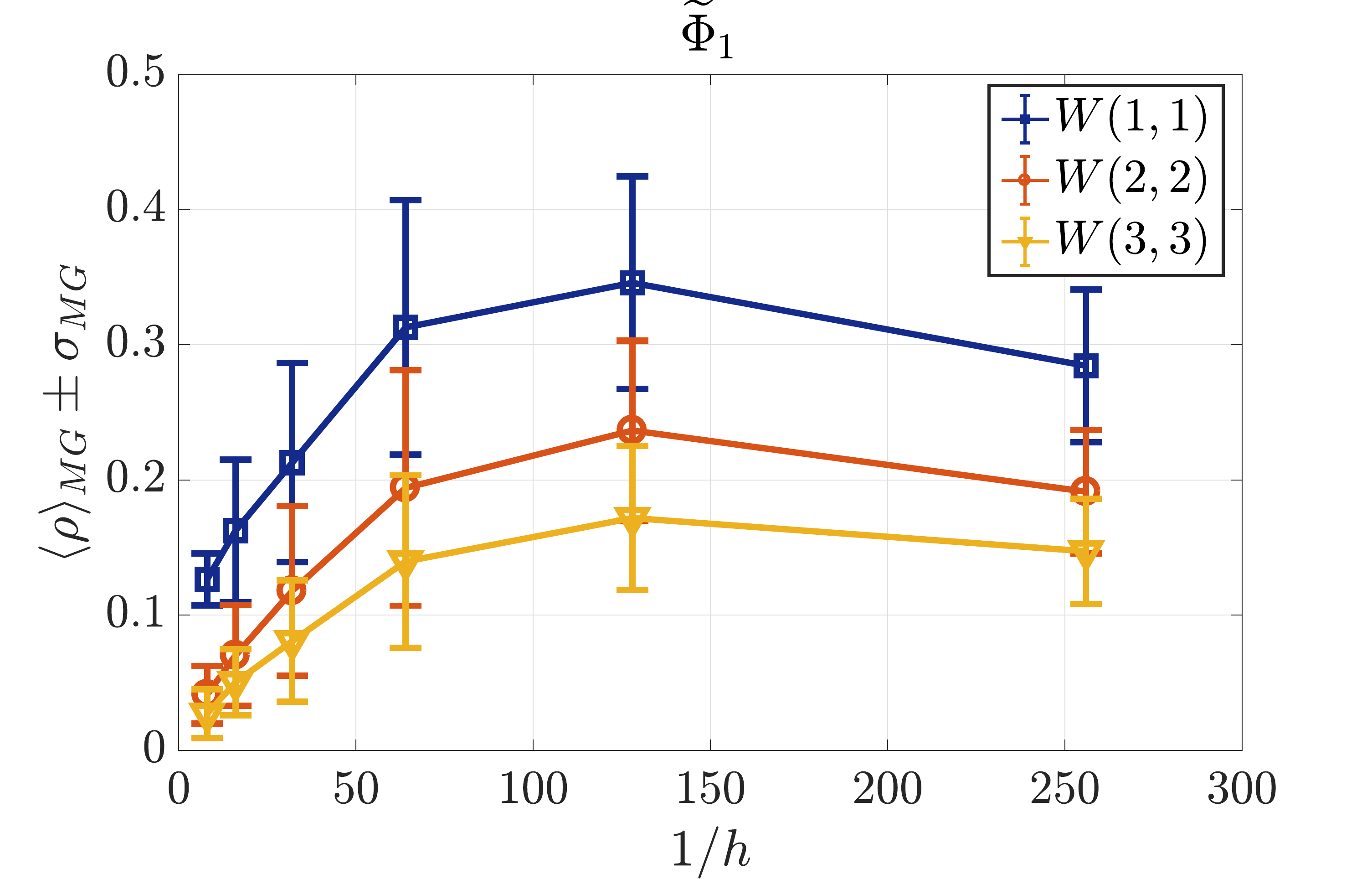}  
\includegraphics[clip, trim=.9cm 0cm 2cm 0cm,scale=.29]{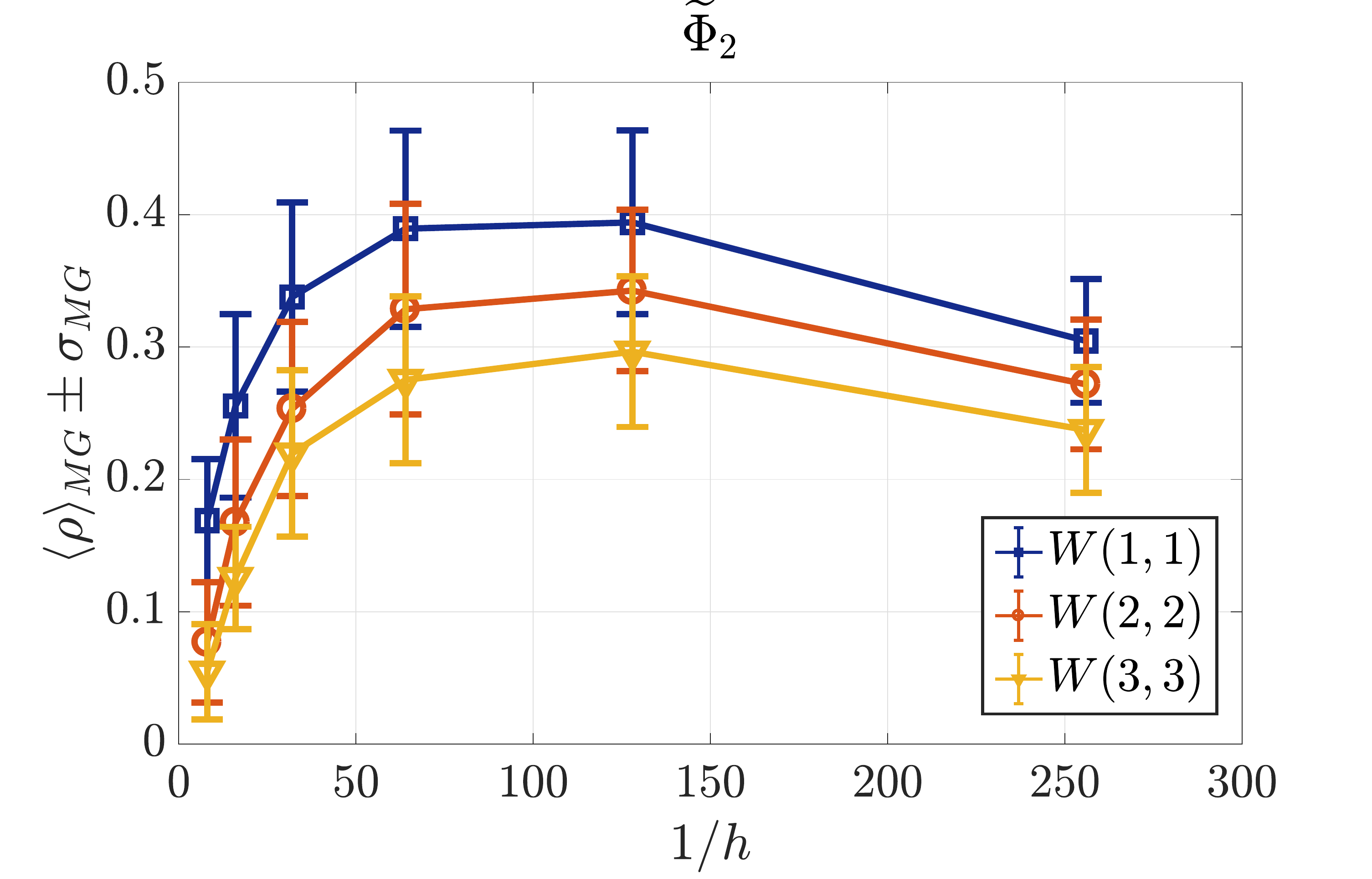}
\end{center}
\caption{Mean and the standard deviation of the asymptotic MG convergence for two reference parameter sets  $\widetilde{\Phi}_1$ (left) and $\widetilde{\Phi}_2$ (right) using $N_{MG} = 100$. }
\label{MG_conv_aniso}
\end{figure}

In Figure \ref{MG_conv_aniso}, we show the mean $\pm$ standard deviation of the asymptotic multigrid convergence for the two parameter sets, $\widetilde{\Phi}_1$ and $\widetilde{\Phi}_2$. Although, we see some deterioration compared to the isotropic case, the convergence rates improve with grid refinement. Here, using the $x$-line smoother for horizontal layering and the alternating line smoother for the non-aligned random field, keeping the other components same, will further improve the convergence.
%\section*{Acknowledgments}
%We would like to acknowledge the assistance of volunteers in putting
%together this example manuscript and supplement.
\newpage
\bibliographystyle{elsarticle-num}

%\bibliography{\myreferences}

\begin{thebibliography}{10}
\expandafter\ifx\csname url\endcsname\relax
  \def\url#1{\texttt{#1}}\fi
\expandafter\ifx\csname urlprefix\endcsname\relax\def\urlprefix{URL }\fi
\expandafter\ifx\csname href\endcsname\relax
  \def\href#1#2{#2} \def\path#1{#1}\fi

\bibitem{bolten}
M.~Bolten, H.~Rittich, \textit{Fourier Analysis of Periodic Stencils in
  Multigrid Methods}, SIAM Journal on Scientific Computing 40~(3) (2018)
  A1642--A1668.
\newblock \href {http://dx.doi.org/10.1137/16M1073959}
  {\path{doi:10.1137/16M1073959}}.

\bibitem{prob3}
J.~P. Delhomme, \textit{Spatial variability and uncertainty in groundwater flow
  parameters: A geostatistical approach}, Water Resources Research 15~(2)
  (1979) 269--280.
\newblock \href {http://dx.doi.org/10.1029/WR015i002p00269}
  {\path{doi:10.1029/WR015i002p00269}}.

\bibitem{WRCR:WRCR1821}
R.~A. Freeze, \textit{A stochastic-conceptual analysis of one-dimensional
  groundwater flow in nonuniform homogeneous media}, Water Resources Research
  11~(5) (1975) 725--741.
\newblock \href {http://dx.doi.org/10.1029/WR011i005p00725}
  {\path{doi:10.1029/WR011i005p00725}}.

\bibitem{WRCR:WRCR3746}
R.~J. Hoeksema, P.~K. Kitanidis, \textit{Analysis of the Spatial Structure of
  Properties of Selected Aquifers}, Water Resources Research 21~(4) (1985)
  563--572.
\newblock \href {http://dx.doi.org/10.1029/WR021i004p00563}
  {\path{doi:10.1029/WR021i004p00563}}.

\bibitem{babuska-collocation-2007}
I.~Babuska, F.~Nobile, R.~Tempone, {A stochastic collocation method for
  elliptic partial differential equations with random input data}, SIAM J.
  Numer. Anal. 45~(3) (2007) 1005--1034.
\newblock \href {http://dx.doi.org/10.1137/050645142}
  {\path{doi:10.1137/050645142}}.

\bibitem{Xiu:2010}
D.~Xiu, \textit{Numerical Methods for Stochastic Computations: A Spectral
  Method Approach}, Princeton University Press, Princeton, NJ, USA, 2010.

\bibitem{SEYNAEVE2007132}
B.~Seynaeve, E.~Rosseel, B.~Nicola, S.~Vandewalle, \textit{Fourier mode
  analysis of multigrid methods for partial differential equations with random
  coefficients }, Journal of Computational Physics 224~(1) (2007) 132 -- 149,
  special Issue Dedicated to Professor Piet Wesseling on the occasion of his
  retirement from Delft University of Technology.
\newblock \href {http://dx.doi.org/10.1016/j.jcp.2006.12.011}
  {\path{doi:10.1016/j.jcp.2006.12.011}}.

\bibitem{MLMC2}
K.~Cliffe, M.~B. Giles, R.~Scheichl, A.~L. Teckentrup, \textit{Multilevel Monte
  Carlo methods and applications to elliptic PDEs with random coefficients},
  Comput. Vis. Sci. 14 (2011) 3--15.
\newblock \href {http://dx.doi.org/10.1007/s00791-011-0160-x}
  {\path{doi:10.1007/s00791-011-0160-x}}.

\bibitem{MLMC1}
M.~B. Giles, \textit{Multilevel Monte Carlo path simulation}, Operations
  Research 256 (2008) 981--986.
\newblock \href {http://dx.doi.org/10.1287/opre.1070.0496}
  {\path{doi:10.1287/opre.1070.0496}}.

\bibitem{ANU:9672986}
M.~B. Giles, \textit{Multilevel Monte Carlo methods}, Acta Numerica 24 (2015)
  259--328.
\newblock \href {http://dx.doi.org/10.1017/S096249291500001X}
  {\path{doi:10.1017/S096249291500001X}}.

\bibitem{brandt2011multigrid}
A.~Brandt, O.~E. Livne, \textit{Multigrid techniques: 1984 guide with
  applications to fluid dynamics}, Vol.~67, SIAM, 2011.

\bibitem{MOLENAAR199625}
J.~Molenaar, \textit{A simple cell-centered multigrid method for 3D interface
  problems}, Computers and Mathematics with Applications 31~(9) (1996) 25 --
  33.
\newblock \href {http://dx.doi.org/10.1016/0898-1221(96)00039-9}
  {\path{doi:10.1016/0898-1221(96)00039-9}}.

\bibitem{Yav98}
I.~Yavneh, \textit{Coarse-Grid Correction for Nonelliptic and Singular
  Perturbation Problems}, SIAM J. Sci. Comput. 19~(5) (1998) 1682--1699.
\newblock \href {http://dx.doi.org/10.1137/S1064827596310998}
  {\path{doi:10.1137/S1064827596310998}}.

\bibitem{hackbusch85}
W.~Hackbusch, {Multigrid Methods and Applications}, Springer, 1985.

\bibitem{Bramble1996}
J.~H. Bramble, R.~E. Ewing, J.~E. Pasciak, J.~Shen, \textit{The analysis of
  multigrid algorithms for cell centered finite difference methods}, Advances
  in Computational Mathematics 5~(1) (1996) 15--29.
\newblock \href {http://dx.doi.org/10.1007/BF02124733}
  {\path{doi:10.1007/BF02124733}}.

\bibitem{Bramble1991}
J.~H. Bramble, J.~E. Pasciak, J.~Xu, \textit{The Analysis of Multigrid
  Algorithms with Nonnested Spaces or Noninherited Quadratic Forms},
  Mathematics of Computation 56~(193) (1991) 1--34.
\newblock \href {http://dx.doi.org/10.2307/2008527}
  {\path{doi:10.2307/2008527}}.

\bibitem{Brezina}
M.~Brezina, R.~Falgout, S.~MacLachlan, T.~Manteuffel, S.~McCormick, J.~Ruge,
  \textit{Adaptive Smoothed Aggregation ($\alpha$SA) Multigrid}, SIAM Review
  47~(2) (2005) 317--346.
\newblock \href {http://dx.doi.org/10.1137/S1064827502418598}
  {\path{doi:10.1137/S1064827502418598}}.

\bibitem{stuben2001introduction}
K.~St{\"u}ben, \textit{Appendix A: An introduction to algebraic multigrid}, in
  Multigrid, U. Trottenberg, C. W. Oosterlee, and A. Sch\"uller, Academic
  Press, San Diego, CA, 2001 (2001) 413--532.

\bibitem{Van2001}
P.~Van{\v{e}}k, M.~Brezina, J.~Mandel, \textit{Convergence of algebraic
  multigrid based on smoothed aggregation}, Numerische Mathematik 88~(3) (2001)
  559--579.
\newblock \href {http://dx.doi.org/10.1007/s211-001-8015-y}
  {\path{doi:10.1007/s211-001-8015-y}}.

\bibitem{Braess1995}
D.~Braess, \textit{Towards algebraic multigrid for elliptic problems of second
  order}, Computing 55~(4) (1995) 379--393.
\newblock \href {http://dx.doi.org/10.1007/BF02238488}
  {\path{doi:10.1007/BF02238488}}.

\bibitem{Van1996}
P.~Van{\v{e}}k, J.~Mandel, M.~Brezina, \textit{Algebraic multigrid by smoothed
  aggregation for second and fourth order elliptic problems}, Computing 56~(3)
  (1996) 179--196.
\newblock \href {http://dx.doi.org/10.1007/BF02238511}
  {\path{doi:10.1007/BF02238511}}.

\bibitem{Alcouffe:1981:MGM}
R.~E. Alcouffe, A.~Brandt, J.~E. {Dendy, Jr.}, J.~W. Painter, \textit{The
  multigrid methods for the diffusion equation with strongly discontinuous
  coefficients}, SIAM Journal on Scientific and Statistical Computing 2~(4)
  (1981) 430--454.
\newblock \href {http://dx.doi.org/10.1137/0902035}
  {\path{doi:10.1137/0902035}}.

\bibitem{dendy0}
J.~Dendy, \textit{Black box multigrid}, Journal of Computational Physics 48~(3)
  (1982) 366 -- 386.
\newblock \href {http://dx.doi.org/10.1016/0021-9991(82)90057-2}
  {\path{doi:10.1016/0021-9991(82)90057-2}}.

\bibitem{DENDY1983261}
J.~Dendy, Black box multigrid for nonsymmetric problems, Applied Mathematics
  and Computation 13~(3) (1983) 261 -- 283.
\newblock \href {http://dx.doi.org/10.1016/0096-3003(83)90016-4}
  {\path{doi:10.1016/0096-3003(83)90016-4}}.

\bibitem{MG2}
M.~{Khalil}, P.~{Wesseling}, \textit{Vertex-centered and cell-centered
  multigrid for interface problems.}, J. Comput. Phys. 98 (1992) 1--10.
\newblock \href {http://dx.doi.org/10.1016/0021-9991(92)90168-X}
  {\path{doi:10.1016/0021-9991(92)90168-X}}.

\bibitem{WESSELING198885}
P.~Wesseling, Cell-centered multigrid for interface problems, Journal of
  Computational Physics 79~(1) (1988) 85 -- 91.
\newblock \href {http://dx.doi.org/10.1016/0021-9991(88)90005-8}
  {\path{doi:10.1016/0021-9991(88)90005-8}}.

\bibitem{dendy2}
J.~E. Dendy, J.~D. Moulton, \textit{Black Box Multigrid with coarsening by a
  factor of three}, Numerical Linear Algebra with Applications 17~(2?3)
  577--598.
\newblock \href {http://dx.doi.org/10.1002/nla.705}
  {\path{doi:10.1002/nla.705}}.

\bibitem{MG1}
U.~{Trottenberg}, C.~W. {Oosterlee}, A.~{Schuller}, \textit{Multigrid},
  Elsevier Academic Press, San Diego, CA, 2000.

\bibitem{knapek}
S.~Knapek, \textit{Matrix-dependent multigrid homogenization for diffusion
  problems}, SIAM Journal on Scientific Computing 20~(2) (1998) 515--533.
\newblock \href {http://dx.doi.org/10.1137/S1064827596304848}
  {\path{doi:10.1137/S1064827596304848}}.

\bibitem{moulton1}
J.~Moulton, J.~E. Dendy, J.~M. Hyman, \textit{The Black Box Multigrid Numerical
  Homogenization Algorithm}, Journal of Computational Physics 142~(1) (1998) 80
  -- 108.
\newblock \href {http://dx.doi.org/10.1006/jcph.1998.5911}
  {\path{doi:10.1006/jcph.1998.5911}}.

\bibitem{scott}
S.~P. MacLachlan, J.~D. Moulton, Multilevel upscaling through variational
  coarsening, Water Resources Research 42~(2).
\newblock \href {http://dx.doi.org/10.1029/2005WR003940}
  {\path{doi:10.1029/2005WR003940}}.

\bibitem{MATERN}
M.~{Handcock}, J.~{Wallis}, \textit{An approach to statistical spatial-temporal
  modeling of meteorological fields (with discussion)}, Journal of the Americal
  Statistical Association 89 (1994) 368--390.
\newblock \href {http://dx.doi.org/10.2307/2290832}
  {\path{doi:10.2307/2290832}}.

\bibitem{Adler_book}
R.~Adler, The Geometry of Random Fields, Society for Industrial and Applied
  Mathematics, 2010.
\newblock \href {http://dx.doi.org/10.1137/1.9780898718980}
  {\path{doi:10.1137/1.9780898718980}}.

\bibitem{Nobile2015}
F.~Nobile, F.~Tesei, \textit{A multilevel Monte Carlo method with control
  variate for elliptic PDEs with log-normal coefficients}, Stochastic Partial
  Differential Equations: Analysis and Computations 3~(3) (2015) 398--444.
\newblock \href {http://dx.doi.org/10.1007/s40072-015-0055-9}
  {\path{doi:10.1007/s40072-015-0055-9}}.

\bibitem{doi:10.1137/110853054}
J.~Charrier, R.~Scheichl, A.~L. Teckentrup, \textit{Finite element error
  analysis of elliptic PDEs with random coefficients and its application to
  multilevel Monte Carlo methods}, SIAM Journal on Numerical Analysis 51~(1)
  (2013) 322--352.
\newblock \href {http://dx.doi.org/10.1137/110853054}
  {\path{doi:10.1137/110853054}}.

\bibitem{doi:10.1137/080717924}
J.~Galvis, M.~Sarkis, \textit{Approximating Infinity-Dimensional Stochastic
  Darcy's Equations without Uniform Ellipticity}, SIAM Journal on Numerical
  Analysis 47~(5) (2009) 3624--3651.
\newblock \href {http://dx.doi.org/10.1137/080717924}
  {\path{doi:10.1137/080717924}}.

\bibitem{RF4}
C.~{Dietrich}, G.~{Newsam}, \textit{Fast and exact simulation of stationary
  Gaussian processes through circulant embedding of the covariance matrix},
  SIAM J. Sci. Comput. 18 (1997) 1088--1107.
\newblock \href {http://dx.doi.org/10.1137/S1064827592240555}
  {\path{doi:10.1137/S1064827592240555}}.

\bibitem{RF2}
A.~{Wood}, G.~{Chan}, \textit{Simulation of stationary Gaussian processes in
  $[0, 1]^d$}, Journal of Computational and Graphical Statistics 3 (1994)
  409--432.
\newblock \href {http://dx.doi.org/10.1080/10618600.1994.10474655}
  {\path{doi:10.1080/10618600.1994.10474655}}.

\bibitem{RF1}
M.~{Ghanem}, P.~{Spanos}, \textit{Stochastic Finite Elements: A Spectral
  Approach.}, Springer, New York, 1991.
\newblock \href {http://dx.doi.org/10.1007/978-1-4612-3094-6}
  {\path{doi:10.1007/978-1-4612-3094-6}}.

\bibitem{RF3}
C.~{Schwab}, R.~{Todor}, \textit{Karhunen-Loeve approximation of random fields
  by generalized fast multipole methods}, J. Comput. Phys. 217(1) (2006)
  100--122.
\newblock \href {http://dx.doi.org/10.1016/j.jcp.2006.01.048}
  {\path{doi:10.1016/j.jcp.2006.01.048}}.

\bibitem{three_grid}
R.~Wienands, C.~W. Oosterlee, \textit{On Three-Grid Fourier Analysis for
  Multigrid}, {SIAM} J. Scientific Computing 23~(2) (2001) 651--671.
\newblock \href {http://dx.doi.org/10.1137/S106482750037367X}
  {\path{doi:10.1137/S106482750037367X}}.

\bibitem{mishra2012sparse}
S.~Mishra, C.~Schwab, \textit{Sparse tensor Multi-level Monte Carlo finite
  volume methods for hyperbolic conservation laws with random initial data},
  Mathematics of Computation 81~(280) (2012) 1979--2018.
\newblock \href {http://dx.doi.org/10.1090/S0025-5718-2012-02574-9}
  {\path{doi:10.1090/S0025-5718-2012-02574-9}}.

\bibitem{Mishra20123365}
S.~Mishra, C.~Schwab, J.~{\v{S}}ukys, \textit{Multi-level Monte Carlo finite
  volume methods for nonlinear systems of conservation laws in
  multi-dimensions}, Journal of Computational Physics 231~(8) (2012) 3365 --
  3388.
\newblock \href {http://dx.doi.org/10.1016/j.jcp.2012.01.011}
  {\path{doi:10.1016/j.jcp.2012.01.011}}.

\bibitem{mishra2016multi}
S.~Mishra, C.~Schwab, J.~{\v{S}}ukys, \textit{Multi-level Monte Carlo finite
  volume methods for uncertainty quantification of acoustic wave propagation in
  random heterogeneous layered medium}, Journal of Computational Physics 312
  (2016) 192--217.
\newblock \href {http://dx.doi.org/10.1016/j.jcp.2016.02.014}
  {\path{doi:10.1016/j.jcp.2016.02.014}}.

\bibitem{doi:10.1137/16M1083591}
D.~Drzisga, B.~Gmeiner, U.~R\"ude, R.~Scheichl, B.~Wohlmuth, \textit{Scheduling
  Massively Parallel Multigrid for Multilevel Monte Carlo Methods}, SIAM
  Journal on Scientific Computing 39~(5) (2017) S873--S897.
\newblock \href {http://dx.doi.org/10.1137/16M1083591}
  {\path{doi:10.1137/16M1083591}}.

\bibitem{10.1007/978-3-642-31464-3_25}
J.~{\v{S}}ukys, S.~, Mishra, C.~Schwab, \textit{Static Load Balancing for
  Multi-level Monte Carlo Finite Volume Solvers}, in: Parallel Processing and
  Applied Mathematics, Springer Berlin Heidelberg, Berlin, Heidelberg, 2012,
  pp. 245--254.
\newblock \href {http://dx.doi.org/10.1007/978-3-642-31464-3_25}
  {\path{doi:10.1007/978-3-642-31464-3_25}}.

\bibitem{Ravalec2000}
M.~L. Ravalec, B.~Noetinger, L.~Y. Hu, \textit{The FFT Moving Average (FFT-MA)
  Generator: An Efficient Numerical Method for Generating and Conditioning
  Gaussian Simulations}, Mathematical Geology 32~(6) (2000) 701--723.
\newblock \href {http://dx.doi.org/10.1023/A:1007542406333}
  {\path{doi:10.1023/A:1007542406333}}.

\bibitem{QMC}
I.~{Graham}, F.~{Kuo}, D.~{Nuyens}, R.~{Scheichl}, I.~{Sloan},
  \textit{Quasi-Monte Carlo methods for elliptic PDEs with random coefficients
  and applications}, J. Comput. Phys. 230(10) (2011) 3668--3694.
\newblock \href {http://dx.doi.org/10.1016/j.jcp.2011.01.023}
  {\path{doi:10.1016/j.jcp.2011.01.023}}.

\bibitem{kumarTransport}
P.~Kumar, P.~Luo, F.~J. Gaspar, C.~W. Oosterlee, \textit{A multigrid multilevel
  Monte Carlo method for transport in the Darcy-Stokes system}, Journal of
  Computational Physics 371 (2018) 382 -- 408.
\newblock \href {http://dx.doi.org/10.1016/j.jcp.2018.05.046}
  {\path{doi:10.1016/j.jcp.2018.05.046}}.

\end{thebibliography}
\end{document}